\documentclass[
  journal=psrm,
  manuscript=research-note,  
  year=2021,
  volume=2,
]{cup-journal}
\DeclareUnicodeCharacter{0301}{\'{e}}
\usepackage{booktabs,microtype,siunitx}
\usepackage{amsfonts,amsmath}
\usepackage{amssymb}
\usepackage{multirow}
\usepackage{tikz}
\usepackage{tcolorbox}
\usepackage{wrapfig}
\usepackage{amsthm}

\newtheorem{theorem}{Theorem}
\newtheorem{remark}{Remark}

\newtheorem{lemma}{Lemma}

\usepackage{amsmath}
\usepackage{amssymb}
\usepackage{mathtools}
\usepackage{amsthm}
\usepackage{amsfonts}       
\usepackage{amsmath}
\usepackage{amssymb}
\usepackage{amsfonts}
\usepackage{nicefrac}       
\usepackage{microtype}      
\usepackage{xcolor}         
\usepackage{graphicx}
\usepackage{enumitem}
\usepackage{multirow}
\usepackage{tikz}
\usepackage{tcolorbox}
\usepackage{wrapfig}
\usepackage{subfigure}
 \usepackage{relsize}
\usepackage{bbm}

\newtheorem{assum}{Assumption}[section]
\newtcbox{\alertinline}[1][red]
  {on line, arc = 0pt, outer arc = 0pt,
    colback = #1!20!white, colframe = #1!50!black,
    boxsep = 0pt, left = 1pt, right = 1pt, top = 2pt, bottom = 2pt,
    boxrule = 0pt, bottomrule = 1pt, toprule = 1pt}

\theoremstyle{plain}
 \usepackage{geometry}

\theoremstyle{definition}

\title{When can Regression-Adjusted Control Variates Help?\\{\Large Rare Events, Sobolev Embedding and Minimax Optimality}}

\author{Jose Blanchet}
\affiliation{alphabetical order}
\alsoaffiliation{Department of Management Science \& Engineering, Stanford University, CA, USA}
\alsoaffiliation{ICME, Stanford University, CA, USA}
\email{jose.blanchet@stanford.edu}

\author{Haoxuan Chen}
\affiliation{alphabetical order}
\alsoaffiliation{ICME, Stanford University, CA, USA}
\email{haoxuanc@stanford.edu}

\author{Yiping Lu}
\affiliation{alphabetical order}
\alsoaffiliation{ICME, Stanford University, CA, USA}
\email{yplu@stanford.edu}

\author{Lexing Ying}
\affiliation{alphabetical order}
\alsoaffiliation{ICME, Stanford University, CA, USA}
\alsoaffiliation{Department of Mathematics, Stanford University, CA, USA}
\email{lexing@stanford.edu}


\begin{document}

\begin{abstract}
This paper studies the use of a machine learning-based estimator as a control variate for mitigating the variance of Monte Carlo sampling. Specifically, we seek to uncover the key factors that influence the efficiency of control variates in reducing variance. We examine a prototype estimation problem that involves simulating the moments of a Sobolev function based on observations obtained from (random) quadrature nodes. Firstly, we establish an information-theoretic lower bound for the problem. We then study a specific quadrature rule that employs a nonparametric regression-adjusted control variate to reduce the variance of the Monte Carlo simulation. We demonstrate that this kind of quadrature rule can improve the Monte Carlo rate and achieve the minimax optimal rate under a sufficient smoothness assumption. Due to the Sobolev Embedding Theorem, the sufficient smoothness assumption eliminates the existence of rare and extreme events. Finally, we show that, in the presence of rare and extreme events, a truncated version of the Monte Carlo algorithm can achieve the minimax optimal rate while the control variate cannot improve the convergence rate.
\end{abstract}

\section{Introduction}

In this paper, we consider a nonparametric quadrature rule on (random) quadrature points based on regression-adjusted control variate \citep{asmussen2007stochastic,davidson1992regression,oates2016control,hickernell2005control}.  To construct the quadrature rule, we partition our available data into two halves. The first half is used to construct a nonparametric estimator, which is then utilized as a control variate to reduce the variance of the Monte Carlo algorithm implemented over the second half of our data. Traditional and well-known results \cite[Chapter 5.2]{asmussen2007stochastic} show that the optimal linear control variate can be obtained via Ordinary Least Squares regression. In this paper, we investigate a similar idea for constructing a quadrature rule \citep{oates2016control, assaraf1999zero,mira2013zero,oates2017control,oates2019convergence,south2018regularised,holzmuller2023convergence}, which uses a non-parametric machine learning-based estimator as a regression-adjusted control variate. We aim to answer the following two questions:

\begin{quotation}
{\emph{Is using optimal nonparametric machine learning algorithms to construct control variates an optimal way to improve Monte Carlo methods? What are the factors that determine the effectiveness of the control variate?}}
\end{quotation}

\begin{figure}
    \centering
    \includegraphics[width=5in]{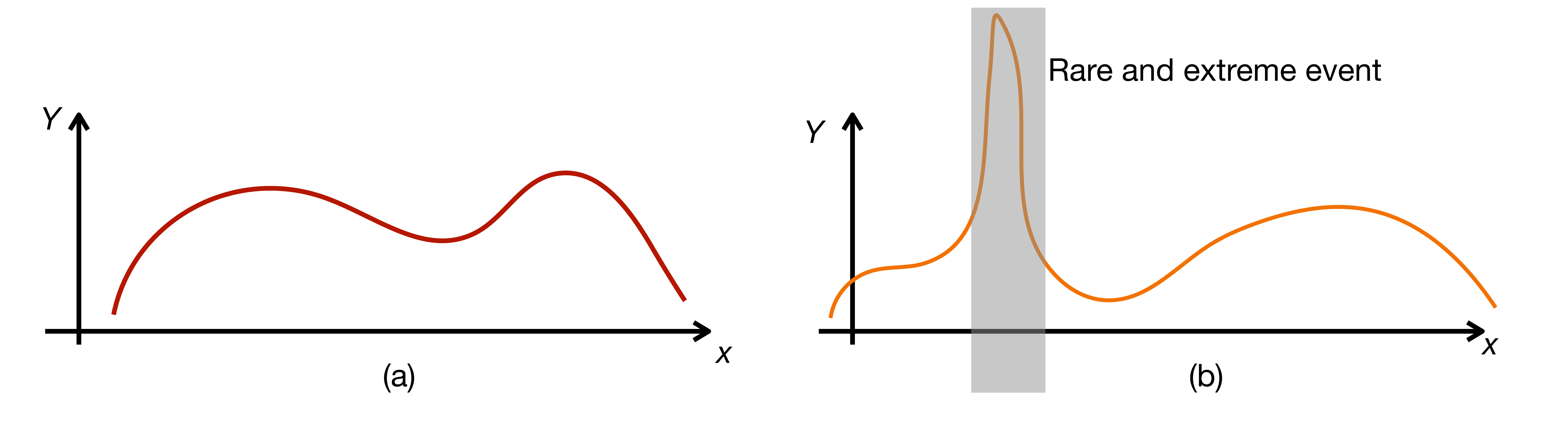}
    \caption{According to the Sobolev Embedding Theorem \citep{adams2003sobolev}, the Sobolev space $W^{s,p}$ can be embedded in $L^{p^\ast}$, where $\frac{1}{p^\ast}=\frac{1}{p}-\frac{s}{d}$. When $s$ is large enough, as shown in (a), the smoothness assumption can rule out the existence of rare and extreme events. When $s$ is not sufficiently large, specifically $s < \frac{2dq-dp}{2pq}$, there may exist a peak (\emph{a.k.a} rare and extreme event) that makes the Monte Carlo simulation hard. Under such circumstances, the function's $2q$-th moment is unbounded.}
    \label{fig:rare_event}
\end{figure}

To understand the two questions, we consider a basic but fundamental prototype problem of estimating moments of a Sobolev function from its values observed on (random) quadrature nodes, which has a wide range of applications in Bayesian inference, the study of complex systems, computational physics, and financial risk management \citep{asmussen2007stochastic}.  Specifically, we estimate the $q$-th moment $\int_{\Omega} f(x)^q dx$ of $f$ based on values $f(x_1),\cdots,f(x_n)$ observed on $n$ (random) quadrature nodes $x_1,\cdots,x_n\in \Omega$ for a function $f$ in the Sobolev space $W^{s,p}(\Omega)$, where $\Omega\subset\mathbb{R}^d$. The parameter $q$ here is introduced to characterize the rare events' extremeness for estimation. To verify the effectiveness of the non-parametric regression adjusted quadrature rule, we first study the statistical limit of the problem by providing a minimax information-theoretic lower bound of magnitude $n^{\max\{(\frac{1}{p}-\frac{s}{d})q-1,-\frac{s}{d}-\frac{1}{2}\}}$. 

\begin{figure}
    \centering

    \includegraphics[width=4.8in]{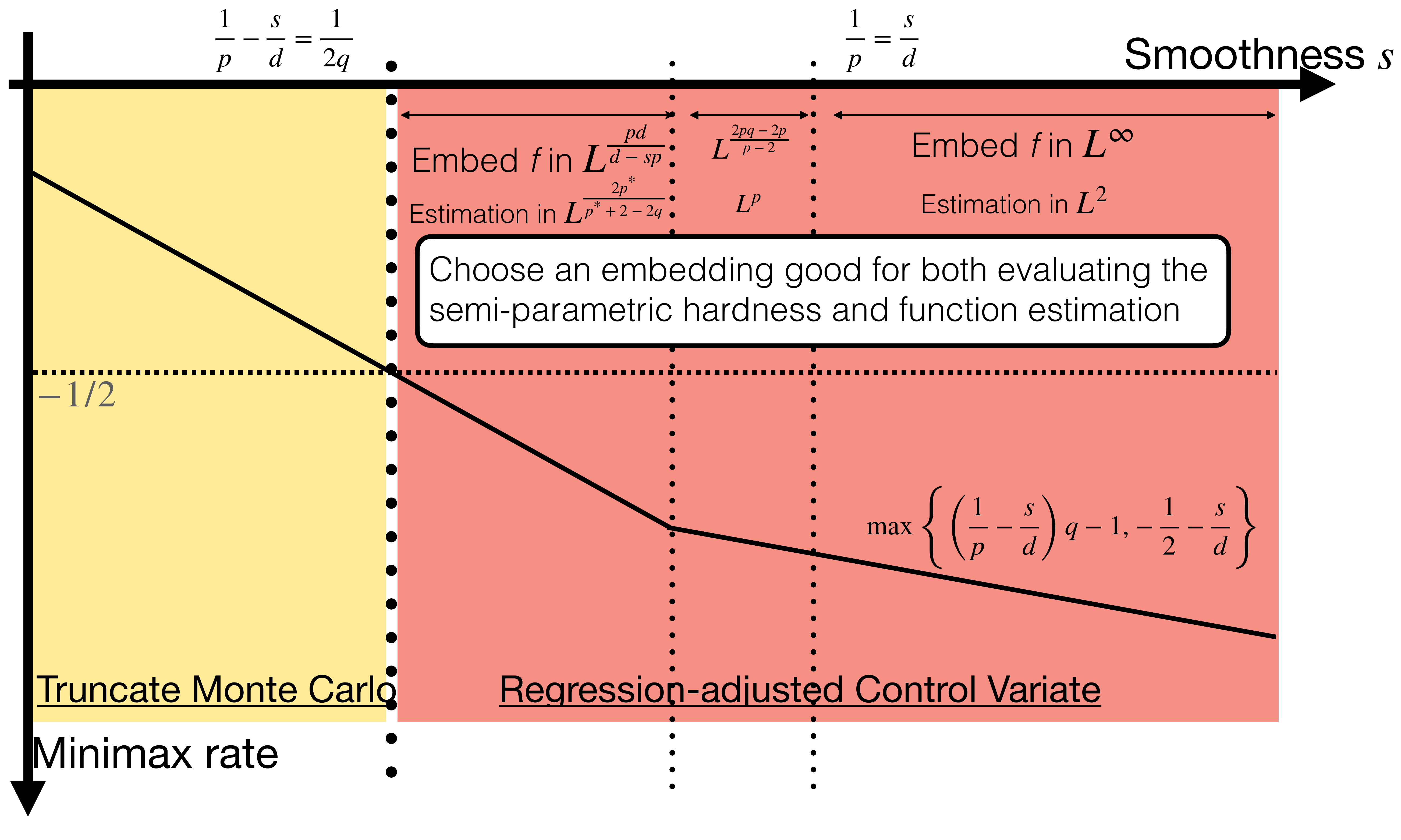}
    \caption{We summarize the minimax optimal rates and the corresponding optimal algorithms with respect to the function smoothness here. When the function is smooth enough, regression-adjusted control variates can improve the Monte Carlo rate. However, when there exist rare and extreme events that are hard to simulate, truncating the Monte Carlo estimate directly yields a minimax optimal algorithm. Above the transition point of algorithm selection is $s=\frac{d(2q-p)}{2pq}$, while the transition point of the optimal convergence rate is $s=\frac{d(2q-p)}{p(2q-2)}$.  To build the optimal convergence guarantee for any algorithm that utilizes a regression-adjusted control variate $\hat{f}$, we need to embed the square of the influence function $(qf^{q-1})^2$ in an appropriate space via the Sobolev Embedding Theorem and evaluate the estimation error $(f-\hat{f})^2$ under the dual norm of the norm associated with the chosen space, which allows us to achieve optimal semi-parametric efficiency. Our selections of the metrics in different regimes are shown in this figure.}
    \label{fig:threeregime}

\end{figure}

We also provide matching upper bounds for different levels of function smoothness. Under the sufficient smoothness assumption that $s>\frac{d(2q-p)}{2pq}$, we find that the non-parametric regression adjusted control variate $\hat{f}$ can improve the rate of classical Monte Carlo algorithm and help us attain a minimax optimal upper bound. In (\ref{eqn: thm 3.1 sketch decomp}) below, we bound variance $\int_{\Omega}(f^q-\hat{f}^q)^2$ of the Monte Carlo target by the sum of the semi-parametric influence part $\int_{\Omega}f^{2q-2}(f-\hat{f})^2$ and the propagated estimation error $\int_{\Omega}(f-\hat{f})^{2q}$. Although the optimal algorithm in this regime remains the same, we need to consider three different cases to derive an upper bound on the semi-parametric influence part, which is the main contribution of our proof. We propose a new proof technique that embeds the square of the influence function $(qf^{q-1})^2$ and estimation error $(f-\hat{f})^2$ in appropriate spaces via the Sobolev Embedding Theorem \citep{adams2003sobolev}. The two norms used for evaluating $(f^{q-1})^2$ and $(f-\hat{f})^2$ should be dual norms of each other. Also, we should select the norm for evaluating $(f-\hat{f})^2$ in a way that it's easy to estimate $f$ under the selected norm, which helps us control the error induced by $(f-\hat{f})^2$. A detailed explanation of how to select the proper norms in different cases via the Sobolev Embedding Theorem is exhibited in Figure \ref{fig:threeregime}. In the first regime when $s>\frac{d}{p}$, we can directly embed $f$ in $L^{\infty}(\Omega)$ and attain a final convergence rate of magnitude $n^{-\frac{s}{d}-\frac{1}{2}}$. For the second regime when $\frac{d(2q-p)}{p(2q-2)}<s<\frac{d}{p}$, the smoothness parameter $s$ is not large enough to ensure that $f \in L^{\infty}(\Omega)$. Thus, we evaluate the estimation error $(f-\hat{f})^2$ under the $L^{\frac{p}{2}}$ norm and embed the square of the influence function $(qf^{q-1})^2$ in the dual space of $L^{\frac{p}{2}}(\Omega)$. Here the validity of such embedding is ensured by the lower bound $\frac{d(2q-p)}{p(2q-2)}$ on $s$. Moreover, the semi-parametric influence part is still dominant in the second regime, so the final convergence rate is the same as that of the first case. In the third regime, when $\frac{d(2q-p)}{2pq}< s < \frac{d(2q-p)}{p(2q-2)}$, the semi-parametric influence no longer dominates and the final converge rate transits from $n^{-\frac{s}{d}-\frac{1}{2}}$ to $n^{q(\frac{1}{p}-\frac{s}{d})-1}$.



When the sufficient smoothness assumption breaks, \emph{i.e.} $s<\frac{d(2q-p)}{2pq}$, according to the Sobolev Embedding Theorem \citep{adams2003sobolev}, the Sobolev space $W^{s,p}$ is embedded in $L^{\frac{dp}{d-sp}}$ and $\frac{dp}{d-sp} < 2q$. This indicates that rare and extreme events might be present, and they are not even guaranteed to have bounded $L^{2q}$ norm, which makes the Monte Carlo estimate of the $q$-th moment have infinite variance. Under this scenario, we consider a truncated version of the Monte Carlo algorithm, which can be proved to attain the minimax optimal rate of magnitude $n^{q(\frac{1}{p}-\frac{s}{d})-1}$. In contrast, the usage of regression-adjusted control variates does not improve the convergence rate under this scenario. Our results reveal how the existence of rare events will change answers to the questions raised at the beginning of the section. 

We also use the estimation of a linear functional as an example to investigate the algorithm's adaptivity to the noise level. In this paper, we provide minimax lower bounds for estimating the integral of a fixed function with a general assumption on the noise level. Specifically, we consider all estimators that have access to observations $\{x_i,f(x_i)+\epsilon_i\}_{i=1}^n$ of some function $f$ that is $s$-H\"older smooth, where $x_i \overset{\text{i.i.d}}{\sim} \text{Uniform}([0,1]^d)$ and $\epsilon_i \overset{\text{i.i.d}}{\sim} n^{-\gamma}\mathcal{N}(0,1)$ for some $\gamma>0$. Based on the method of two fuzzy hypotheses, we present a lower bound of magnitude {$n^{\max\{-\frac{1}{2}-\gamma,-\frac{1}{2}-\frac{s}{d}\}}$}, which exhibits a smooth transition from the Monte Carlo rate to the Quasi-Monte Carlo rate. At the same time, our information-theoretic lower bound also matches the upper bound built for quadrature rules taking use of non-parametric regression-adjusted control variates.




\vspace{-0.05in}
\subsection{Related Work}
\vspace{-0.05in}
\paragraph{Regression-adjusted Control Variate}  Regression-adjusted control variates have shown both theoretical and empirical improvements in a broad range of applications, 
including the construction of confidence intervals \citep{angelopoulos2023prediction,romano2019conformalized}, randomized trace-estimation, \citep{meyer2021hutch++,sobczyk2022approximate,lin2017randomized}, dimension reduction \citep{sobczyk2022approximate}, causal inference \citep{liu2020regression}, estimation of the normalizing factor \citep{holzmuller2023convergence} and gradient estimation \citep{shi2022gradient,liu2017action}. It is also used as a technique used for proving the approximation bounds on two-layer neural networks in the Barron space \citep{siegel2022high}.

In connection to the related literature to our work, we mention \citep{oates2016control,oates2017control,oates2019convergence,holzmuller2023convergence}, which also study the use of nonparametric control variate estimator. However, the theoretical analysis in \citep{oates2016control,oates2017control} does not provide a specific convergence rate in the Reproducing Kernel Hilbert Space, which requires a high level of smoothness for the underlying function. In contrast to prior work, our research delves into the effectiveness of a non-parametric regression-adjusted control variate in boosting convergence rates across various degrees of smoothness assumptions and identifies the key factor that determines the effectiveness of these control variates. 

\paragraph{Quadrature Rule} There is a long literature on building quadrature rules in the Reproducing Kernel Hilbert Space, including Bayes–Hermite quadrature \citep{o1991bayes,kanagawa2016convergence,bach2017equivalence,karvonen2018fully,kanagawa2019convergence}, determinantal point processes \citep{belhadji2019kernel,belhadji2021analysis,bardenet2020monte,gautier2019two}, Nystr\"{o}m approximation \citep{hayakawa2021positively,hayakawa2023sampling}, kernel herding\citep{chen2012super,lacoste2015sequential,huszar2012optimally} and kernel thinning \citep{chen2018stein,dwivedi2021kernel,dwivedi2021generalized}. Nevertheless, the quadrature points chosen in these studies all have the ability to reconstruct the function's information, which results in a suboptimal rate for estimating the moments.

\paragraph{Functional Estimation} There are also lines of research that investigated the optimal rates of estimating both linear \citep{oates2019convergence, novak2006deterministic, traub1994information, novak2008tractability1,
novak2008tractability2, bakhvalov2015approximate, hinrichs2014curse, novak2016some, hinrichs2020power, hinrichs2022lower,krieg2020random, krieg2022recovery} and nonlinear \citep{birge1995estimation, donoho1990minimax,donoho1988one,donoho1991geometrizing2, donoho1991geometrizing3, robins2008higher,jiao2015minimax,krishnamurthy2014nonparametric, mathe1991random, heinrich2009randomized2, heinrich2009randomized3, han2020estimation,lepski1999estimation, heinrich2018complexity} functionals, such as integrals and the $L^q$ norm. However, as far as the authors know, previous works on this topic have assumed sufficient smoothness, which rules out the existence of rare and extreme events that are hard to simulate. Additionally, existing proof techniques are only applicable in scenarios where there is either no noise or a constant level of noise present. We have developed a novel and unified proof technique that leverages the method of two fuzzy hypotheses, which allows us to account for not only rare and extreme events but also different levels of noise. 

\subsection{Contribution}
\begin{itemize}
    \item We determine all the regimes when a quadrature rule utilizing a nonparametric estimator as a control variate to reduce the Monte Carlo estimate's variance can boost the convergence rate of estimating the moments of a Sobolev function. Under sufficient smoothness assumption, which rules out the existence of rare and extreme events due to the Sobolev Embedding Theorem, the regression-adjusted control variate improves the convergence rate and achieves the minimax optimal rate. The major technical difficulty in building the convergence guarantee in this regime is determining the right evaluation metric for function estimation. In our work, we bring a new proof technique to  select such a metric by embedding the influence function into an appropriate space via the Sobolev Embedding Theorem and evaluating the function estimation in the corresponding dual norm to achieve optimal semi-parametric efficiency. The selection of the metric is shown in Figure \ref{fig:threeregime}.
    \item Without the sufficient smoothness assumption, however, there may exist rare and extreme events that are hard to simulate. In this circumstance, we discover that a truncated version of the Monte Carlo method is minimax optimal, while regression-adjusted control variate can't improve the convergence rate.  As far as the authors know, our paper is the first work that considers this problem beyond the sufficient smoothness regime. 
    \item To study how the regression adjusted control variate adapts to the noise level, we examine the linear functionals, \emph{i.e.} the definite integral. We prove that this method is minimax optimal regardless of the level of noise present in the observed data. 
    
    
    
\end{itemize}

\subsection{Notations}
Let $||\cdot||$ be the standard Euclidean norm and $\Omega = [0,1]^d$ be the unit cube in $\mathbb{R}^d$ for any fixed $d \in \mathbb{N}$. Also, let $\mathbbm{1}=\mathbbm{1}\{\cdot\}$ denote the indicator function, \emph{i.e}, for any event $A$ we have $\mathbbm{1}\{A\}=1$ if $A$ is true and $\mathbbm{1}\{A\} = 0$ otherwise. For any region $R \subseteq \Omega$, we use $V(R) := \int_{\Omega}\mathbbm{1}\{x \in R\}dx$ to denote the volume of $R$. Let $C(\Omega)$ denote the space of all continuous functions $f:\Omega \rightarrow \mathbb{R}$ and $\lfloor \cdot\rfloor$ be the rounding function. For any $s>0$ and $f \in C(\Omega)$, we define the H\"older norm $||\cdot||_{C^{s}(\Omega)}$ by
\begin{equation}
\label{defn: holder norm}
||f||_{C^{s}(\Omega)} := \max_{|k| \leq \lfloor s \rfloor}||D^{k}f||_{L^{\infty}(\Omega)} + \max_{|k|=\lfloor s \rfloor}\sup_{x,y \in \Omega, x \neq y}\frac{|D^{k}f(x)-D^{k}f(y)|}{||x-y||^{s-\lfloor s \rfloor}}.    
\end{equation}
The corresponding H\"older space is defined as $C^s(\Omega):=\Big\{f\in C(\Omega): ||f||_{C^s(\Omega)} < \infty\Big\}$. When $s=0$, we have that the two norms $||\cdot||_{C^0(\Omega)}$ and $||\cdot||_{L^{\infty}(\Omega)}$ are equivalent and $C^{0}(\Omega) = L^{\infty}(\Omega)$. Let $\mathbb{N}_{0} := \mathbb{N} \cup \{0\}$ be the set of all non-negative integers. For any $s \in \mathbb{N}_{0}$ and $1 \leq p \leq \infty$, we define the Sobolev space $W^{s,p}(\Omega)$ by
\begin{equation}
\label{defn: sobolev space}
W^{s,p}(\Omega) := \Big\{f \in L^p(\Omega): D^{\alpha}f \in L^{p}(\Omega), \forall \ \alpha \in \mathbb{N}^{d}_{0} \text{ satisfying } |\alpha| \leq s\Big\}.    
\end{equation}
Let $(c)_{+}$ denote $\max\{c,0\}$ for any $c \in \mathbb{R}$. Fix any two non-negative sequences $\{a_n\}_{n=1}^{\infty}$ and $\{b_n\}_{n=1}^{\infty}$. We write $a_n \lesssim b_n$, or $a_n = O(b_n)$, to denote that $a_n \leq Cb_n$ for some constant $C$ independent of $n$. Similarly, we write $a_n \gtrsim b_n$, or $a_n = \omega(b_n)$, to denote that $a_n \geq cb_n$ for some constant $c$ independent of $n$. We use $a_n = \Theta(b_n)$ to denote that $a_n = O(b_n)$ and $a_n = \omega(b_n)$.  
\vspace{-0.1in}
\section{Information-Theoretic Lower Bound on Moment Estimation}
\label{sec 2: lower bound on q-moment estimation}
\paragraph{Problem Setup} To understand how the non-parametric regression-adjusted control variate improves the Monte Carlo estimator's convergence rate, we consider a prototype problem that estimates a function's $q$-th moment. For any fixed $q \in \mathbb{N}$ and $f \in W^{s,p}(\Omega)$, we want to estimate the $q$-th moment $I_f^q :=\int_{\Omega}f^q(x)dx$ with $n$ random quadrature points $\{x_i\}_{i=1}^{n} \subset \Omega$. On each quadrature point $x_i \ (i=1,\cdots,n)$, we can observe the function value $y_i :=f(x_i) $. 

In this section, we study the information-theoretic limit for the problem above via the method of two fuzzy hypotheses \citep{tsybakov2004introduction}. We have the following information-theoretic lower bound on the class $\mathcal{H}_{n}^{f,q}$ that contains all estimators $\hat{H}^{q}: \Omega^n \times \mathbb{R}^n \rightarrow \mathbb{R}$ of the $q$-th moment $I_{f}^{q}$.  
\begin{theorem}[Lower Bound on Estimating the Moment]
\label{thm: lower bound q-th moment}
When $p>2$ and $q<p<2q$, let $\mathcal{H}_n^{f}$ denote the class of all the estimators that use $n$ quadrature points $\{x_{i}\}_{i=1}^{n}$ and observed function values $\{y_i=f(x_i)\}_{i=1}^{n}$ to estimate the $q$-th moment of $f$, where $\{x_i\}_{i=1}^{n}$ are independently and identically sampled from the uniform distribution on $\Omega$. Then we have
\begin{equation}
\label{eq: lower bound q-th moment}
\inf_{\hat{H}^q \in \mathcal{H}_{n}^{f,q}}\sup_{f \in W^{s,p}(\Omega)}\mathlarger{\mathlarger{\mathbb{E}}}_{\substack{\{x_i\}_{i=1}^{n},\{y_i\}_{i=1}^{n}}}\Bigg[\left|\hat{H}^q\Big(\{x_i\}_{i=1}^{n},\{y_i\}_{i=1}^{n}\Big) - I_f^q\right|\Bigg] \gtrsim n^{\max\left\{-q\left(\frac{s}{d}-\frac{1}{p}\right)-1,-\frac{1}{2}-\frac{s}{d}\right\}}.  
\end{equation}
\end{theorem}

\paragraph{Proof Sketch} Here we give a sketch for our proof of Theorem \ref{thm: lower bound q-th moment}. Our proof is based on the method of two fuzzy hypotheses, which is a generalization of the traditional Le Cam's two-point method. In fact, each hypothesis in the generalized method is constructed via a prior distribution.  In order to attain a lower bound of magnitude $\Delta$ via the method of two fuzzy hypotheses, one needs to pick two prior distributions $\mu_0,\mu_1$ on the Sobolev space $W^{s,p}(\Omega)$ such that the following two conditions hold. Firstly, the estimators $I_{f}^q$ differ by $\Delta$ with constant probability under the two priors. Secondly, the TV distance between the two corresponding distributions $\mathbb{P}_0$ and $\mathbb{P}_1$ of data generated by $\mu_0$ and $\mu_1$ is of constant magnitude. In order to prove the two lower bounds given in (\ref{eq: lower bound q-th moment}), we pick two different pairs of prior distributions as follows:

Below we set $m=\Theta(n^{\frac{1}{d}})$ and divide the domain $\Omega$ into $m^d$ small cubes $\Omega_1,\Omega_2,\cdots,\Omega_{m^d}$, each of which has side length $m^{-1}$. For any $p \in (0,1)$, we use $v_{p},w_{p}$ to denote the discrete random variables satisfying $\mathbb{P}(v_p = 0) = \mathbb{P}(w_p = -1) = p$ and $\mathbb{P}(v_p = 1)= \mathbb{P}(w_p = 1)=1-p$.

(I) For the first lower bound in (\ref{eq: lower bound q-th moment}), we construct some bump function $g \in W^{s,p}(\Omega)$ satisfying $\text{supp}(g) \subseteq \Omega_1$ and $I_{g}^{q} = \int_{\Omega_1}g(x)dx = \Theta(m^{q(-s+\frac{d}{p})-d})$. Now let's take some sufficiently small constant $\epsilon \in (0,1)$ and pick $\mu_0,\mu_1$ to be discrete measures supported on the two finite sets $\Big\{v_{\frac{1+\epsilon}{2}}g\Big\}$ and  $\Big\{v_{\frac{1-\epsilon}{2}}g \Big\}$. On the one hand, the difference between the $q$-th moments under $\mu_0$ and $\mu_1$ can be lower bounded by $\Theta(n^{q(\frac{1}{p}-\frac{s}{d})-1})$ with constant probability. On the other hand, $KL(\mathbb{P}_0 || \mathbb{P}_1)$ can be upper bounded by the KL divergence between $v_{\frac{1+\epsilon}{2}}$ and $v_{\frac{1-\epsilon}{2}}$, which is of constant magnitude.

(II) For the second lower bound in (\ref{eq: lower bound q-th moment}), we set $M > 0$ to be some sufficiently large constant and $\kappa =\Theta(\frac{1}{\sqrt{n}})$. For any $1 \leq j \leq m^d$, we construct bump functions $f_j \in W^{s,p}(\Omega)$ satisfying $\text{supp}(f_j) \subseteq \Omega_j$ and $I_{f_j}^{k} = \int_{\Omega_j}f_j(x)dx = \Theta(m^{-ks-d})$ for any $1 \leq j \leq m^d$ and $1 \leq k \leq s$. Now let's pick $\mu_0,\mu_1$ to be discrete measures supported on the two finite sets $\Big\{M+\sum_{j=1}^{m^d}w_j^{(0)}f_j \Big\}$ and  $\Big\{M+\sum_{j=1}^{m^d}w_j^{(1)}f_j \Big\}$, where $\{w_{j}^{(0)}\}_{j=1}^{m^d}$ and $\{w_{j}^{(1)}\}_{j=1}^{m^d}$ are independent and identical copies of $w_{\frac{1+\kappa}{2}}$ and $w_{\frac{1-\kappa}{2}}$ respectively. On the one hand, applying Hoeffding's inequality yields that the $q$-th moments under $\mu_0$ and $\mu_1$ differ by $\Theta(n^{-\frac{s}{d}-\frac{1}{2}})$ with constant probability. On the other hand, note that $KL(\mathbb{P}_0 || \mathbb{P}_1)$ can be bounded by the KL divergence between two multivariate discrete distributions $(w_{j_1}^{(0)},\cdots,w_{j_n}^{(0)})$ and $(w_{j_1}^{(1)},\cdots,w_{j_n}^{(1)})$, where $\{w_{j_i}^{(0)}\}_{i=1}^{n}$ and $\{w_{j_i}^{(1)}\}_{i=1}^{n}$ are independent and identical copies of $w_{\frac{1+\kappa}{2}}$ and $w_{\frac{1-\kappa}{2}}$ respectively. Hence, $KL(\mathbb{P}_0 || \mathbb{P}_1)$ is of constant magnitude.

Combining the two cases above gives us the minimax lower bound in (\ref{eq: lower bound q-th moment}). We defer a complete proof of Theorem \ref{thm: lower bound q-th moment} to Appendix \ref{app: proof of lower bound in sec 2}.

\section{Minimax Optimal Estimators for Moment Estimation}
\label{sec 3: upper bound on q-moment estimation}
This section is devoted to constructing minimax optimal estimators of the $q$-th moment. We show that under the sufficient smoothness assumption, a regression-adjusted control variate is essential for building minimax optimal estimators. However, when the given function is not sufficiently smooth, we demonstrate that a truncated version of the Monte Carlo algorithm is minimax optimal, and control variates cannot give any improvement. 

\subsection{Sufficient Smoothness Regime: Non-parametric Regression-Adjusted Control Variate}

This subsection is devoted to building a minimax optimal estimator of the $q$-th moment under the assumption that $\frac{s}{d} > \frac{1}{p}-\frac{1}{2q}$, which guarantees that functions in the space $W^{s,p}$ are sufficiently smooth. From the Sobolev Embedding theorem, we know that the sufficient smoothness assumption implies $W^{s,p}(\Omega)\subset L^{p^\ast} (\Omega) \subset L^{2q}(\Omega)$, where $\frac{1}{p^\ast}=\frac{1}{p}-\frac{s}{d}$.
Given any function $f \in W^{s,p}(\Omega)$ along with $n$ uniformly sampled quadrature points $\{x_i\}_{i=1}^{n}$ and corresponding observations $\{y_i=f(x_i)\}_{i=1}^{n}$ of $f$, the key idea behind the construction of our estimator $\hat{H}^{q}_{C}$ is to build a nonparametric estimation $\hat{f}$ of $f$ based on a sub-dataset and use $\hat{f}$ as a control variate for Monte Carlo simulation. Consequently, it takes three steps to compute the numerical estimation of $I_f^q$ for any estimator $\hat{H}_{C}^{q}:\Omega^{n} \times \mathbb{R}^n \rightarrow \mathbb{R}$. The first step is to divide the observed data into two subsets $\mathcal{S}_1:=\{(x_i,y_i)\}_{i=1}^{\frac{n}{2}},\mathcal{S}_2:=\{(x_i,y_i)\}_{i=\frac{n}{2}+1}^{n}$ of equal size and use a machine learning algorithm to compute a nonparametric estimation $\hat{f}_{1:\frac{n}{2}}$ of $f$ based on $\mathcal{S}_1$. Without loss of generality, we may assume that the number of data points is even. Secondly, we treat $\hat{f}_{1:\frac{n}{2}}$ as a control variate and compute the $q$-th moment $I_{\hat{f}}^{q}$. Using the other dataset $\mathcal{S}_2$, we may obtain a Monte Carlo estimate of $I_{f}^q-I_{\hat{f}_{1:\frac{n}{2}}}^q$ as follows: 
$I_{f}^q-I_{\hat{f}_{1:\frac{n}{2}}}^q \approx \frac{2}{n}\sum_{i=\frac{n}{2}+1}^{n}\Big(y_i^q-\hat{f}_{1:\frac{n}{2}}^q(x_i)\Big). $ Finally, combining the estimation of the $q$-th moment $I_{\hat{f}_{1:\frac{n}{2}}}^{q}=\int_{\Omega}\hat{f}_{1:\frac{n}{2}}^{q}(x)dx$ with the estimation of $I_{f}^q-I_{\hat{f}_{1:\frac{n}{2}}}^q$ gives us the numerical estimation returned by $\hat{H}_{C}^{q}$:
\begin{equation}
\label{eqn: regression adjusted estimator}
\hat{H}_{C}^q\Big(\{x_i\}_{i=1}^{n},\{y_i\}_{i=1}^{n}\Big) := \int_{\Omega}\hat{f}_{1:\frac{n}{2}}^q(x)dx + \frac{2}{n}\sum_{i=\frac{n}{2}+1}^{n}\Big(y_i^q-\hat{f}_{1:\frac{n}{2}}^q(x_i)\Big). 
\end{equation}
We assume that our function estimation $\hat f$ is obtained from an $\frac{n}{2}$-oracle $K_{\frac{n}{2}}: \Omega^{\frac{n} {2}} \times \mathbb{R}^{\frac{n}{2}} \rightarrow W^{s,p}(\Omega)$ satisfying Assumption \ref{assumption: existence of oracle}.  For example, there are lines of research \citep{krieg2020random,krieg2022recovery,mathe1991random,heinrich2009randomized2,heinrich2009randomized3} considering how the moving least squares method \citep{wendland2001local,wendland2004scattered} can achieve the convergence rate in (\ref{oracle: assumption}). 
\begin{assum}[Optimal Function Estimator as an Oracle]
\label{assumption: existence of oracle}
Given any function $f\in W^{s,p}(\Omega)$ and $n \in \mathbb{N}$, let $\{x_i\}_{i=1}^{n}$ be $n$ data points sampled independently and identically from the uniform distribution on $\Omega$. Assume that there exists an oracle $K_n:\Omega^{n} \times \mathbb{R}^{n} \rightarrow W^{s,p}(\Omega)$ that estimates $f$ based on the $n$ points $\{x_{i}\}_{i=1}^{n}$ along with the $n$ observed function values $\{f(x_i)\}_{i=1}^{n}$ and satisfies the following bound for any $r$ satisfying $\frac{1}{r} \in (\frac{d-sp}{pd},1]$:
\begin{equation}
\label{oracle: assumption}
\left(\mathlarger{\mathbb{E}}_{\{x_{i}\}_{i=1}^{n}}\Big[||K_n(\{x_i\}_{i=1}^{n},\{f(x_{i})\}_{i=1}^{n})-f||_{L^r(\Omega)}^{r}\Big]\right)^{\frac{1}{r}} \lesssim n^{-\frac{s}{d}+(\frac{1}{p}-\frac{1}{r})_{+}}.
\end{equation}
\end{assum}

Based on the oracle above, we can obtain the following upper bound that matches the information-theoretic lower bound in Theorem \ref{thm: lower bound q-th moment}.

\begin{theorem}[Upper Bound on Moment Estimation with Sufficient Smoothness]
\label{thm: upper bound control variate}
Assume that $p > 2$, $q < p < 2q$ and $s>\frac{2dq-dp}{2pq}$. Let $\{x_{i}\}_{i=1}^{n}$ be $n$ quadrature points independently and identically sampled from the uniform distribution on $\Omega$ and $\{y_i:= f(x_i) \}_{i=1}^n$ be the corresponding $n$ observations of $f \in W^{s,p}(\Omega)$. Then the estimator $\hat{H}_{C}^q$ constructed in (\ref{eqn: regression adjusted estimator}) above satisfies
\begin{equation}
\label{eq: upper bound (control variate)}
\mathlarger{\mathlarger{\mathbb{E}}}_{\substack{\{x_i\}_{i=1}^{n},\{y_i\}_{i=1}^{n}}}\Bigg[\left|\hat{H}_{C}^q\Big(\{x_i\}_{i=1}^{n},\{y_i\}_{i=1}^{n}\Big) - I_f^q\right|\Bigg] \lesssim n^{\max\{-q(\frac{s}{d}-\frac{1}{p})-1,-\frac{s}{d}-\frac{1}{2}\}},    
\end{equation}
\end{theorem}

\paragraph{Proof Sketch} Given a non-parametric estimator $\hat f$ of the function $f$, we may bound the variance of the Monte Carlo process by $(f^q -\hat{f}^q)^2$ and further upper bound it by the sum of the following two terms:
\begin{equation}
\label{eqn: thm 3.1 sketch decomp}
\begin{aligned}
|f^q-\hat f^q|^2\lesssim \underbrace{|f^{q-1}(f-\hat f)|^2}_{\text{semi-parametric influnce}} + \underbrace{|(f-\hat f)^{q}|^2}_{\text{estimation error  propagation }}.
\end{aligned}
\end{equation}
The first term above represents the semi-parametric influence part of the problem, as $qf^{q-1}$ is the influence function for the estimation of the $q$-th moment $f^q$. The second term characterizes how function estimation affects functional estimation. If we consider the special case of estimating the mean instead of a general $q$-th moment, \emph{i.e,} $q=1$, the semi-parametric influence term will disappear. Consequently, the convergence rate won't transit from $n^{-\frac{1}{2}-\frac{s}{d}}$ to $n^{-q(\frac{s}{d}-\frac{1}{p})-1}$ in the special case.

Although the algorithm remains unchanged in the sufficient smooth regime, we need to consider three separate cases to obtain an upper bound on the integral of the semi-parametric influence term $|f^{q-1}(f-\hat{f})|^2$ in (\ref{eqn: thm 3.1 sketch decomp}). An illustration of the three cases is given in Figure \ref{fig:threeregime}.

From H\"older's inequality, we know that  $\int_{\Omega} {f^{2q-2}(x)(f(x)-\hat{f}(x))^2}dx$ can be upper bounded by $||f^{2q-2}||_{L^{r'}(\Omega)}||(f-\hat{f})^2||_{L^{r_\ast}(\Omega)}$, where $||\cdot||_{L^{r'}(\Omega)}$ and $||||_{L^{r_\ast}(\Omega)}$ are dual norms. Therefore, the main difficulty here is to embed the function $f$ in different spaces via the Sobolev Embedding Theorem under different assumptions on the smoothness parameter $s$. When the function is smooth enough, \emph{i.e.} $s>\frac{d}{p}$, we embed the function $f$ in $L^{\infty}(\Omega)$ and evaluate the estimation error $f-\hat{f}$ under the $L^2$ norm. Then our assumption on the oracle (\ref{oracle: assumption}) gives us an upper bound of magnitude $n^{-\frac{2s}{d}}$ on $||f-\hat{f}||_{L^2(\Omega)}^2$, which helps us further upper bound the semi-parametric influence part $\int_{\Omega} {f^{2q-2}(x)(f(x)-\hat{f}(x))^2}dx$ by $n^{-\frac{2s}{d}}$ up to constants. When $\frac{d(2q-p)}{p(2q-2)}<s<\frac{d}{p}$, we embed the function $f$ in $L^{\frac{2pq-2p}{{p}-2}}(\Omega) \subseteq L^{\frac{pd}{d-sp}}(\Omega)$ and evaluate the estimation error $f-\hat{f}$ under the $L^{p}$ norm. Applying our assumption on the oracle (\ref{oracle: assumption}) again implies that the semi-parametric influence part $\int_{\Omega} {f^{2q-2}(x)(f(x)-\hat{f}(x))^2}dx$ can be upper bounded by $n^{-\frac{2s}{d}}$ up to constants. When $\frac{d(2q-p)}{2pq}<s<\frac{d(2q-p)}{p(2q-2)}$, we embed the function $f$ in $L^{p^\ast}$ and evaluate the error of the oracle in $L^{\frac{2p^\ast}{p^\ast+2-2q}}$, where $\frac{1}{p^\ast} =\frac{1}{p}-\frac{s}{d}$. Similarly,  we can use (\ref{oracle: assumption}) to upper bound the semi-parametric influence part $\int_{x\in\Omega} {f^{2q-2}(x)(f(x)-\hat{f}(x))^2}dx$ by $n^{2q(\frac{1}{p}-\frac{s}{d})-1}$.



The upper bound on the propagated estimation error $\int_{x\in\Omega} {(f(x)-\hat{f}(x))^{2q}}dx$ in (\ref{eqn: thm 3.1 sketch decomp}) can be derived by evaluating the error of the oracle under the $L^{2q}$ norm. \emph{i.e,} by picking $r=2q$ in (\ref{oracle: assumption}) above, which yields an upper bound of magnitude $n^{2q(\frac{1}{p}-\frac{s}{d})-1}$. 

The obtained upper bounds on the semi-parametric influence part and the propagated estimation error above provide us with a clear view of the upper bound on the variance of $f^q-\hat f^q$, which is the random variable we aim to simulate via Monte-Carlo in the second stage. Using the standard Monte-Carlo algorithm to simulate the expectation of $f^q-\hat f^q$ then gives us an extra $n^{-\frac{1}{2}}$ factor for the convergence rate, which helps us attain the final upper bounds given in (\ref{eq: upper bound (control variate)}).  A complete proof of Theorem \ref{thm: upper bound control variate} is given in Appendix \ref{app: proof of upper bound in sec 3.1}.

\subsection{Beyond the Sufficient Smoothness Regime: Truncated Monte Carlo}

In this subsection, we study the case when the sufficient smoothness assumption breaks, \emph{i.e.} $\frac{s}{d} < \frac{1}{p}-\frac{1}{2q}$. According to the Sobolev Embedding theorem, we have that $W_s^{p}$ is embedded in $L^{\frac{dp}{d-sp}}$. Since $\frac{1}{p}-\frac{s}{d}>\frac{1}{2q}$ implies $\frac{dp}{d-sp} <2q$, the underlying function $f$ is not guaranteed to have bounded $L^{2q}$ norm, which indicates the existence of rare and extreme events. Consequently, the Monte Carlo estimate of $f$'s $q$-th moment must have infinite variance, which makes it hard to simulate. Here we present a truncated version of the Monte Carlo algorithm that can achieve the minimax optimal convergence rate. For any fixed parameter $M > 0$, our estimator is designed as follows:
\begin{equation}
\label{eqn: truncated monte carlo}
\hat{H}_{M}^{q}\Big(\{x_i\}_{i=1}^{n},\{y_i\}_{i=1}^{n}\Big) := \frac{1}{n}\sum_{i=1}^{n}\max\Big\{\min\{y_i,M\},-M\Big\}^q.
\end{equation}
In Theorem \ref{thm: upper bound truncate MC}, we provide the convergence rate of the estimator (\ref{eqn: truncated monte carlo}) by choosing the truncation parameter $M$ in an optimal way. 
\begin{theorem}[Upper Bound on Moment Estimation without Sufficient Smoothness]
\label{thm: upper bound truncate MC}
Assuming that $p > 2$, $q < p < 2q$ and $s<\frac{2dq-dp}{2pq}$, we pick $M=\Theta(n^{\frac{1}{p}-\frac{s}{d}})$. Let $\{x_{i}\}_{i=1}^{n}$ be $n$ quadrature points independently and identically sampled from the uniform distribution on $\Omega$ and $\{y_i:= f(x_i) \}_{i=1}^n$ be the corresponding $n$ observations of $f \in W^{s,p}(\Omega)$. Then we have that the estimator $\hat{H}_{M}^q$ constructed in (\ref{eqn: truncated monte carlo}) above satisfies
\begin{equation}
\label{eq: upper bound (truncate MC)}
\mathlarger{\mathlarger{\mathbb{E}}}_{\substack{\{x_i\}_{i=1}^{n},\{y_i\}_{i=1}^{n}}}\Bigg[\left|\hat{H}_{M}^q\Big(\{x_i\}_{i=1}^{n},\{y_i\}_{i=1}^{n}\Big) - I_f^q\right|\Bigg] \lesssim n^{-q(\frac{s}{d}-\frac{1}{p})-1}.    
\end{equation}
\end{theorem}

\paragraph{Proof Sketch}  

The error can be decomposed into bias and variance parts. The bias part is caused by the truncation in our algorithm, which is controlled by the parameter $M$ and can be bounded by $\int_{\{x:|f(x)|>M\}} |f|^qdx$. According to the Sobolev Embedding Theorem, $W^{s,p}(\Omega)$ can be embedded in the space $L^{p^\ast}$, where $\frac{1}{p^\ast}=\frac{1}{p}-\frac{s}{d}$. As $|f(x)| > M$ implies $|f(x)|^q<M^{q-p^\ast}|f(x)|^{p^\ast}$, the bias can be upper bounded by $M^{q-p^\ast}$. Similarly, the variance is controlled by $M$ and can be upper bounded by $M^{q-\frac{p^\ast}{2}}$. Combining the bias and variance bound, we can bound the final error as $M^{q-p^\ast}+\frac{M^{q-\frac{p^\ast}{2}}}{\sqrt{n}}$. By selecting $M=\Theta(n^{\frac{1}{p^\ast}})=\Theta(n^{\frac{1}{p}-\frac{s}{d}})$, we obtain the final convergence rate $ n^{-q(\frac{s}{d}-\frac{1}{p})-1}$. A complete proof of Theorem \ref{thm: upper bound truncate MC} is given in Appendix \ref{app: proof of upper bound in sec 3.2}.

\begin{remark}
\citep{heinrich2018complexity} has shown that the convergence rate of the optimal non-parametric regression-based estimation is $n^{-\frac{s}{d}+\frac{1}{p}-\frac{1}{q}}$, which is slower than the convergence rate of the truncated Monte Carlo estimator that we show above. 
\end{remark}

\section{Adapting to the Noise Level: a Case Study for Linear Functional}
\label{sec 4: case study}
In this section, we study how the regression-adjusted control variate adapts to different noise levels. Here we consider the linear functional, \emph{i.e.} estimating a function's definite integral via low-noise observations at random points.

\paragraph{Problem Setup} We consider estimating $I_f=\int_{\Omega}f(x)dx$, the integral of $f$ over $\Omega$, for a fixed function $f \in C^{s}(\Omega)$ with uniformly sampled quadrature points $\{x_i\}_{i=1}^{n} \subset \Omega$. On each quadrature point $x_i \ (i=1,\cdots,n)$, we have a noisy observation $y_i :=f(x_i) + \epsilon_i$. Here the $\epsilon_i$'s are independently and identically distributed Gaussian noises sampled from $\mathcal{N}(0,n^{-2\gamma})$, where $\gamma \in [0,\infty]$.

\subsection{Information-Theoretic Lower Bound on Mean Estimation}
In this subsection, we present a minimax lower bound (Theorem \ref{thm: lower bound for integral}) for all estimators $\hat{H}:\Omega^{n} \times \mathbb{R}^n\rightarrow\mathbb{R}$ of the integral $I_f$ of a function $f\in C^{s}(\Omega)$ when one can only access noisy observations.

\begin{theorem}[Lower Bound for Integral Estimation]
\label{thm: lower bound for integral}
Let $\mathcal{H}_n^{f}$ denote the class of all the estimators that use $n$ quadrature points $\{x_{i}\}_{i=1}^{n}$ and noisy observations $\{y_i = f(x_i)+\epsilon_i\}_{i=1}^{n}$ to estimate the integral of $f$, where $\{x_{i}\}_{i=1}^{n}$ and $\{\epsilon_i\}_{i=1}^{n}$ are independently and identically sampled from the uniform distribution on $\Omega$ and the normal distribution $\mathcal{N}(0,n^{-2\gamma})$ respectively. Assuming that $\gamma \in [0,\infty]$ and $s > 0$, we have
\begin{equation}
\label{eqn: lower bound for integral estimation}
\inf_{\hat{H} \in \mathcal{H}_{n}^{f}}\sup_{f \in C^{s}(\Omega)}\mathlarger{\mathlarger{\mathbb{E}}}_{\substack{\{x_i\}_{i=1}^{n},\{y_i\}_{i=1}^{n}}}\Bigg[\left|\hat{H}\Big(\{x_i\}_{i=1}^{n},\{y_i\}_{i=1}^{n}\Big) - I_f\right|\Bigg] \gtrsim n^{\max\{-\frac{1}{2}-\gamma, -\frac{1}{2}-\frac{s}{d}\}}.   
\end{equation}
\end{theorem}

\begin{remark}
Functional estimation is a well-studied problem in the literature of nonparametric statistics. However, current information-theoretic lower bounds for functional estimation \citep{birge1995estimation,donoho1990minimax,donoho1988one,robins2008higher,jiao2015minimax,krishnamurthy2014nonparametric,tsybakov2004introduction,han2020optimal} assume a constant level of noise on the observed function values. One essential idea for proving these lower bounds is to leverage the existence of the observational noise, which enables us to upper bound the amount of information required to distinguish between two reduced hypotheses. In contrast, we provide a minimax lower bound that is applicable for noises at any level by constructing two priors with overlapping support and assigning distinct probabilities to the corresponding Bernoulli random variables, which separates the two hypotheses. A comprehensive proof of Theorem \ref{thm: lower bound for integral} is given in Appendix \ref{app: proof of lower bound in thm 4.1}.
\end{remark}

\subsection{Optimal Nonparametric Regression-Adjusted Quadrature Rule}
\label{sec 4.2: k nearest neighbor is optimal for regression adjusted CV}

In the discussion below, we use the nearest-neighbor method as an example. For any $k \in \{1,2,\cdots,\frac{n}{2}\}$, the $k$-nearest neighbor estimator $\hat{f}_{k\text{-NN}}$ of $f$ is given by $\hat{f}_{k\text{-NN}}(z) := \frac{1}{k}\sum_{j=1}^{k}y_{i^{(z)}_j} $, where $\{x_{i^{(z)}_j}\}_{j=1}^{\frac{n}{2}}$ is a permutation of the quadrature points $\{x_{i}\}_{i=1}^{\frac{n}{2}}$ such that $||x_{i^{(z)}_1}-z|| \leq ||x_{i^{(z)}_2}-z|| \leq \cdots \leq||x_{i^{(z)}_{\frac{n}{2}}} - z||$ holds for any $z \in \Omega$. Moreover, we use $\mathcal{T}_{k,z} := \{x_{i^{(z)}_{j}}\}_{j=1}^{k}$ to denote the collection of the $k$ nearest neighbors of $z$ among $\{x_{i}\}_{i=1}^{\frac{n}{2}}$ for any $z \in \Omega$. For any $1 \leq i \leq \frac{n}{2}$, we take $D_i \subset \Omega$ to be the region formed by all the points whose $k$ nearest neighbors contain $x_i$, \emph{i.e,} $D_i := \Big\{z \in \Omega: x_i \in \mathcal{T}_{k,z}\Big\}$. Our estimator $\hat{H}_{k\text{-NN}}$ can be formally represented as 
$$
\hat{H}_{k\text{-NN}}\Big(\{x_i\}_{i=1}^{n},\{y_i\}_{i=1}^{n}\Big) = \underbrace{\sum_{i=1}^{\frac{n}{2}}\frac{V(D_i)}{k}y_{i}}_{\int_{\Omega}\hat{f}_{k\text{-NN}}(x)dx}+  \underbrace{\frac{2}{n}\sum_{i=\frac{n}{2}+1}^{n}y_i - \frac{2}{n}\sum_{i=\frac{n}{2}+1}^{n}\Big(\frac{1}{k}\sum_{j=1}^{\frac{n}{2}}\mathbbm{1}\{x_i \in D_j\}y_j\Big)}_{\frac{2}{n}\sum_{i=\frac{n}{2}+1}^{n}\Big(y_i-\hat{f}_{k\text{-NN}}(x_i)\Big)}.
$$

In the following theorem, we present an upper bound on the expected risk of the estimator $\hat{H}_{k\text{-NN}}$:

\begin{theorem}[Matching Upper Bound for Integral Estimation]
\label{thm: upper bound for integral}
Let $\{x_{i}\}_{i=1}^{n}$ be $n$ quadrature points independently and identically sampled from the uniform distribution on $\Omega$ and $\{y_i:= f(x_i) +\epsilon_i\}_{i=1}^n$ be the corresponding $n$ noisy observations of $f \in C^s(\Omega)$, where $\{\epsilon_i\}_{i=1}^{n}$ are independently and identically sampled from the normal distribution $\mathcal{N}(0,n^{-2\gamma})$. Assuming that $\gamma \in [0,\infty]$ and $s\in(0,1)$, we have that there exists $k \in \mathbb{N}$ such that the estimator $\hat{H}_{k\text{-NN}}$ constructed above satisfies 
\begin{equation}
\label{eq: upper bound (final)}
\mathlarger{\mathlarger{\mathbb{E}}}_{\substack{\{x_i\}_{i=1}^{n},\{y_i\}_{i=1}^{n}}}\Bigg[\left|\hat{H}_{k\text{-NN}}\Big(\{x_i\}_{i=1}^{n},\{y_i\}_{i=1}^{n}\Big) - I_f\right|\Bigg] \lesssim n^{\max\{-\frac{1}{2}-\gamma, -\frac{1}{2}-\frac{s}{d}\}}.    
\end{equation}
\end{theorem}

\begin{remark}
Our upper bound in Theorem \ref{thm: upper bound for integral} matches our minimax lower bound in Theorem \ref{thm: lower bound for integral}, which indicates that the regression-adjusted quadrature rule associated with the nearest neighbor estimator is minimax optimal.  When the noise level is high ($\gamma<\frac{s}{d}$), the control variate helps to improve the rate from $n^{-\frac{1}{2}}$ (the Monte Carlo rate) to $n^{-\frac{1}{2}-\gamma}$ via eliminating \textbf{\emph{all}} the effects of simulating the smooth function.  When the noise level is low ($\gamma>\frac{s}{d}$), we show that our estimator $\hat{H}_{k\text{-NN}}$ can achieve the optimal rate of quadrature rules \citep{novak2016some}. We defer a complete proof of Theorem \ref{thm: upper bound for integral} to Appendix \ref{app: proof of upper bound in thm 4.2}. 

\end{remark}

\section{Discussion and Conclusion}

In this paper, we have investigated whether a non-parametric regression-adjusted control variate can improve the rate of estimating functionals and if it is minimax optimal. Using the Sobolev Embedding Theorem, we discover that the existence of rare and extreme events will change the answer to this question. We show that when rare and extreme events are present, using a non-parametric machine learning algorithm as a control variate does not help, and truncated Monte Carlo is minimax optimal. Investigating how to apply importance sampling under this scenario may be of future interest. Also, the study of how regression-adjusted control variates adapt to the noise level for non-linear functionals \citep{han2020estimation, lepski1999estimation} is left as future work. Another interesting direction is to analyze how to use the data distribution's information \citep{oates2016control,oates2017control} to achieve both better computational trackability and convergence rate \citep{oates2019convergence}.

\begin{acknowledgement}
 Yiping Lu is supported by the Stanford Interdisciplinary Graduate
Fellowship (SIGF). Jose Blanchet is supported in part by the Air Force Office of Scientific Research under award number FA9550-20-1-0397. Lexing Ying is supported is supported by National Science Foundation under award DMS-2208163.
\end{acknowledgement}


\newpage
\printbibliography

@article{angelopoulos2023prediction,
  title={Prediction-Powered Inference},
  author={Angelopoulos, Anastasios N and Bates, Stephen and Fannjiang, Clara and Jordan, Michael I and Zrnic, Tijana},
  journal={arXiv preprint arXiv:2301.09633},
  year={2023}
}

@article{siegel2022high,
  title={High-order approximation rates for shallow neural networks with cosine and ReLU$^k$ activation functions},
  author={Siegel, Jonathan W and Xu, Jinchao},
  journal={Applied and Computational Harmonic Analysis},
  volume={58},
  pages={1--26},
  year={2022},
  publisher={Elsevier}
}

@inproceedings{meyer2021hutch++,
  title={Hutch++: Optimal stochastic trace estimation},
  author={Meyer, Raphael A and Musco, Cameron and Musco, Christopher and Woodruff, David P},
  booktitle={Symposium on Simplicity in Algorithms (SOSA)},
  pages={142--155},
  year={2021},
  organization={SIAM}
}

@article{sobczyk2022approximate,
  title={Approximate Euclidean lengths and distances beyond Johnson-Lindenstrauss},
  author={Sobczyk, Aleksandros and Luisier, Mathieu},
  journal={arXiv preprint arXiv:2205.12307},
  year={2022}
}

@article{liu2020regression,
  title={Regression-adjusted average treatment effect estimates in stratified randomized experiments},
  author={Liu, Hanzhong and Yang, Yuehan},
  journal={Biometrika},
  volume={107},
  number={4},
  pages={935--948},
  year={2020},
  publisher={Oxford University Press}
}

@article{bakhvalov2015approximate,
  title={On the approximate calculation of multiple integrals},
  author={Bakhvalov, Nikolai Sergeevich},
  journal={Journal of Complexity},
  volume={31},
  number={4},
  pages={502--516},
  year={2015},
  publisher={Elsevier}
}

@book{asmussen2007stochastic,
  title={Stochastic simulation: algorithms and analysis},
  author={Asmussen, S{\o}ren and Glynn, Peter W},
  volume={57},
  year={2007},
  publisher={Springer}
}

@article{davidson1992regression,
  title={Regression-based methods for using control variates in Monte Carlo experiments},
  author={Davidson, Russell and MacKinnon, James G},
  journal={Journal of Econometrics},
  volume={54},
  number={1-3},
  pages={203--222},
  year={1992},
  publisher={Elsevier}
}

@article{han2020optimal,
  title={Optimal rates of entropy estimation over Lipschitz balls},
  author={Han, Yanjun and Jiao, Jiantao and Weissman, Tsachy and Wu, Yihong},
  volume = {48},
journal = {The Annals of Statistics},
number = {6},
publisher = {Institute of Mathematical Statistics},
pages = {3228 -- 3250},
  year={2020}
}

@article{jiao2015minimax,
  title={Minimax estimation of functionals of discrete distributions},
  author={Jiao, Jiantao and Venkat, Kartik and Han, Yanjun and Weissman, Tsachy},
  journal={IEEE Transactions on Information Theory},
  volume={61},
  number={5},
  pages={2835--2885},
  year={2015},
  publisher={IEEE}
}

@article{lepski1999estimation,
  title={On estimation of the $L_{r}$ norm of a regression function},
  author={Lepski, Oleg and Nemirovski, Arkady and Spokoiny, Vladimir},
  journal={Probability theory and related fields},
  volume={113},
  pages={221--253},
  year={1999},
  publisher={Springer}
}

@article{han2020estimation,
  title={On estimation of $L_{r}$-norms in Gaussian white noise models},
  author={Han, Yanjun and Jiao, Jiantao and Mukherjee, Rajarshi},
  journal={Probability Theory and Related Fields},
  volume={177},
  number={3-4},
  pages={1243--1294},
  year={2020},
  publisher={Springer}
}

@book{biau2015lectures,
  title={Lectures on the nearest neighbor method},
  author={Biau, G{\'e}rard and Devroye, Luc},
  volume={246},
  year={2015},
  publisher={Springer}
}

@article{tsybakov2004introduction,
  title={Introduction to nonparametric estimation, 2009},
  author={Tsybakov, Alexandre B},
  journal={URL https://doi. org/10.1007/b13794. Revised and extended from the},
  volume={9},
  number={10},
  year={2004}
}

@article{robins2008higher,
  title={Higher order influence functions and minimax estimation of nonlinear functionals},
  author={Robins, James and Li, Lingling and Tchetgen, Eric and van der Vaart, Aad and others},
  journal={Probability and statistics: essays in honor of David A. Freedman},
  volume={2},
  pages={335--421},
  year={2008},
  publisher={Institute of Mathematical Statistics Beachwood, OH}
}

@article{lin2017randomized,
  title={Randomized estimation of spectral densities of large matrices made accurate},
  author={Lin, Lin},
  journal={Numerische Mathematik},
  volume={136},
  pages={183--213},
  year={2017},
  publisher={Springer}
}

@article{birge1995estimation,
  title={Estimation of integral functionals of a density},
  author={Birg{\'e}, Lucien and Massart, Pascal},
  journal={The Annals of Statistics},
  volume={23},
  number={1},
  pages={11--29},
  year={1995},
  publisher={Institute of Mathematical Statistics}
}

@inproceedings{krishnamurthy2014nonparametric,
  title={Nonparametric estimation of renyi divergence and friends},
  author={Krishnamurthy, Akshay and Kandasamy, Kirthevasan and Poczos, Barnabas and Wasserman, Larry},
  booktitle={International Conference on Machine Learning},
  pages={919--927},
  year={2014},
  organization={PMLR}
}

@article{bach2017equivalence,
  title={On the equivalence between kernel quadrature rules and random feature expansions},
  author={Bach, Francis},
  journal={The Journal of Machine Learning Research},
  volume={18},
  number={1},
  pages={714--751},
  year={2017},
  publisher={JMLR. org}
}

@article{karvonen2018fully,
  title={Fully symmetric kernel quadrature},
  author={Karvonen, Toni and Sarkka, Simo},
  journal={SIAM Journal on Scientific Computing},
  volume={40},
  number={2},
  pages={A697--A720},
  year={2018},
  publisher={SIAM}
}

@article{kanagawa2016convergence,
  title={Convergence guarantees for kernel-based quadrature rules in misspecified settings},
  author={Kanagawa, Motonobu and Sriperumbudur, Bharath K and Fukumizu, Kenji},
  journal={Advances in Neural Information Processing Systems},
  volume={29},
  year={2016}
}

@article{belhadji2019kernel,
  title={Kernel quadrature with DPPs},
  author={Belhadji, Ayoub and Bardenet, R{\'e}mi and Chainais, Pierre},
  journal={Advances in Neural Information Processing Systems},
  volume={32},
  year={2019}
}

@article{bardenet2020monte,
  title={Monte Carlo with Determinantal Point Processes},
  author={Bardenet, R{\'e}mi and Hardy, Adrien},
  journal={Annals of Applied Probability},
  year={2020}
}

@article{belhadji2021analysis,
  title={An analysis of Ermakov-Zolotukhin quadrature using kernels},
  author={Belhadji, Ayoub},
  journal={Advances in Neural Information Processing Systems},
  volume={34},
  pages={27278--27289},
  year={2021}
}

@article{gautier2019two,
  title={On two ways to use determinantal point processes for Monte Carlo integration},
  author={Gautier, Guillaume and Bardenet, R{\'e}mi and Valko, Michal},
  journal={Advances in Neural Information Processing Systems},
  volume={32},
  year={2019}
}

@article{hayakawa2021positively,
  title={Positively weighted kernel quadrature via subsampling},
  author={Hayakawa, Satoshi and Oberhauser, Harald and Lyons, Terry},
  journal={arXiv preprint arXiv:2107.09597},
  year={2021}
}

@article{hayakawa2023sampling,
  title={Sampling-based Nystr\"{o}m Approximation and Kernel Quadrature},
  author={Hayakawa, Satoshi and Oberhauser, Harald and Lyons, Terry},
  journal={arXiv preprint arXiv:2301.09517},
  year={2023}
}

@article{chen2012super,
  title={Super-samples from kernel herding},
  author={Chen, Yutian and Welling, Max and Smola, Alex},
  journal={arXiv preprint arXiv:1203.3472},
  year={2012}
}

@inproceedings{lacoste2015sequential,
  title={Sequential kernel herding: Frank-Wolfe optimization for particle filtering},
  author={Lacoste-Julien, Simon and Lindsten, Fredrik and Bach, Francis},
  booktitle={Artificial Intelligence and Statistics},
  pages={544--552},
  year={2015},
  organization={PMLR}
}

@article{holzmuller2023convergence,
  title={Convergence Rates for Non-Log-Concave Sampling and Log-Partition Estimation},
  author={Holzm{\"u}ller, David and Bach, Francis},
  journal={arXiv preprint arXiv:2303.03237},
  year={2023}
}

@article{huszar2012optimally,
  title={Optimally-weighted herding is Bayesian quadrature},
  author={Husz{\'a}r, Ferenc and Duvenaud, David},
  journal={arXiv preprint arXiv:1204.1664},
  year={2012}
}

@article{o1991bayes,
  title={Bayes--hermite quadrature},
  author={O'Hagan, Anthony},
  journal={Journal of statistical planning and inference},
  volume={29},
  number={3},
  pages={245--260},
  year={1991},
  publisher={Elsevier}
}

@article{kanagawa2019convergence,
  title={Convergence guarantees for adaptive Bayesian quadrature methods},
  author={Kanagawa, Motonobu and Hennig, Philipp},
  journal={Advances in Neural Information Processing Systems},
  volume={32},
  year={2019}
}

@inproceedings{oates2016control,
  title={Control functionals for quasi-Monte Carlo integration},
  author={Oates, Chris and Girolami, Mark},
  booktitle={Artificial Intelligence and Statistics},
  pages={56--65},
  year={2016},
  organization={PMLR}
}

@article{oates2017control,
  title={Control functionals for Monte Carlo integration},
  author={Oates, Chris J and Girolami, Mark and Chopin, Nicolas},
  journal={Journal of the Royal Statistical Society. Series B (Statistical Methodology)},
  pages={695--718},
  year={2017},
  publisher={JSTOR}
}

@article{shi2022gradient,
  title={Gradient estimation with discrete Stein operators},
  author={Shi, Jiaxin and Zhou, Yuhao and Hwang, Jessica and Titsias, Michalis and Mackey, Lester},
  journal={Advances in Neural Information Processing Systems},
  volume={35},
  pages={25829--25841},
  year={2022}
}

@article{assaraf1999zero,
  title={Zero-variance principle for Monte Carlo algorithms},
  author={Assaraf, Roland and Caffarel, Michel},
  journal={Physical review letters},
  volume={83},
  number={23},
  pages={4682},
  year={1999},
  publisher={APS}
}

@article{south2018regularised,
  title={Regularised zero-variance control variates},
  author={South, Leah F and Oates, CJ and Mira, A and Drovandi, C},
  journal={arXiv preprint arXiv:1811.05073},
  year={2018}
}

@article{mira2013zero,
  title={Zero variance markov chain monte carlo for bayesian estimators},
  author={Mira, Antonietta and Solgi, Reza and Imparato, Daniele},
  journal={Statistics and Computing},
  volume={23},
  pages={653--662},
  year={2013},
  publisher={Springer}
}

@article{liu2017action,
  title={Action-depedent Control Variates for Policy Optimization via Stein's Identity},
  author={Liu, Hao and Feng, Yihao and Mao, Yi and Zhou, Dengyong and Peng, Jian and Liu, Qiang},
  journal={arXiv preprint arXiv:1710.11198},
  year={2017}
}

@article{oates2019convergence,
  title={Convergence rates for a class of estimators based on Stein’s method},
  author={Oates, Chris J and Cockayne, Jon and Briol, Fran{\c{c}}ois-Xavier and Girolami, Mark},
  volume = {25},
journal = {Bernoulli},
number = {2},
publisher = {Bernoulli Society for Mathematical Statistics and Probability},
pages = {1141 -- 1159},
  year={2019}
}

@article{romano2019conformalized,
  title={Conformalized quantile regression},
  author={Romano, Yaniv and Patterson, Evan and Candes, Emmanuel},
  journal={Advances in neural information processing systems},
  volume={32},
  year={2019}
}

@article{krieg2022recovery,
  title={Recovery of Sobolev functions restricted to iid sampling},
  author={Krieg, David and Novak, Erich and Sonnleitner, Mathias},
  journal={Mathematics of Computation},
  volume={91},
  number={338},
  pages={2715--2738},
  year={2022}
}

@article{hinrichs2022lower,
  title={Lower bounds for integration and recovery in $L_2$},
  author={Hinrichs, Aicke and Krieg, David and Novak, Erich and Vyb{\'\i}ral, Jan},
  journal={Journal of Complexity},
  volume={72},
  pages={101662},
  year={2022},
  publisher={Elsevier}
}

@article{dwivedi2021generalized,
  title={Generalized kernel thinning},
  author={Dwivedi, Raaz and Mackey, Lester},
  journal={arXiv preprint arXiv:2110.01593},
  year={2021}
}

@inproceedings{chen2018stein,
  title={Stein points},
  author={Chen, Wilson Ye and Mackey, Lester and Gorham, Jackson and Briol, Fran{\c{c}}ois-Xavier and Oates, Chris},
  booktitle={International Conference on Machine Learning},
  pages={844--853},
  year={2018},
  organization={PMLR}
}

@article{dwivedi2021kernel,
  title={Kernel thinning},
  author={Dwivedi, Raaz and Mackey, Lester},
  journal={arXiv preprint arXiv:2105.05842},
  year={2021}
}

@article{hinrichs2014curse,
  title={The curse of dimensionality for numerical integration of smooth functions},
  author={Hinrichs, Aicke and Novak, Erich and Ullrich, Mario and Wo{\'z}niakowski, H},
  journal={Mathematics of Computation},
  volume={83},
  number={290},
  pages={2853--2863},
  year={2014}
}

@article{donoho1990minimax,
  title={Minimax quadratic estimation of a quadratic functional},
  author={Donoho, David L and Nussbaum, Michael},
  journal={Journal of Complexity},
  volume={6},
  number={3},
  pages={290--323},
  year={1990},
  publisher={Elsevier}
}

@article{donoho1988one,
  title={One-sided inference about functionals of a density},
  author={Donoho, David L},
  journal={The Annals of Statistics},
  pages={1390--1420},
  year={1988},
  publisher={JSTOR}
}

@article{donoho1991geometrizing3,
  title={Geometrizing rates of convergence, III},
  author={Donoho, David L and Liu, Richard C},
  journal={The Annals of Statistics},
  pages={668--701},
  year={1991},
  publisher={JSTOR}
}

@article{donoho1991geometrizing2,
  title={Geometrizing rates of convergence, II},
  author={Donoho, David L and Liu, Richard C},
  journal={The Annals of Statistics},
  pages={633--667},
  year={1991},
  publisher={JSTOR}
}

@article{heinrich2018complexity,
  title={On the complexity of computing the $L_{q}$ norm},
  author={Heinrich, Stefan},
  journal={Journal of Complexity},
  volume={49},
  pages={1--26},
  year={2018},
  publisher={Elsevier}
}

@article{heinrich2009randomized2,
  title={Randomized approximation of Sobolev embeddings, II},
  author={Heinrich, Stefan},
  journal={Journal of Complexity},
  volume={25},
  number={5},
  pages={455--472},
  year={2009},
  publisher={Elsevier}
}

@article{heinrich2009randomized3,
  title={Randomized approximation of Sobolev embeddings, III},
  author={Heinrich, Stefan},
  journal={Journal of Complexity},
  volume={25},
  number={5},
  pages={473--507},
  year={2009},
  publisher={Elsevier}
}

@article{mathe1991random,
  title={Random approximation of Sobolev embeddings},
  author={Math{\'e}, Peter},
  journal={Journal of Complexity},
  volume={7},
  number={3},
  pages={261--281},
  year={1991},
  publisher={Elsevier}
}

@book{novak2008tractability2,
  title={Tractability of Multivariate Problems: Standard Information for Functionals},
  author={Novak, Erich and Wo{\'z}niakowski, Henryk},
  volume={2},
  year={2008},
  publisher={European Mathematical Society}
}

@article{krieg2020random,
  title={Random points are optimal for the approximation of Sobolev functions},
  author={Krieg, David and Sonnleitner, Mathias},
  journal={arXiv preprint arXiv:2009.11275},
  year={2020}
}

@article{hinrichs2020power,
  title={On the power of random information},
  author={Hinrichs, Aicke and Krieg, David and Novak, Erich and Prochno, Joscha and Ullrich, Mario},
  journal={Multivariate Algorithms and information-based complexity},
  volume={27},
  pages={43--64},
  year={2020},
  publisher={De Gruyter Berlin/Boston}
}

@article{novak2016some,
  title={Some results on the complexity of numerical integration},
  author={Novak, Erich},
  journal={Monte Carlo and Quasi-Monte Carlo Methods: MCQMC, Leuven, Belgium, April 2014},
  pages={161--183},
  year={2016},
  publisher={Springer}
}

@article{novak2008tractability1,
  title={Tractability of Multivariate Problems, Volume I: Linear Information, European Math},
  author={Novak, E and Wozniakowski, H},
  journal={Soc., Z{\"u}rich},
  volume={2},
  number={3},
  year={2008}
}

@article{wendland2001local,
  title={Local polynomial reproduction and moving least squares approximation},
  author={Wendland, Holger},
  journal={IMA Journal of Numerical Analysis},
  volume={21},
  number={1},
  pages={285--300},
  year={2001},
  publisher={OUP}
}

@book{wendland2004scattered,
  title={Scattered data approximation},
  author={Wendland, Holger},
  volume={17},
  year={2004},
  publisher={Cambridge university press}
}

@article{hickernell2005control,
  title={Control variates for quasi-Monte Carlo},
  author={Hickernell, Fred J and Lemieux, Christiane and Owen, Art B},
  volume = {20},
  journal = {Statistical Science},
  number = {1},
  publisher = {Institute of Mathematical Statistics},
  pages = {1 -- 31},
  year={2005}
}

@article{traub1994information,
  title={Information-based complexity},
  author={Traub, Joseph F and Wasilkowski, GW and Wozniakowski, H and Novak, Erich},
  journal={SIAM Review},
  volume={36},
  number={3},
  pages={514--514},
  year={1994},
  publisher={Philadelphia, Society for Industrial and Applied Mathematics.}
}

@book{adams2003sobolev,
  title={Sobolev spaces},
  author={Adams, Robert A and Fournier, John JF},
  year={2003},
  publisher={Elsevier}
}

@book{novak2006deterministic,
  title={Deterministic and stochastic error bounds in numerical analysis},
  author={Novak, Erich},
  volume={1349},
  year={2006},
  publisher={Springer}
}
\appendix
\newpage
\section*{Appendix}
The appendix is organized as follows: 
\begin{itemize}
    \item In Appendix A, we list some notations and standard lemmas used in our proofs.
    \item Appendix B contains a comprehensive proof of the information-theoretic lower bound on the estimation of $q$-th moments, which is established in Theorem \ref{thm: lower bound q-th moment}.
    
    \item In Appendix C, we provide a detailed proof of Theorem \ref{thm: upper bound control variate} and \ref{thm: upper bound truncate MC}, which gives us the minimax optimal upper bound on estimating $q$-th moments.
    
    \item Appendix D consists of our proof for the information-theoretic lower bounds and minimax optimal upper bounds on integral estimation and function estimation, which are listed in Theorem \ref{thm: lower bound for integral} and \ref{thm: upper bound for integral}. 
\end{itemize}

\section{Preliminaries and Basic Tools}
\subsection{Preliminaries}
This subsection is devoted to presenting some basic notations used in our proofs. For any fixed convex function $f:\mathbb{R}^{+} \rightarrow \mathbb{R}$ satisfying $f(1) = 0$, we use $D_f(\cdot || \cdot)$ to denote the corresponding $f$-divergence, \emph{i.e,} $D_f(P || Q) = \int_{\mathcal{Y}}f\Big(\frac{dP}{dQ}\Big)dQ$ for any two probability distributions $P$ and $Q$ over some fixed space $\mathcal{Y}$. In particular, when $f(x) = \frac{1}{2}|x-1|$, $D_f(\cdot || \cdot)$ is the total variation (TV) distance $TV(\cdot || \cdot)$. When $f(x)=x\log x$, $D_f(\cdot || \cdot)$ coincides with the Kullback–Leibler (KL) divergence $KL(\cdot || \cdot)$. Moreover, for any $a \in \mathbb{R}$, we use $\delta_a(\cdot)$ to denote the Dirac delta distribution at point $a$, \emph{i.e,} $\int_{-\infty}^{\infty}f(x)\delta_{a}(x) = f(a)$ for any function $f:\mathbb{R} \rightarrow \mathbb{R}$. 

\subsection{Basic Lemmas}
In this subsection, we list some basic lemmas that serve as essential tools in our proofs.  
\begin{lemma}[Sobolev Embedding Theorem \citep{adams2003sobolev}]
\label{lem: sobolev embedding}
For some fixed dimension $d \in \mathbb{N}$, we have that\\
(I) For any $s,t \in \mathbb{N}_0$ and $p,q \in \mathbb{R}$ satisfying $s > t$, $p <d$ and $1 \leq p < q \leq \infty$, we have $W^{s,p}(\mathbb{R}^d) \subseteq W^{t,q}(\mathbb{R}^d)$ when the relation $\frac{1}{p}-\frac{s}{d} = \frac{1}{q} - \frac{t}{d}$ holds. In the special case when $t=0$, we have $W^{s,p}(\mathbb{R}^d) \subseteq L^q(\mathbb{R}^d)$ for any $s \in \mathbb{N}$ and $p,q \in \mathbb{R}$ satisfying $1 \leq p < q \leq \infty$ and $\frac{1}{p}-\frac{s}{d} \leq \frac{1}{q}$.\\
(II) For any $\alpha \in (0,1)$, let $\beta = \frac{d}{1-\alpha} \in (d,\infty]$. Then we have $C^{1}(\mathbb{R}^d) \cap W^{1,\beta}(\mathbb{R}^d) \subseteq C^{\alpha}(\mathbb{R}^d)$.
\end{lemma}

\begin{lemma}[H\"older's Inequality]
\label{lem: holder ineq}
For any fixed domain $\Omega$ and $p,q \in [1,\infty]$ satisfying $\frac{1}{p}+\frac{1}{q} = 1$, we have that $||fg||_{L^1(\Omega)} \leq ||f||_{L^p(\Omega)} ||g||_{L^q(\Omega)}$ holds for any $f \in L^p(\Omega), g \in L^q(\Omega)$.
\end{lemma}

\begin{lemma}[Hoeffding's Inequality]
\label{lem: hoeffding ineq}
Let $X_1,X_2,\cdots,X_n$ be independent random variables satisfying $X_i \in [a_i,b_i]$ for any $1 \leq i \leq n$. Then for any $\epsilon >0$, the sum $S_n := \sum_{i=1}^{n}X_i$ of these $n$ random variables satisfies the following inequality:
\begin{equation}
\begin{aligned}
\mathbb{P}(S_n \geq \mathbb{E}[S_n] + t) &\leq \exp\Big(-\frac{2t^2}{\sum_{i=1}^{n}(b_i-a_i)^2}\Big)\\
\mathbb{P}(S_n \leq \mathbb{E}[S_n] - t) &\leq \exp\Big(-\frac{2t^2}{\sum_{i=1}^{n}(b_i-a_i)^2}\Big)
\end{aligned}
\end{equation}
\end{lemma}

\begin{lemma}[Data Processing Inequality]
\label{lem: data processing ineq}
Given some Markov Chain $X \rightarrow Z$, where $X$ and $Z$ are two random variables the measurable spaces $(\mathcal{X},\mu)$ and $(\mathcal{Z},\nu)$ respectively. Let $K$ be the transition kernel of the Markov Chain above, \emph{i.e,} for any $x \in \mathcal{X}$, the probability distribution of $Z$ is given by $K(\cdot,x)$ when conditioned on $X=x$. For any two fixed two distributions $P,Q$ over $\mathcal{X}$ with probability density functions $p,q$, we use $K_P(\cdot)$ and $K_Q(\cdot)$ to denote the corresponding marginal distributions respectively, \emph{i.e,} $K_P(\cdot) := \int_{\mathcal{X}} K(\cdot,x)p(x)d\mu(x)$ and $K_Q(\cdot) = \int_{\mathcal{X}} K(\cdot,x)q(x)d\mu(x)$. Then we have $D_f(K_P || K_Q) \leq D_f(P || Q)$ holds for any $f$-divergence $D_f(\cdot || \cdot)$.
\end{lemma}

\section{Proof of Lower Bounds in Section \ref{sec 2: lower bound on q-moment estimation}}

\subsection{A Key Lemma for Building Minimax Optimal Lower Bounds}
In this subsection, we firstly present the method of two fuzzy hypotheses, which turns out to be the most essential tool for establishing all the minimax optimal lower bounds in our paper, before giving our complete proof of Theorem \ref{thm: lower bound q-th moment}. 
\begin{lemma}[Method of Two Fuzzy Hypotheses: Theorem 2.15 (i), \citep{tsybakov2004introduction}]
\label{lem: method of two fuzzy hypo}
Let $F: \mathit{\Theta} \rightarrow \mathbb{R}$ be some continuous functional defined on the measurable space $(\mathit{\Theta},\mathcal{U})$ and taking values in $(\mathbb{R},\mathcal{B}(\mathbb{R}))$, where $\mathcal{B}(\mathbb{R})$ denotes the Borel $\sigma$-algebra on $\mathbb{R}$. Suppose that each parameter $\theta \in \mathit{\Theta}$ is associated with a distribution $\mathbb{P}_{\theta}$, which together form a collection $\{\mathbb{P}_{\theta}:\theta \in \mathit{\Theta}\}$ of distributions. 

For any fixed $\theta \in \mathit{\Theta}$, assume that our observation $\mathbf{X}$ is distributed as $\mathbb{P}_{\theta}$. Let $\hat{F}$ be an arbitrary estimator of $F(\theta)$ based on $\mathbf{X}$. Let $\mu_0,\mu_1$ be two prior measures on $\mathit{\Theta}$. Assume that there exist constants $c \in \mathbb{R}, \Delta \in (0,\infty)$ and $ \beta_0,\beta_1 \in [0,1)$, such that:
\begin{equation}
\label{eq: separated hypothesis (two prior)}
\begin{aligned}
\mu_0(\theta \in \mathit{\Theta}:F(\theta) \leq c-\Delta) &\geq 1-\beta_0, \\ 
\mu_1(\theta \in \mathit{\Theta}: F(\theta) \geq c+\Delta) &\geq 1-\beta_1.
\end{aligned}    
\end{equation}
For $j \in \{0,1\}$, we use $\mathbb{P}_j(\cdot) := \int \mathbb{P}_{\theta}(\cdot)\mu_j(d\theta)$ to denote the marginal distribution $\mathbb{P}_j$ associated with the prior distribution $\mu_j$. Then we have the following lower bound: 
\begin{equation}
\label{eq: bound on misclassification probability (two prior)}
\inf_{\hat{F}}\sup_{\theta \in \Theta}\mathbb{P}_{\theta}(|\hat{F}-F(\theta)| \geq \Delta) \geq \frac{1-TV(\mathbb{P}_0 || \mathbb{P}_1)-\beta_0-\beta_1}{2}.    
\end{equation}
\end{lemma}

\subsection{Proof of Theorem \ref{thm: lower bound q-th moment} (Information-Theoretic Lower Bound on Moment Estimation)}
\label{app: proof of lower bound in sec 2}
In this subsection, we give a detailed proof of the two minimax lower bounds established in Theorem \ref{thm: lower bound q-th moment} above via the method of two fuzzy hypotheses (Lemma \ref{lem: method of two fuzzy hypo}). We start off by introducing some preliminary tools used in our proof. Consider the function $K_0$ defined as follows:
\begin{equation}
\label{eqn: def of compact C_infty function K_0 on [-1,1]^d}
K_0(x) := \prod_{i=1}^{d}\exp\Big(-\frac{1}{1-x_{i}^2}\Big)\mathbbm{1}(|x_{i}| \leq 1), \forall \ x=(x_1,x_2,\cdots,x_d) \in \mathbb{R}^d.    
\end{equation}
Moreover, we pick some function $K$ satisfying  
\begin{equation}
\label{eqn: def of K, rescaled K_0 on [-1/2,1/2]^d}
K(x) := K_0(2x), \ \forall \ x \in \mathbb{R}^d,    
\end{equation} 
From our construction of $K$ and $K_0$ above, we have that $K_0$ is in $C^{\infty}(\mathbb{R}^d)$ and compactly supported on $[-\frac{1}{2},\frac{1}{2}]^d$. Furthermore, we set $m= (200n)^{\frac{1}{d}}$ and divide the domain $\Omega$ into $m^d$ small cubes $\Omega_1,\Omega_2,\cdots,\Omega_{m^d}$, each of which has side length $m^{-1}$. For any $1 \leq j \leq m^d$, we use $c_j$ to denote the center of the cube $\Omega_j$. Similar to the proof sketch of Theorem \ref{thm: lower bound q-th moment}, below we again use $w_{p}$ to denote the discrete random variable satisfying $\mathbb{P}(w_p = -1) = p$ and $\mathbb{P}(w_p = 1)=1-p$ for any $p \in (0,1)$. Furthermore, we use $\vec{x}:=(x_1,x_2,\cdots,x_n)$ and $\vec{y}:=(y_1,y_2,\cdots,y_n)$ to denote the two $n$-dimensional vectors formed by the quadrature points and observed function values, After introducing all preliminaries above, let's present the essential parts of our proof. Given that our lower bound in Theorem \ref{thm: lower bound q-th moment} consists of two terms, our proof is also divided into two parts:

(Case I) For the first lower bound in (\ref{eq: lower bound q-th moment}), let's consider two functions $g_0$ and $g_1$ defined as follows:
\begin{equation}
\begin{aligned}
g_0(x) &\equiv 0 \ (\forall \ x \in \Omega), \\
g_1(x) &= 
\begin{cases}
m^{-s+\frac{d}{p}}K(m(x-c_1)) \ (x \in \Omega_1), \\
0 \ (\text{otherwise}).
\end{cases}    
\end{aligned}
\end{equation}
Clearly we have $g_0 \in W^{s,p}(\Omega)$ and $I_{g_0}^q =0$. Now let's verify that $g_1 \in W^{s,p}(\Omega)$ for any $m$. Note that the following bound holds for any $t \in \mathbb{N}_{0}^d$ satisfying $|t| \leq s$:
\begin{align*}
||D^{t}g_1||_{L^p(\Omega)}&= \Bigg(\int_{\Omega_1}\Big|m^{-s+\frac{d}{p}}m^{|t|}(D^{t}K)(m(x-c_1))\Big|^p dx\Bigg)^{\frac{1}{p}}\\
&= m^{|t|-s+\frac{d}{p}}\Bigg(\int_{[-\frac{1}{2},\frac{1}{2}]^d}\Big|(D^{t}K)(y)\Big|^p \frac{1}{m^d}dy\Bigg)^{\frac{1}{p}} = m^{|t|-s}||D^tK||_{L^p([-\frac{1}{2},\frac{1}{2}]^d)} \lesssim 1.   
\end{align*}
This implies $g_1 \in W^{s,p}(\Omega)$ for any $m$, as desired. Moreover, computing the $q$-th moment of $g_1$ yields
\begin{equation}
\label{eqn: thm 2.1 distance case 1}
\begin{aligned}
I_{g_1}^q &= \int_{\Omega}g_1^q(x)dx = \int_{\Omega_1}(m^{-s+\frac{d}{p}}K(m(x-c_1)))^qdx \\
&= m^{-q(s-\frac{d}{p})}\int_{[-\frac{1}{2},\frac{1}{2}]^d}(K(y))^q\frac{1}{m^d}dy = m^{-q(s-\frac{d}{p})-d}||K||_{L^q([-\frac{1}{2},\frac{1}{2}]^d)}^q.  
\end{aligned}
\end{equation}

Now let us take $\epsilon = \frac{1}{2}$ and pick two discrete measures $\mu_0,\mu_1$ supported on the finite set $\{g_0,g_1\} \subset W^{s,p}(\Omega)$ as below:
\begin{equation}
\label{eqn: def of two priors_1 in thm 2.1 case 1}
\begin{aligned}
\mu_0(\{g_0\}) &= \frac{1+\epsilon}{2}, \mu_0(\{g_1\}) = \frac{1-\epsilon}{2},\\  
\mu_1(\{g_0\}) &= \frac{1-\epsilon}{2}, \mu_1(\{g_1\}) = \frac{1+\epsilon}{2}.
\end{aligned}    
\end{equation}
On the one hand, by taking $c=\Delta =\frac{1}{2}I_{g_1}^q$ and $\beta_0=\beta_1=\frac{1-\epsilon}{2}$, we may use (\ref{eqn: def of two priors_1 in thm 2.1 case 1}) to deduce that 
\begin{equation}
\begin{aligned}
\mu_0(f \in W^{s,p}(\Omega): I_{f}^q &\leq c-\Delta) = \mu_0(I_{f}^{q} \leq 0) \geq \frac{1+\epsilon}{2}= 1-\beta_0, \\
\mu_1(f \in W^{s,p}(\Omega): I_{f}^q &\geq c+\Delta) =\mu_1(I_{f}^q \geq I_{g_1}^q) \geq \frac{1+\epsilon}{2}= 1-\beta_1.     
\end{aligned}    
\end{equation}
Hence, we have that (\ref{eq: separated hypothesis (two prior)}) holds true. On the other hand, recall that the quadrature points $\{x_1,\cdots,x_n\}$ are identical and independent samples from the uniform distribution on $\Omega$, which enables us to write the marginal distributions in an explicit form as follows: 
\begin{equation}
\label{eqn: marginal distributions thm 2.1 case 1}
\begin{aligned}
\mathbb{P}_0(\vec{x},\vec{y}) &= \Big(\frac{1+\epsilon}{2}\prod_{i:x_i \in \Omega_1}\delta_{0}(y_i) + \frac{1-\epsilon}{2}\prod_{i:x_i \in \Omega_1}\delta_{g_1(x_i)}(y_i)\Big) \cdot \prod_{j=2}^{m^d}\Big(\prod_{i:x_i \in \Omega_j}\delta_{0}(y_i)\Big),\\
\mathbb{P}_1(\vec{x},\vec{y}) &= \Big(\frac{1-\epsilon}{2}\prod_{i:x_i \in \Omega_1}\delta_{0}(y_i) + \frac{1+\epsilon}{2}\prod_{i:x_i \in \Omega_1}\delta_{g_1(x_i)}(y_i)\Big) \cdot \prod_{j=2}^{m^d}\Big(\prod_{i:x_i \in \Omega_j}\delta_{0}(y_i)\Big).   
\end{aligned}    
\end{equation}
In particular, we have $\mathbb{P}_0 = \mathbb{P}_1$ when the set $\{i:x_i \in \Omega_1\}$ is empty. Combing this fact with (\ref{eqn: marginal distributions thm 2.1 case 1}) above allows us to compute the KL divergence between $\mathbb{P}_0$ and $\mathbb{P}_1$ as below
\begin{equation}
\label{eqn: thm 2.1 case 1 upper bound on KL}
\begin{aligned}
&KL(\mathbb{P}_{0} || \mathbb{P}_{1}) = \int_{\Omega} \cdots\int_{\Omega}\Big(\int_{-\infty}^{\infty}\cdots\int_{-\infty}^{\infty}\log\Big(\frac{\mathbb{P}_0(\vec{x},\vec{y})}{\mathbb{P}_1(\vec{x},\vec{y})}\Big)\mathbb{P}_0(\vec{x},\vec{y})dy_1 \cdots dy_n\Big)dx_1 \cdots dx_n\\
&= \int_{\Omega} \cdots\int_{\Omega}\Bigg(\int_{-\infty}^{\infty}\cdots\int_{-\infty}^{\infty}\log\Big(\frac{\frac{1+\epsilon}{2}\prod_{i:x_i \in \Omega_1}\delta_{0}(y_i) + \frac{1-\epsilon}{2}\prod_{i:x_i \in \Omega_1}\delta_{g_1(x_i)}(y_i)}{\frac{1-\epsilon}{2}\prod_{i:x_i \in \Omega_1}\delta_{0}(y_i) + \frac{1+\epsilon}{2}\prod_{i:x_i \in \Omega_1}\delta_{g_1(x_i)}(y_i)}\Big)\\
&\cdot \Big(\frac{1+\epsilon}{2}\prod_{i:x_i \in \Omega_1}\delta_{0}(y_i) + \frac{1-\epsilon}{2}\prod_{i:x_i \in \Omega_1}\delta_{g_1(x_i)}(y_i)\Big) \cdot \Big(\prod_{j=2}^{m^d}\prod_{i:x_i \in \Omega_j}\delta_{0}(y_i)\Big)\prod_{i=1}^{n}dy_i\Bigg)\prod_{i=1}^{n}dx_i\\
&=\int_{\Omega} \cdots\int_{\Omega}\Bigg(\int_{-\infty}^{\infty}\cdots\int_{-\infty}^{\infty}\log\Big(\frac{\frac{1+\epsilon}{2}\prod_{i:x_i \in \Omega_1}\delta_{0}(y_i) + \frac{1-\epsilon}{2}\prod_{i:x_i \in \Omega_1}\delta_{g_1(x_i)}(y_i)}{\frac{1-\epsilon}{2}\prod_{i:x_i \in \Omega_1}\delta_{0}(y_i) + \frac{1+\epsilon}{2}\prod_{i:x_i \in \Omega_1}\delta_{g_1(x_i)}(y_i)}\Big)\\
&\cdot \Big(\frac{1+\epsilon}{2}\prod_{i:x_i \in \Omega_1}\delta_{0}(y_i) + \frac{1-\epsilon}{2}\prod_{i:x_i \in \Omega_1}\delta_{g_1(x_i)}(y_i)\Big) \prod_{i:x_i \in \Omega_1}dy_i\Bigg)\prod_{i=1}^{n}dx_i\\
&=\Big(\log\Big(\frac{1+\epsilon}{1-\epsilon}\Big)\frac{1+\epsilon}{2} + \log\Big(\frac{1-\epsilon}{1+\epsilon}\Big)\frac{1-\epsilon}{2}\Big)\mathbb{P}\Big(\{i:x_i \in \Omega_1\} \neq \varnothing\Big) \\
&= \epsilon\log\Big(\frac{1+\epsilon}{1-\epsilon}\Big)\mathbb{P}\Big(\{i:x_i \in \Omega_1\} \neq \varnothing\Big).
\end{aligned}    
\end{equation}

Moreover, since the probability that $\{i:x_i \in \Omega_1\} = \varnothing$ equals to $(\frac{m^d-1}{m^d})^n = (\frac{m^d-1}{m^d})^{\frac{m^d}{200}}$, we have
\begin{equation}
\label{eqn: thm 2.1 case 1 bound on probability}
\begin{aligned}
\mathbb{P}\Big(\{i:x_i \in \Omega_1\} \neq \varnothing \Big) = 1- \Big(1-\frac{1}{m^d}\Big)^{\frac{m^d}{200}} \leq 1-\Bigg(\frac{1}{e}\Big(1-\frac{1}{m^d}\Big)\Bigg)^{\frac{1}{200}} \leq 1-(2e)^{-\frac{1}{200}}.   
\end{aligned}    
\end{equation}

Now we may combine (\ref{eqn: thm 2.1 case 1 upper bound on KL}), (\ref{eqn: thm 2.1 case 1 bound on probability}) and Pinkser's inequality to upper bound the TV distance between $\mathbb{P}_0$ and $\mathbb{P}_1$ as below:
\begin{equation}
\label{thm 2.1 case 1: bound on TV}
TV(\mathbb{P}_0 || \mathbb{P}_1) \leq  \sqrt{\frac{1}{2}KL(\mathbb{P}_0 || \mathbb{P}_1)} \leq \sqrt{\frac{1-(2e)^{-\frac{1}{200}}}{2}\epsilon\log\Big(\frac{1+\epsilon}{1-\epsilon}\Big)} \leq \sqrt{\frac{3}{100}}\epsilon= \frac{\sqrt{3}}{10}\epsilon.   
\end{equation}

Finally, by substituting (\ref{eqn: thm 2.1 distance case 1}), (\ref{thm 2.1 case 1: bound on TV}), $\Delta =\frac{1}{2}I_{g_1}^q$ and $\beta_0=\beta_1=\frac{1-\epsilon}{2} = \frac{1}{4}$ into (\ref{eq: bound on misclassification probability (two prior)}) and applying Markov's inequality, we obtain the final lower bound
\begin{equation}
\label{eqn: thm 2.1 case 1 final bound}
\begin{aligned}
&\inf_{\hat{H}^q \in \mathcal{H}_{n}^{f,q}}\sup_{f \in W^{s,p}(\Omega)}\mathlarger{\mathlarger{\mathbb{E}}}_{\substack{\{x_i\}_{i=1}^{n},\{y_i\}_{i=1}^{n}}}\Bigg[\left|\hat{H}^q\Big(\{x_i\}_{i=1}^{n},\{y_i\}_{i=1}^{n}\Big) - I_f^q\right|\Bigg]\\
&\geq \Delta \inf_{\hat{H}^q \in \mathcal{H}_{n}^{f,q}}\sup_{f \in W^{s,p}(\Omega)}\mathlarger{\mathlarger{\mathbb{P}}}_{\substack{\{x_i\}_{i=1}^{n},\{y_i\}_{i=1}^{n}}}\Bigg[\left|\hat{H}^q\Big(\{x_i\}_{i=1}^{n},\{y_i\}_{i=1}^{n}\Big) - I_f^q\right| \geq \Delta\Bigg]\\
&\geq \frac{1}{2}I_{g_1}^q\frac{1-TV(\mathbb{P}_0 || \mathbb{P}_1)-\beta_0-\beta_1}{2} \geq \frac{1}{4}\Big(1-\frac{\sqrt{3}}{10}\Big)\epsilon I_{g_1}^q \\
&= \frac{1}{8}\Big(1-\frac{\sqrt{3}}{10}\Big)(200n)^{-\frac{q}{d}(s-\frac{d}{p})-1}||K||_{L^q([-\frac{1}{2},\frac{1}{2}]^d)}^q \gtrsim n^{-q(\frac{s}{d}-\frac{1}{p})-1}, 
\end{aligned}    
\end{equation}
which is exactly the first term in the RHS of (\ref{eq: lower bound q-th moment}).

(Case II) Now let us proceed to prove the second lower bound in (\ref{eq: lower bound q-th moment}). For any $1 \leq j \leq m^d$, consider first some function $f_j$ defined as follows
\begin{equation}
\begin{aligned}
f_j(x) &= 
\begin{cases}
m^{-s}K(m(x-c_j)) \ (x \in \Omega_j), \\
0 \ (\text{otherwise}),
\end{cases}    
\end{aligned}
\end{equation}
which satisfies $\text{supp}(f_j) \subseteq \Omega_j$, $f_j \in C^{\infty}(\Omega)$ and $f_j(x) \geq 0 \ (\forall \ x \in \Omega)$. We further pick two constants $\alpha, M$ satisfying $\alpha := ||K||_{L^{\infty}([-\frac{1}{2},\frac{1}{2}]^d)}$ and $M= 3\alpha$. Now consider the following finite set of $2^{m^d}$ functions:
\begin{equation}
\mathcal{S} := \Big\{M+\sum_{j=1}^{m^d}\eta_j f_j: \eta_j \in \{\pm 1\}, \ \forall \ 1 \leq j \leq m^d\Big\}.    
\end{equation}
We will proceed to verify that any element in $\mathcal{S}$ must be in $W^{s,p}(\Omega)$ for any $m$. Note that for any $\eta_j \in \{\pm 1\} \ (1 \leq j \leq m^d)$ and  any $t \in \mathbb{N}_{0}^d$ satisfying $|t| \leq s$, we have
\begin{align*}
\Bigg||D^{t}\Big(M+\sum_{j=1}^{m^d}\eta_j f_j\Big)\Bigg||_{L^p(\Omega)}^p &\leq \Bigg(M + \Big||\sum_{j=1}^{m^d}\eta_j (D^{t} f_j)\Big||_{L^p(\Omega)}\Bigg)^p \\
&\leq 2^p\Bigg(M^p+\Big||\sum_{j=1}^{m^d}\eta_j (D^{t} f_j)\Big||_{L^p(\Omega)}^p\Bigg)\lesssim M^p + \sum_{j=1}^{m^d}||D^t f_j||_{L^p(\Omega_j)}^p \\
&= M^p + \sum_{j=1}^{m^d}\int_{\Omega_j}\Big|m^{-s+|t|}(D^{t}K)(m(x-c_j))\Big|^p dx\\
&= M^p + m^{(|t|-s)p}\sum_{j=1}^{m^d}\int_{[-\frac{1}{2},\frac{1}{2}]^d}\Big|(D^{t}K)(y)\Big|^p \frac{1}{m^d}dy \\
&\leq M^p + ||D^tK||_{L^p([-\frac{1}{2},\frac{1}{2}]^d)}^p  \lesssim 1.   
\end{align*}
This gives us that $\mathcal{S} \subset W^{s,p}(\Omega)$ for any $m$, as desired. Now let's pick $\kappa = \frac{1}{3}\sqrt{\frac{2}{3n}}$ and take $\{w_{j}^{(0)}\}_{j=1}^{m^d}$ and $\{w_{j}^{(1)}\}_{j=1}^{m^d}$ to be independent and identical copies of $w_{\frac{1+\kappa}{2}}$ and $w_{\frac{1-\kappa}{2}}$ respectively. Then we define $\mu_0,\mu_1$ to be two discrete measures supported on the finite set $\mathcal{S}$ such that the following condition holds for any $\eta_j \in \{\pm 1\} \ (1 \leq j \leq m^d)$:
\begin{equation}
\label{eqn: def of two priors_1 in thm 2.1 case 2}
\mu_k\Big(\Big\{M+\sum_{j=1}^{m^d}\eta_j f_j\Big\}\Big) = \prod_{j=1}^{m^d}\mathbb{P}(w_{j}^{(k)} = \eta_j), \ k \in \{0,1\}.
\end{equation}

In order to determine the separation distance $\Delta$ between the two priors $\mu_0$ and $\mu_1$, we need to define two quantities $A:= \int_{\Omega_j}(M+f_j(x))^qdx $ and $B:= \int_{\Omega_j}(M-f_j(x))^qdx$, which both remain the same for any $1 \leq j \leq m^d$. Now consider deriving a lower bound on the quantity $\Delta':=A-B > 0$. Note that for any fixed $j \in \{1,2,\cdots,m^d\}$, we have $M > 2\alpha \geq 2m^{-s}||K||_{L^{\infty}([-\frac{1}{2},\frac{1}{2}]^d)} = 2||f_j||_{L^{\infty}(\Omega_j)}$, which implies $M+y > \frac{1}{2}M > 0$ for any $y \in [-||f_j||_{L^{\infty}(\Omega_j)},||f_j||_{L^{\infty}(\Omega_j)}]$. This helps us obtain the following lower bound on $\Delta'$:
\begin{equation}
\label{eqn: thm 2.1 case 2 lower bound on grid sep distance}
\begin{aligned}
\Delta' &= \int_{\Omega_j}(M+f_j(x))^qdx - \int_{\Omega_j}(M-f_j(x))^qdx = \int_{\Omega_j}\Big(\int_{-f_j(x)}^{f_j(x)}q(M+y)^{q-1}dy\Big)dx\\
&\geq \int_{\Omega_j}\Big(\int_{-f_j(x)}^{f_j(x)}q(\frac{1}{2}M)^{q-1}dy\Big)dx = \frac{q}{2^{q-1}}M^{q-1}\int_{\Omega_j}\Big(2f_j(x)\Big)dx \gtrsim \int_{\Omega_j}f_j(x)dx\\
&= \int_{\Omega_j}m^{-s}K(m(x-c_j))dx = m^{-s}\int_{[-\frac{1}{2},\frac{1}{2}]^d}K(y)\frac{1}{m^d}dy = m^{-s-d}||K||_{L^1([-\frac{1}{2},\frac{1}{2}]^d)}. 
\end{aligned}    
\end{equation}
Moreover, let us pick $\lambda = \frac{1}{2}$ and apply Hoeffding's Inequality (Lemma \ref{lem: hoeffding ineq}) to the bounded random variables $\{w_{j}^{(0)}\}_{j=1}^{m^d}$ and $\{w_{j}^{(1)}\}_{j=1}^{m^d}$ to deduce that
\begin{equation}
\label{eqn: thm 2.1 case 2 hoeffding bound}
\begin{aligned}
\mathbb{P}\Big(\sum_{j=1}^{m^d}w_{j}^{(0)} \geq -(1-\lambda)m^d \kappa\Big) &\leq \exp\Big(-\frac{2(\lambda m^d \kappa)^2}{4m^d}\Big) = \exp\Big(-\frac{1}{2}\lambda^2 \kappa^2 m^d\Big), \\
\mathbb{P}\Big(\sum_{j=1}^{m^d}w_{j}^{(1)} \leq (1-\lambda)m^d \kappa\Big) &\leq \exp\Big(-\frac{2(\lambda m^d \kappa)^2}{4m^d}\Big) = \exp\Big(-\frac{1}{2}\lambda^2 \kappa^2 m^d\Big).
\end{aligned}    
\end{equation}
By taking $c:= \frac{m^d}{2}(A+B), \Delta:=(1-\lambda)\kappa m^d(A-B)=(1-\lambda)\kappa m^d\Delta'$ and $\beta_0 =\beta_1 =\exp\Big(-\frac{1}{2}\lambda^2 \kappa^2 m^d\Big)$, we may combine (\ref{eqn: thm 2.1 case 2 lower bound on grid sep distance}) and (\ref{eqn: thm 2.1 case 2 hoeffding bound}) justified above to get that
\begin{equation}
\begin{aligned}
&\mu_0(f \in W^{s,p}(\Omega): I_{f}^q \leq c-\Delta) = \mathbb{P}\Big(\sum_{j=1}^{m^d}I_{M+w_{j}^{(0)}f_j}^{q} \leq \frac{1-(1-\lambda)\kappa}{2}m^d A + \frac{1+(1-\lambda)\kappa}{2}m^d B\Big) \\
&\geq \mathbb{P}\Big(\sum_{j=1}^{m^d}w_{j}^{(0)} \leq -(1-\lambda)m^d \kappa\Big)=1-\mathbb{P}\Big(\sum_{j=1}^{m^d}w_{j}^{(0)} \geq -(1-\lambda)m^d \kappa\Big)\\
&\geq 1-\exp\Big(-\frac{1}{2}\lambda^2 \kappa^2 m^d\Big) = 1-\beta_0, \\
&\mu_1(f \in W^{s,p}(\Omega): I_{f}^q \geq c+\Delta) = \mathbb{P}\Big(\sum_{j=1}^{m^d}I_{M+w_{j}^{(1)}f_j}^{q} \geq \frac{1+(1-\lambda)\kappa}{2}m^d A + \frac{1-(1-\lambda)\kappa}{2}m^d B\Big) \\
&\geq \mathbb{P}\Big(\sum_{j=1}^{m^d}w_{j}^{(1)} \geq (1-\lambda)m^d \kappa\Big)=1-\mathbb{P}\Big(\sum_{j=1}^{m^d}w_{j}^{(0)} \leq (1-\lambda)m^d \kappa\Big)\\
&\geq 1-\exp\Big(-\frac{1}{2}\lambda^2 \kappa^2 m^d\Big) = 1-\beta_1,   
\end{aligned}    
\end{equation}
which indicates that (\ref{eq: separated hypothesis (two prior)}) holds true. Now let's consider bounding the KL divergence between the two marginal distributions $\mathbb{P}_0,\mathbb{P}_1$ associated with $\mu_0,\mu_1$, respectively. Using the fact that $\{x_1,\cdots,x_n\}$ are identical and independent samples from the uniform distribution on $\Omega$ again allows us to write the marginal distributions in an explicit form as follows: 
\begin{equation}
\label{eqn: marginal distributions thm 2.1 case 2}
\begin{aligned}
\mathbb{P}_0(\vec{x},\vec{y}) &= \prod_{j=1}^{m^d}\Big(\frac{1+\kappa}{2}\prod_{i:x_i \in \Omega_j}\delta_{M-f_j(x_i)}(y_i) + \frac{1-\kappa}{2}\prod_{i:x_i \in \Omega_j}\delta_{M+f_j(x_i)}(y_i)\Big),\\
\mathbb{P}_1(\vec{x},\vec{y}) &= \prod_{j=1}^{m^d}\Big(\frac{1-\kappa}{2}\prod_{i:x_i \in \Omega_j}\delta_{M-f_j(x_i)}(y_i) + \frac{1+\kappa}{2}\prod_{i:x_i \in \Omega_j}\delta_{M+f_j(x_i)}(y_i)\Big). 
\end{aligned}    
\end{equation}
Furthermore, for any $n$ quadrature points $\{x_i\}_{i=1}^{n}$, we use $\mathcal{J}_n$ to denote the set of all indices $j$ satisfying that $\Omega_j$ contains at least one of the points in $\{x_i\}_{i=1}^{n}$, \emph{i.e,}
\begin{equation}
\mathcal{J}_{n}:=\mathcal{J}_n(x_1,\cdots,x_n) = \Big\{j: 1 \leq j \leq m^d \text{ and } \Omega_j \cap \{x_1,\cdots,x_n\} \neq \varnothing\Big\} 
\end{equation}
Given that $m^d=200n > n$, we have $|\mathcal{J}_n| \leq n$ for any $n$ quadrature points $\{x_i\}_{i=1}^{n}$. Using this upper bound on $|\mathcal{J}_n|$ allows us to bound the KL divergence between $\mathbb{P}_0$ and $\mathbb{P}_1$ in the following way:
\begin{equation}
\label{eqn: thm 2.1 case 2 upper bound on KL}
\begin{aligned}
&KL(\mathbb{P}_{0} || \mathbb{P}_{1}) = \int_{\Omega} \cdots\int_{\Omega}\Big(\int_{-\infty}^{\infty}\cdots\int_{-\infty}^{\infty}\log\Big(\frac{\mathbb{P}_0(\vec{x},\vec{y})}{\mathbb{P}_1(\vec{x},\vec{y})}\Big)\mathbb{P}_0(\vec{x},\vec{y})dy_1 \cdots dy_n\Big)dx_1 \cdots dx_n\\
&= \int_{\Omega} \cdots\int_{\Omega}\Bigg(\int_{-\infty}^{\infty}\cdots\int_{-\infty}^{\infty}\log\Big(\prod_{j=1}^{m^d}\frac{\frac{1+\kappa}{2}\prod_{i:x_i \in \Omega_j}\delta_{M-f_j(x_i)}(y_i) + \frac{1-\kappa}{2}\prod_{i:x_i \in \Omega_j}\delta_{M+f_j(x_i)}(y_i)}{\frac{1-\kappa}{2}\prod_{i:x_i \in \Omega_j}\delta_{M-f_j(x_i)}(y_i) + \frac{1+\kappa}{2}\prod_{i:x_i \in \Omega_j}\delta_{M+f_j(x_i)}(y_i)}\Big)\\
&\cdot \prod_{j=1}^{m^d}\Big(\frac{1+\kappa}{2}\prod_{i:x_i \in \Omega_j}\delta_{M-f_j(x_i)}(y_i) + \frac{1-\kappa}{2}\prod_{i:x_i \in \Omega_j}\delta_{M+f_j(x_i)}(y_i)\Big)\prod_{i=1}^{n}dy_i\Bigg)\prod_{i=1}^{n}dx_i\\
&= \int_{\Omega} \cdots\int_{\Omega}\Bigg(\int_{-\infty}^{\infty}\cdots\int_{-\infty}^{\infty}\log\Big(\prod_{j \in \mathcal{J}_n}\frac{\frac{1+\kappa}{2}\prod_{i:x_i \in \Omega_j}\delta_{M-f_j(x_i)}(y_i) + \frac{1-\kappa}{2}\prod_{i:x_i \in \Omega_j}\delta_{M+f_j(x_i)}(y_i)}{\frac{1-\kappa}{2}\prod_{i:x_i \in \Omega_j}\delta_{M-f_j(x_i)}(y_i) + \frac{1+\kappa}{2}\prod_{i:x_i \in \Omega_j}\delta_{M+f_j(x_i)}(y_i)}\Big)\\
&\cdot \prod_{j \in \mathcal{J}_n}\Big(\frac{1+\kappa}{2}\prod_{i:x_i \in \Omega_j}\delta_{M-f_j(x_i)}(y_i) + \frac{1-\kappa}{2}\prod_{i:x_i \in \Omega_j}\delta_{M+f_j(x_i)}(y_i)\Big)\prod_{i=1}^{n}dy_i\Bigg)\prod_{i=1}^{n}dx_i\\
&= \int_{\Omega} \cdots\int_{\Omega}\Bigg(\sum_{j \in \mathcal{J}_n}\int_{-\infty}^{\infty}\cdots\int_{-\infty}^{\infty}\log\Big(\frac{\frac{1+\kappa}{2}\prod_{i:x_i \in \Omega_j}\delta_{M-f_j(x_i)}(y_i) + \frac{1-\kappa}{2}\prod_{i:x_i \in \Omega_j}\delta_{M+f_j(x_i)}(y_i)}{\frac{1-\kappa}{2}\prod_{i:x_i \in \Omega_j}\delta_{M-f_j(x_i)}(y_i) + \frac{1+\kappa}{2}\prod_{i:x_i \in \Omega_j}\delta_{M+f_j(x_i)}(y_i)}\Big)\\
&\cdot \Big(\frac{1+\kappa}{2}\prod_{i:x_i \in \Omega_j}\delta_{M-f_j(x_i)}(y_i) + \frac{1-\kappa}{2}\prod_{i:x_i \in \Omega_j}\delta_{M+f_j(x_i)}(y_i)\Big)\prod_{i:x_i \in \Omega_j}dy_i\Bigg)\prod_{i=1}^{n}dx_i\\
&=\int_{\Omega} \cdots\int_{\Omega}|\mathcal{J}_n|\Big(\log\Big(\frac{1+\kappa}{1-\kappa}\Big)\frac{1+\kappa}{2} + \log\Big(\frac{1-\kappa}{1+\kappa}\Big)\frac{1-\kappa}{2}\Big)\prod_{i=1}^{n}dx_i \leq n\kappa\log\Big(\frac{1+\kappa}{1-\kappa}\Big).
\end{aligned}    
\end{equation}

Now we may combine (\ref{eqn: thm 2.1 case 2 upper bound on KL}) and Pinkser's inequality to upper bound the TV distance between $\mathbb{P}_0$ and $\mathbb{P}_1$ as below:
\begin{equation}
\label{thm 2.1 case 2: bound on TV}
TV(\mathbb{P}_0 || \mathbb{P}_1) \leq  \sqrt{\frac{1}{2}KL(\mathbb{P}_0 || \mathbb{P}_1)} \leq \sqrt{\frac{n\kappa}{2}\log\Big(\frac{1+\kappa}{1-\kappa}\Big)} \leq \sqrt{\frac{3n}{2}}\kappa= \frac{1}{3}.   
\end{equation}

Finally, by substituting (\ref{eqn: thm 2.1 case 2 lower bound on grid sep distance}), (\ref{thm 2.1 case 2: bound on TV}), $\Delta =(1-\lambda)\kappa m^d\Delta'$ and $\beta_0=\beta_1=\exp\Big(-\frac{1}{2}\lambda^2 \kappa^2 m^d\Big) = \exp(-\frac{50}{27}) < \frac{1}{6}$ into (\ref{eq: bound on misclassification probability (two prior)}) and applying Markov's inequality, we obtain the final lower bound
\begin{equation}
\label{eqn: thm 2.1 case 2 final bound}
\begin{aligned}
&\inf_{\hat{H}^q \in \mathcal{H}_{n}^{f,q}}\sup_{f \in W^{s,p}(\Omega)}\mathlarger{\mathlarger{\mathbb{E}}}_{\substack{\{x_i\}_{i=1}^{n},\{y_i\}_{i=1}^{n}}}\Bigg[\left|\hat{H}^q\Big(\{x_i\}_{i=1}^{n},\{y_i\}_{i=1}^{n}\Big) - I_f^q\right|\Bigg]\\
&\geq \Delta \inf_{\hat{H}^q \in \mathcal{H}_{n}^{f,q}}\sup_{f \in W^{s,p}(\Omega)}\mathlarger{\mathlarger{\mathbb{P}}}_{\substack{\{x_i\}_{i=1}^{n},\{y_i\}_{i=1}^{n}}}\Bigg[\left|\hat{H}^q\Big(\{x_i\}_{i=1}^{n},\{y_i\}_{i=1}^{n}\Big) - I_f^q\right| \geq \Delta\Bigg]\\
&\geq (1-\lambda)\kappa m^d\Delta'\frac{1-TV(\mathbb{P}_0 || \mathbb{P}_1)-\beta_0-\beta_1}{2} \geq \frac{1}{2}\frac{\sqrt{2}}{3\sqrt{3n}}\cdot(200n) \cdot \frac{\Delta'}{6} \\
&\gtrsim \sqrt{n}\Delta' \gtrsim \sqrt{n}(200n)^{-\frac{s+d}{d}}||K||_{L^1([-\frac{1}{2},\frac{1}{2}]^d)} \gtrsim n^{-\frac{s}{d}-\frac{1}{2}}, 
\end{aligned}    
\end{equation}
which is exactly the second term in the RHS of (\ref{eq: lower bound q-th moment}). Combining the two lower bounds proved in (\ref{eqn: thm 2.1 case 1 final bound}) and (\ref{eqn: thm 2.1 case 2 final bound}) concludes our proof of Theorem \ref{thm: lower bound q-th moment}.

\section{Proof of Upper Bounds in Section \ref{sec 3: upper bound on q-moment estimation}}
\label{app: proof of upper bound in sec 3}
\subsection{Proof of Theorem \ref{thm: upper bound control variate} (Regression-Adjusted Control Variate)}
\label{app: proof of upper bound in sec 3.1}
In this subsection, we present a detailed proof of Theorem \ref{thm: upper bound control variate}. With the first half of the quadrature points $\{x_i\}_{i=1}^{\frac{n}{2}}$ and observed function values $\{y_i\}_{i=1}^{\frac{n}{2}}$ as inputs, we pick the regression adjusted control variate $\hat{f}_{1:\frac{n}{2}}$ to be the estimator returned by the oracle $K_{\frac{n}{2}}$ specified in Assumption \ref{assumption: existence of oracle}. Moreover, we use the following expression to denote the variance of the function $\hat{f}_{1:\frac{n}{2}}^q(x) - f^q(x)$ with respect to the uniform distribution on $\Omega$:
\begin{equation}
\label{eqn def: variance of difference in f^q}
\text{Var}(\hat{f}_{1:\frac{n}{2}}^q-f^q) := \int_{\Omega}(f^q(x)-\hat{f}_{1:\frac{n}{2}}^q(x))^2dx - \Big(\int_{\Omega}(f^q(x)-\hat{f}_{1:\frac{n}{2}}^q(x))dx\Big)^2.     
\end{equation}
By plugging in the expression of $\hat{H}^q_{C}, I_{f}^q$ and using the fact that $\{x_i\}_{i=1}^{n}$ are identical and independent copies of the uniform random variable over $\Omega$, we have 
\begin{equation}
\label{eqn: variance decomposition C2}
\begin{aligned}
&\mathlarger{\mathlarger{\mathbb{E}}}_{\substack{\{x_i\}_{i=1}^{n},\{y_i\}_{i=1}^{n}}}\Bigg[\left|\hat{H}_{C}^q\Big(\{x_i\}_{i=1}^{n},\{y_i\}_{i=1}^{n}\Big) - I_f^q\right|^2\Bigg] \\
&= \mathlarger{\mathlarger{\mathbb{E}}}_{\substack{\{x_i\}_{i=1}^{n}}}\Bigg[\Big|\int_{\Omega}\hat{f}_{1:\frac{n}{2}}^q(x)dx + \frac{2}{n}\sum_{i=\frac{n}{2}+1}^{n}\Big(f^q(x_i)-\hat{f}_{1:\frac{n}{2}}^q(x_i)\Big)-\int_{\Omega}f^q(x) dx\Big|^2\Bigg]\\
&= \mathlarger{\mathlarger{\mathbb{E}}}_{\substack{\{x_i\}_{i=1}^{\frac{n}{2}}}}\Bigg[\mathlarger{\mathlarger{\mathbb{E}}}_{\substack{\{x_i\}_{i=\frac{n}{2}+1}^{n}}}\Bigg[\Big|\frac{1}{\frac{n}{2}}\sum_{i=\frac{n}{2}+1}^{n}\Big(f^q(x_i)-\hat{f}_{1:\frac{n}{2}}^q(x_i)-\int_{\Omega}(f^q(x)-\hat{f}_{1:\frac{n}{2}}^q(x))dx\Big)\Big|^2\Bigg]\Bigg]\\
&= \mathlarger{\mathlarger{\mathbb{E}}}_{\substack{\{x_i\}_{i=1}^{\frac{n}{2}}}}\Bigg[\frac{4}{n^2}\sum_{i=\frac{n}{2}+1}^{n}\mathlarger{\mathlarger{\mathbb{E}}}_{x_i}\Bigg[\Big|\Big(f^q(x_i)-\hat{f}_{1:\frac{n}{2}}^q(x_i)-\int_{\Omega}(f^q(x)-\hat{f}_{1:\frac{n}{2}}^q(x))dx\Big)\Big|^2\Bigg]\Bigg]\\
&= \mathlarger{\mathlarger{\mathbb{E}}}_{\substack{\{x_i\}_{i=1}^{\frac{n}{2}}}}\Big[\frac{4}{n^2}\sum_{i=\frac{n}{2}+1}^{n}\text{Var}(\hat{f}_{1:\frac{n}{2}}^q-f^q)\Big] = \frac{2}{n}\mathlarger{\mathlarger{\mathbb{E}}}_{\substack{\{x_i\}_{i=1}^{\frac{n}{2}}}}\Big[\text{Var}(\hat{f}_{1:\frac{n}{2}}^q-f^q)\Big].    
\end{aligned}    
\end{equation}
From the identity above, we know that it suffices to upper bound the term $\mathlarger{\mathlarger{\mathbb{E}}}_{\substack{\{x_i\}_{i=1}^{\frac{n}{2}}}}\Big[\text{Var}(\hat{f}_{1:\frac{n}{2}}^q-f^q)\Big]$. Let $g_{1:\frac{n}{2}}:=\hat{f}_{1:\frac{n}{2}}-f$ denote the difference between the estimator $\hat{f}_{1:\frac{n}{2}}$ and underlying function $f$. Then we may further upper bound the expression $\mathlarger{\mathlarger{\mathbb{E}}}_{\substack{\{x_i\}_{i=1}^{\frac{n}{2}}}}\Big[\text{Var}(\hat{f}_{1:\frac{n}{2}}^q-f^q)\Big]$ as follows:
\begin{equation}
\label{eqn: upper bound on the variance by two integrals}
\begin{aligned}
&\mathlarger{\mathlarger{\mathbb{E}}}_{\substack{\{x_i\}_{i=1}^{\frac{n}{2}}}}\Big[\text{Var}(\hat{f}_{1:\frac{n}{2}}^q-f^q)\Big] = \mathlarger{\mathlarger{\mathbb{E}}}_{\substack{\{x_i\}_{i=1}^{\frac{n}{2}}}}\Big[\int_{\Omega}(f^q(x)-\hat{f}_{1:\frac{n}{2}}^q(x))^2dx - \Big(\int_{\Omega}(f^q(x)-\hat{f}_{1:\frac{n}{2}}^q(x))dx\Big)^2\Big]\\
&\leq \mathlarger{\mathlarger{\mathbb{E}}}_{\substack{\{x_i\}_{i=1}^{\frac{n}{2}}}}\Big[\int_{\Omega}\Big(f^q(x)-\hat{f}_{1:\frac{n}{2}}^q(x)\Big)^2dx\Big] = \mathlarger{\mathlarger{\mathbb{E}}}_{\substack{\{x_i\}_{i=1}^{\frac{n}{2}}}}\Bigg[\int_{\Omega}\Big((f(x) + g_{1:\frac{n}{2}}(x))^q-f^q(x)\Big)^2dx\Bigg]\\
&= \mathlarger{\mathlarger{\mathbb{E}}}_{\substack{\{x_i\}_{i=1}^{\frac{n}{2}}}}\Bigg[\int_{\Omega}\Big(\int_{0}^{g_{1:\frac{n}{2}}(x)}q(f(x)+y)^{q-1}dy\Big)^2dx\Bigg] \\
&\leq \mathlarger{\mathlarger{\mathbb{E}}}_{\substack{\{x_i\}_{i=1}^{\frac{n}{2}}}}\Bigg[\int_{\Omega}\Big|\int_{0}^{g_{1:\frac{n}{2}}(x)}1dy\Big|\Big|\int_{0}^{g_{1:\frac{n}{2}}(x)}q^2(|f(x)+y|^2)^{q-1}dy\Big|dx\Bigg]\\
&\lesssim \mathlarger{\mathlarger{\mathbb{E}}}_{\substack{\{x_i\}_{i=1}^{\frac{n}{2}}}}\Bigg[\int_{\Omega}|g_{1:\frac{n}{2}}(x)| \cdot |g_{1:\frac{n}{2}}(x)|\max\Big\{|f^{2q-2}(x)|,|g^{2q-2}_{1:\frac{n}{2}}(x)|\Big\}dx\Bigg] \\
&\lesssim \mathlarger{\mathlarger{\mathbb{E}}}_{\substack{\{x_i\}_{i=1}^{\frac{n}{2}}}}\Bigg[\int_{\Omega}|g_{1:\frac{n}{2}}^{2q}(x)|dx\Bigg] + \mathlarger{\mathlarger{\mathbb{E}}}_{\substack{\{x_i\}_{i=1}^{\frac{n}{2}}}}\Bigg[\int_{\Omega}|g_{1:\frac{n}{2}}^2(x)f^{2q-2}(x)|dx\Bigg]. 
\end{aligned}    
\end{equation}
Now let's proceed to bound from above the two expected integrals in the last line of (\ref{eqn: upper bound on the variance by two integrals}). For the first expected integral, since $s > \frac{2dq-dp}{2pq} \Rightarrow \frac{1}{2q} > \frac{d-sp}{pd}$, we may apply (\ref{oracle: assumption}) in Assumption \ref{assumption: existence of oracle} to deduce that
\begin{equation}
\label{eqn: App C bound part 1}
\begin{aligned}
\mathlarger{\mathlarger{\mathbb{E}}}_{\substack{\{x_i\}_{i=1}^{\frac{n}{2}}}}\Bigg[\int_{\Omega}|g_{1:\frac{n}{2}}^{2q}(x)|dx\Bigg] &= \mathlarger{\mathbb{E}}_{\{x_{i}\}_{i=1}^{\frac{n}{2}}}\Big[||\hat{f}_{1:\frac{n}{2}}-f||_{L^{2q}(\Omega)}^{2q}\Big] \\
&\lesssim ((\frac{n}{2})^{-\frac{s}{d}+(\frac{1}{p}-\frac{1}{2q})_{+}})^{2q} \lesssim n^{2q(-\frac{s}{d}+\frac{1}{p}-\frac{1}{2q})}=n^{2q(\frac{1}{p}-\frac{s}{d})-1},    
\end{aligned}  
\end{equation}
where the last equality above follows from the given assumption that $p<2q$. Now let's proceed to bound from above the second expected integral in (\ref{eqn: upper bound on the variance by two integrals}). Here we define $p^\ast = (\max\{\frac{1}{p}-\frac{s}{d},0\})^{-1}$, \emph{i.e,} $p^\ast = \frac{pd}{d-sp}$ when $s<\frac{d}{p}$ and $p^\ast=\infty$ otherwise. From Sobolev Embedding Theorem (Lemma \ref{lem: sobolev embedding}), we have that $W^{s,p}(\Omega) \subseteq L^{p^\ast}(\Omega)$. Based on the value of the smoothness parameter $s$, we have three separate cases as below:

(Case I) When $s \in (\frac{d}{p},\infty)$, we have $p^\ast = \infty$ and $f \in W^{s,p}(\Omega) \subset L^{\infty}(\Omega)$. Since $\hat{f}_{1:\frac{n}{2}}$ and $f$ are both in the Sobolev space $W^{s,p}(\Omega) \subseteq L^{\infty}(\Omega)$, we may further deduce that $g_{1:\frac{n}{2}} = \hat{f}_{1:\frac{n}{2}} - f \in W^{s,p}(\Omega) \subseteq L^{\infty}(\Omega) \subseteq L^{2}(\Omega)$. By picking $r=2$ in in (\ref{oracle: assumption}) of Assumption \ref{assumption: existence of oracle}, we may use the facts that $p > 2$ and $f \in L^{\infty}(\Omega)$ to deduce that
\begin{equation}
\label{eqn: App C final upper bound on C3_2, Case 0}
\begin{aligned}
&\mathlarger{\mathlarger{\mathbb{E}}}_{\substack{\{x_i\}_{i=1}^{\frac{n}{2}}}}\Bigg[\int_{\Omega}|g_{1:\frac{n}{2}}^2(x)f^{2q-2}(x)|dx\Bigg] \lesssim \mathlarger{\mathlarger{\mathbb{E}}}_{\substack{\{x_i\}_{i=1}^{\frac{n}{2}}}}\Bigg[\int_{\Omega}|g_{1:\frac{n}{2}}^2(x)|dx\Bigg] \\
&= \mathlarger{\mathbb{E}}_{\{x_{i}\}_{i=1}^{\frac{n}{2}}}\Big[||\hat{f}_{1:\frac{n}{2}}-f||_{L^{2}(\Omega)}^{2}\Big]
\lesssim \Big(n^{-\frac{s}{d} + (\frac{1}{p}-\frac{1}{2})_{+}}\Big)^2 = n^{-\frac{2s}{d}},    
\end{aligned}   
\end{equation}
which is our final upper bound on the second expected integral in (\ref{eqn: upper bound on the variance by two integrals}) under the assumption that $s \in (\frac{d}{p},\infty)$. 

(Case II) When $s \in (\frac{d(2q-p)}{p(2q-2)}, \frac{d}{p})$, we have $p^\ast = \frac{pd}{d-sp} > \frac{pd}{d-p\frac{d(2q-p)}{p(2q-2)}}=\frac{p(2q-2)}{p-2}$, which implies $f \in W^{s,p}(\Omega) \subseteq L^{p^\ast}(\Omega) \subseteq L^{\frac{p(2q-2)}{p-2}}(\Omega) \subseteq L^p(\Omega)$. Given that $\frac{p}{p-2} > 1$, we can further deduce that $f^{2q-2} \in L^{\frac{p}{p-2}}(\Omega)$ . Moreover, since $\hat{f}_{1:\frac{n}{2}} \in W^{s,p}(\Omega) \subseteq L^p(\Omega)$, we have that $g_{1:\frac{n}{2}} = \hat{f}_{1:\frac{n}{2}}-f \in L^p(\Omega)$. Given that $p>2$, we can further deduce that $g_{1:\frac{n}{2}}^2 \in L^{\frac{p}{2}}(\Omega)$. Then we may apply H\"older's inequality (Lemma \ref{lem: holder ineq}) to $g_{1:\frac{n}{2}}^2 \in L^{\frac{p}{2}}(\Omega)$ and $f^{2q-2} \in L^{\frac{p}{p-2}}(\Omega)$ to obtain that 
\begin{equation}
\label{eqn: App C equation C5}
\begin{aligned}
&\mathlarger{\mathlarger{\mathbb{E}}}_{\substack{\{x_i\}_{i=1}^{\frac{n}{2}}}}\Bigg[\int_{\Omega}|g_{1:\frac{n}{2}}^2(x)f^{2q-2}(x)|dx\Bigg] = \mathlarger{\mathlarger{\mathbb{E}}}_{\substack{\{x_i\}_{i=1}^{\frac{n}{2}}}}\Big[||g_{1:\frac{n}{2}}^2f^{2q-2}||_{L^1(\Omega)}\Big] \\
&\leq \mathlarger{\mathlarger{\mathbb{E}}}_{\substack{\{x_i\}_{i=1}^{\frac{n}{2}}}}\Big[\Big||g_{1:\frac{n}{2}}^2\Big||_{L^{\frac{p}{2}}(\Omega)}\Big||f^{2q-2}\Big||_{L^{\frac{p}{p-2}}(\Omega)}\Big] \leq \Big||f\Big||_{L^{\frac{p(2q-2)}{p-2}}(\Omega)}^{2q-2}\mathlarger{\mathlarger{\mathbb{E}}}_{\substack{\{x_i\}_{i=1}^{\frac{n}{2}}}}\Big[\Big||g_{1:\frac{n}{2}}\Big||_{L^{p}(\Omega)}^2\Big].     
\end{aligned}    
\end{equation}
Note that the function $h(t) = t^{\frac{2}{p}}$ is concave and $\frac{1}{p} \in (\frac{d-sp}{pd},1]$ when $p > 2$. Hence, applying Jensen's inequality and picking $r=p$ in (\ref{oracle: assumption}) of Assumption \ref{assumption: existence of oracle} further allows us to upper bound the last term in (\ref{eqn: App C equation C5}) as follows:
\begin{equation}
\label{eqn: App C equation C6}
\begin{aligned}
\mathlarger{\mathbb{E}}_{\substack{\{x_i\}_{i=1}^{\frac{n}{2}}}}\Big[||g_{1:\frac{n}{2}}||_{L^{p}(\Omega)}^2\Big] &= \mathlarger{\mathbb{E}}_{\substack{\{x_i\}_{i=1}^{\frac{n}{2}}}}\Big[\Big(||g_{1:\frac{n}{2}}||_{L^{p}(\Omega)}^{p}\Big)^{\frac{2}{p}}\Big] \\
&\leq \mathlarger{\mathbb{E}}_{\substack{\{x_i\}_{i=1}^{\frac{n}{2}}}}\Big[||g_{1:\frac{n}{2}}||_{L^{p}(\Omega)}^p\Big]^{\frac{2}{p}} = \mathlarger{\mathbb{E}}_{\{x_{i}\}_{i=1}^{\frac{n}{2}}}\Big[||\hat{f}_{1:\frac{n}{2}}-f||_{L^{p}(\Omega)}^{p}\Big]^{\frac{2}{p}} \\
&\lesssim \Big((\frac{n}{2})^{-\frac{s}{d}+(\frac{1}{p}-\frac{1}{p})_{+}}\Big)^{2} \lesssim n^{-\frac{2s}{d}}.
\end{aligned}
\end{equation}
Substituting (\ref{eqn: App C equation C6}) into (\ref{eqn: App C equation C5}) then gives us the final upper bound on the second expected integral in (\ref{eqn: upper bound on the variance by two integrals}) under the assumption that $s \in (\frac{d(2q-p)}{p(2q-2)}, \infty)$:
\begin{equation}
\label{eqn: App C final upper bound on C3_2, Case 1}
\mathlarger{\mathlarger{\mathbb{E}}}_{\substack{\{x_i\}_{i=1}^{\frac{n}{2}}}}\Bigg[\int_{\Omega}|g_{1:\frac{n}{2}}^2(x)f^{2q-2}(x)|dx\Bigg] \lesssim n^{-\frac{2s}{d}}.    
\end{equation}

(Case III) When $s \in (\frac{d(2q-p)}{2pq},\frac{d(2q-p)}{p(2q-2)})$, we have that $s < \frac{d}{p}$, which indicates that $p^\ast = \frac{pd}{d-sp}$ satisfies $2q < p^\ast < \frac{p(2q-2)}{p-2}$. Given that $p^\ast > 2q>2q-2$ and $f \in W^{s,p}(\Omega) \subseteq L^{p^\ast}(\Omega)$, we can deduce that $f^{2q-2} \in L^{\frac{p^\ast}{2q-2}}(\Omega)$. Furthermore, note that $p^\ast > 2q$ implies $\frac{2p^\ast}{p^\ast + 2-2q} < p^\ast$ and $p^\ast < \frac{p(2q-2)}{p-2}$ implies $\frac{2p^\ast}{p^\ast + 2-2q} > p$. Since $\hat{f}_{1:\frac{n}{2}}$ and $f$ are both in the Sobolev space $W^{s,p}(\Omega) \subseteq L^{p^\ast}(\Omega)$, we may further deduce that $g_{1:\frac{n}{2}} = \hat{f}_{1:\frac{n}{2}} - f \in W^{s,p}(\Omega) \subseteq L^{p^\ast}(\Omega) \subseteq L^{\frac{2p^\ast}{p^\ast+2-2q}}(\Omega)$. Given that $q \geq 1 \Rightarrow \frac{p^\ast}{p^\ast + 2-2q} \geq 1$, we have $g_{1:\frac{n}{2}}^{2} \in L^{\frac{p^\ast}{p^\ast+2-2q}}(\Omega)$. Then we may apply H\"older's inequality (Lemma \ref{lem: holder ineq}) to $g_{1:\frac{n}{2}}^{2} \in L^{\frac{p^\ast}{p^\ast+2-2q}}(\Omega)$ and $f^{2q-2} \in L^{\frac{p^\ast}{2q-2}}(\Omega)$, which yields the following upper bound:  

\begin{equation}
\label{eqn: App C equation C8}
\begin{aligned}
\mathlarger{\mathlarger{\mathbb{E}}}_{\substack{\{x_i\}_{i=1}^{\frac{n}{2}}}}\Bigg[\int_{\Omega}|g_{1:\frac{n}{2}}^2(x)f^{2q-2}(x)|dx\Bigg] &= \mathlarger{\mathlarger{\mathbb{E}}}_{\substack{\{x_i\}_{i=1}^{\frac{n}{2}}}}\Big[||g_{1:\frac{n}{2}}^2f^{2q-2}||_{L^1(\Omega)}\Big] \\
&\leq \mathlarger{\mathlarger{\mathbb{E}}}_{\substack{\{x_i\}_{i=1}^{\frac{n}{2}}}}\Bigg[\Big||g_{1:\frac{n}{2}}^2\Big||_{L^{\frac{p^\ast}{p^\ast+2-2q}}(\Omega)}\Big||f^{2q-2}\Big||_{L^{\frac{p^\ast}{2q-2}}(\Omega)}\Bigg] \\
&\leq \Big||f\Big||_{L^{p^\ast}(\Omega)}^{2q-2}\mathlarger{\mathlarger{\mathbb{E}}}_{\substack{\{x_i\}_{i=1}^{\frac{n}{2}}}}\Bigg[\Big||g_{1:\frac{n}{2}}\Big||_{L^{\frac{2p^\ast}{p^\ast+2-2q}}(\Omega)}^2\Bigg].     
\end{aligned}    
\end{equation}
Note that the function $\omega(t) = t^{\frac{p^\ast+2-2q}{p^\ast}}$ is concave since $q \geq 1$. Moreover, using the given assumption $s \in (\frac{d(2q-p)}{2pq},\frac{d(2q-p)}{p(2q-2)})$ we get that $\frac{pd}{d-sp} > 2q$, which further yields 
$$\frac{p^\ast+2-2q}{2p^\ast} = \frac{\frac{pd}{d-sp}+2-2q}{2\frac{pd}{d-sp}} > \frac{2}{2\frac{pd}{d-sp}} = \frac{d-sp}{pd},$$  
\emph{i.e,} $\frac{(p^\ast+2-2q)}{2p^\ast} \in (\frac{d-sp}{pd},1]$. Hence, we may apply Jensen's inequality and (\ref{oracle: assumption}) in Assumption \ref{assumption: existence of oracle} to upper-bound the last term in (\ref{eqn: App C equation C8}) as follows:
\begin{equation}
\label{eqn: App C equation C9}
\begin{aligned}
\mathlarger{\mathlarger{\mathbb{E}}}_{\substack{\{x_i\}_{i=1}^{\frac{n}{2}}}}\Bigg[\Big||g_{1:\frac{n}{2}}\Big||_{L^{\frac{2p^\ast}{p^\ast+2-2q}}(\Omega)}^2\Bigg] &= \mathlarger{\mathlarger{\mathbb{E}}}_{\substack{\{x_i\}_{i=1}^{\frac{n}{2}}}}\Bigg[\Bigg(\Big||g_{1:\frac{n}{2}}\Big||_{L^{\frac{2p^\ast}{p^\ast+2-2q}}(\Omega)}^{\frac{2p^\ast}{p^\ast+2-2q}}\Bigg)^{\frac{p^\ast+2-2q}{p^\ast}}\Bigg] \\
&\leq \mathlarger{\mathlarger{\mathbb{E}}}_{\substack{\{x_i\}_{i=1}^{\frac{n}{2}}}}\Bigg[\Big||g_{1:\frac{n}{2}}\Big||_{L^{\frac{2p^\ast}{p^\ast+2-2q}}(\Omega)}^{\frac{2p^\ast}{p^\ast+2-2q}}\Bigg]^{\frac{(p^\ast+2-2q)}{p^\ast}}  \\
&= \mathlarger{\mathlarger{\mathbb{E}}}_{\substack{\{x_i\}_{i=1}^{\frac{n}{2}}}}\Bigg[\Big||\hat{f}_{1:\frac{n}{2}}-f\Big||_{L^{\frac{2p^\ast}{p^\ast+2-2q}}(\Omega)}^{\frac{2p^\ast}{p^\ast+2-2q}}\Bigg]^{\frac{(p^\ast+2-2q)}{p^\ast}}\\
&\lesssim \Big((\frac{n}{2})^{-\frac{s}{d}+(\frac{1}{p}-\frac{p^\ast+2-2q}{2p^\ast})_{+}}\Big)^{2} \\
&\lesssim n^{-\frac{2s}{d}+2(\frac{1}{p}-\frac{p^\ast+2-2q}{2p^\ast})_{+}}.
\end{aligned}
\end{equation}
In order to simplify the last expression in (\ref{eqn: App C equation C9}), let's recall the fact that $p^\ast \in (2q,\frac{p(2q-2)}{p-2})$ proved above. This gives us that $p^\ast(p-2) < p(2q-2) \Rightarrow 2p^\ast > p(p^\ast+2-2q)$, \emph{i.e,} $\frac{1}{p}> \frac{p^\ast+2-2q}{2p^\ast}$. Then we may simplify the power term in the last expression of (\ref{eqn: App C equation C9}) as follows:
$$-\frac{2s}{d}+2\Big(\frac{1}{p}-\frac{p^\ast+2-2q}{2p^\ast}\Big)_{+} = -\frac{2s}{d}+\frac{2}{p}-\Big(1+\frac{2}{p^\ast}-\frac{2q}{p^\ast}\Big) = \frac{2q}{p^\ast}-1 = 2q\Big(\frac{1}{p}-\frac{s}{d}\Big)-1.$$
Now let's substitute (\ref{eqn: App C equation C9}) into (\ref{eqn: App C equation C8}), which gives us the final upper bound on the second expected integral in (\ref{eqn: upper bound on the variance by two integrals}) under the assumption that $s \in (\frac{d(2q-p)}{2pq},\frac{d(2q-p)}{p(2q-2)})$:
\begin{equation}
\label{eqn: App C final upper bound on C3_2, Case 2}
\mathlarger{\mathlarger{\mathbb{E}}}_{\substack{\{x_i\}_{i=1}^{\frac{n}{2}}}}\Bigg[\int_{\Omega}|g_{1:\frac{n}{2}}^2(x)f^{2q-2}(x)|dx\Bigg] \lesssim n^{2q(\frac{1}{p}-\frac{s}{d})-1}.    
\end{equation}
Combining the upper bounds derived in (\ref{eqn: App C bound part 1}), (\ref{eqn: App C final upper bound on C3_2, Case 0}), (\ref{eqn: App C final upper bound on C3_2, Case 1}) and (\ref{eqn: App C final upper bound on C3_2, Case 2}) finally allows us to upper bound the expected variance $\mathlarger{\mathlarger{\mathbb{E}}}_{\substack{\{x_i\}_{i=1}^{\frac{n}{2}}}}\Big[\text{Var}(\hat{f}_{1:\frac{n}{2}}^q-f^q)\Big]$ as below:
\begin{equation}
\label{eqn: App C equation C11}
\mathlarger{\mathlarger{\mathbb{E}}}_{\substack{\{x_i\}_{i=1}^{\frac{n}{2}}}}\Big[\text{Var}(\hat{f}_{1:\frac{n}{2}}^q-f^q)\Big] \lesssim n^{2q(\frac{1}{p}-\frac{s}{d})-1} + \max\{n^{-\frac{2s}{d}},n^{2q(\frac{1}{p}-\frac{s}{d})-1}\}.
\end{equation}
Finally, substituting (\ref{eqn: App C equation C11}) into \ref{eqn: variance decomposition C2}) derived at the beginning gives us the final upper bound:
\begin{equation}
\begin{aligned}
&\mathlarger{\mathlarger{\mathbb{E}}}_{\substack{\{x_i\}_{i=1}^{n},\{y_i\}_{i=1}^{n}}}\Bigg[\left|\hat{H}_{C}^q\Big(\{x_i\}_{i=1}^{n},\{y_i\}_{i=1}^{n}\Big) - I_f^q\right|\Bigg]\\
&\leq \sqrt{\mathlarger{\mathlarger{\mathbb{E}}}_{\substack{\{x_i\}_{i=1}^{n},\{y_i\}_{i=1}^{n}}}\Bigg[\left|\hat{H}_{C}^q\Big(\{x_i\}_{i=1}^{n},\{y_i\}_{i=1}^{n}\Big) - I_f^q\right|^2\Bigg]} = \sqrt{\frac{2}{n}\mathlarger{\mathlarger{\mathbb{E}}}_{\substack{\{x_i\}_{i=1}^{\frac{n}{2}}}}\Big[\text{Var}(\hat{f}_{1:\frac{n}{2}}^q-f^q)\Big]}\\
&\lesssim n^{-\frac{1}{2}}\sqrt{n^{2q(\frac{1}{p}-\frac{s}{d})-1} + \max\{n^{-\frac{2s}{d}},n^{2q(\frac{1}{p}-\frac{s}{d})-1}\}} \lesssim \max\{n^{-\frac{s}{d}-\frac{1}{2}},n^{-q(\frac{s}{d}-\frac{1}{p})-1}\}. 
\end{aligned}
\end{equation}
This concludes our proof of Theorem \ref{thm: upper bound control variate}. 

\subsection{Proof of Theorem \ref{thm: upper bound truncate MC} (Truncated Monte Carlo)}
\label{app: proof of upper bound in sec 3.2}
In this subsection, we provide a complete proof of Theorem \ref{thm: upper bound truncate MC}. For any fixed parameter $M>0$, we may divide $\Omega$ into the following two regions:
\begin{equation}
\label{eqn: App C equation C13, def of truncate regions}
\Omega_{M}^{+} := \{x \in \Omega: |f(x)| \geq M\}, \ \Omega_{M}^{-}:= \{x \in \Omega: |f(x)| < M\},    
\end{equation}
where $\Omega_{M}^{+} \cap \Omega_{M}^{-} = \varnothing$ and $\Omega_{M}^{+} \cup \Omega_{M}^{-} = \Omega$. Let $f_{M}(x) := \max\Big\{\min\{f(x),M\},-M\Big\} \ (\forall \ x \in \Omega)$ denote a truncated version of the given function $f$, where $M$ is the threshold. Also, we use the following expression to denote the expectation of the $q$-th power of the truncated function $f_{M}$ with respect to the uniform distribution on $\Omega$:
\begin{equation}
\label{eqn: App C equation C14, def of truncate function's expectation}
\begin{aligned}
\mathbb{E}(f_{M}^q(x)) &= \int_{\Omega}\max\Big\{\min\{f(x),M\},-M\Big\}^qdx = \int_{\Omega_{M}^{+}}M^qdx + \int_{\Omega_{M}^{-}}f(x)^qdx, 
\end{aligned} 
\end{equation}
where the last identity in (\ref{eqn: App C equation C14, def of truncate function's expectation}) above follows from our definition of the two regions defined in (\ref{eqn: App C equation C13, def of truncate regions}). In a similar way, we can define the variance of the function $f_{M}^q$ as below:
\begin{equation}
\label{eqn: App C equation C15, def of truncate function's variance}
\begin{aligned}
&\text{Var}(f_{M}^q(x)) = \mathbb{E}(f_{M}^{2q}(x))-\mathbb{E}(f_{M}^q(x))^2\\
&= \int_{\Omega}\max\Big\{\min\{f(x),M\},-M\Big\}^{2q}dx - \Big(\int_{\Omega}\max\Big\{\min\{f(x),M\},-M\Big\}^{q}dx\Big)^2. 
\end{aligned} 
\end{equation}
Furthermore, as $\{x_i\}_{i=1}^{n}$ are identical and independent samples of the uniform distribution on $\Omega$, we have that for any $1 \leq i \leq n$, the following identity holds
\begin{equation}
\label{eqn: App C equation C16}
\begin{aligned}
&\mathlarger{\mathbb{E}}_{\substack{\{x_i\}_{i=1}^{n}, \{y_i\}_{i=1}^{n}}}\Big[\hat{H}_{M}^q\Big(\{x_i\}_{i=1}^{n},\{y_i\}_{i=1}^{n}\Big)\Big]\\
&= \mathlarger{\mathbb{E}}_{\substack{\{x_i\}_{i=1}^{n}, \{y_i\}_{i=1}^{n}}}\Big[\frac{1}{n}\sum_{i=1}^{n}\max\Big\{\min\{y_i,M\},-M\Big\}^q\Big] \\
&= \mathlarger{\mathbb{E}}_{x_i}\Big[\max\Big\{\min\{f(x_i),M\},-M\Big\}^q\Big] =\mathbb{E}_{x_i}[f_{M}^q(x_i)] = \mathbb{E}(f_{M}^q(x)).   
\end{aligned}    
\end{equation}
Now we may use (\ref{eqn: App C equation C16}) and the bias-variance decomposition to derive an upper bound on the squared expected risk of the estimator $\hat{H}_{M}^q$ as follows:
\begin{equation}
\label{eqn: App C equation C17, bias variance decomposition}
\begin{aligned}
&\mathlarger{\mathlarger{\mathbb{E}}}_{\substack{\{x_i\}_{i=1}^{n},\{y_i\}_{i=1}^{n}}}\Bigg[\left|\hat{H}_{M}^q\Big(\{x_i\}_{i=1}^{n},\{y_i\}_{i=1}^{n}\Big) - I_f^q\right|^2\Bigg] \\
&= \mathlarger{\mathlarger{\mathbb{E}}}_{\substack{\{x_i\}_{i=1}^{n},\{y_i\}_{i=1}^{n}}}\Bigg[\left|\hat{H}_{M}^q\Big(\{x_i\}_{i=1}^{n},\{y_i\}_{i=1}^{n}\Big) -\mathbb{E}(f_{M}^q(x)) + \mathbb{E}(f_{M}^q(x)) - I_f^q\right|^2\Bigg]\\
&\leq 2\mathlarger{\mathlarger{\mathbb{E}}}_{\substack{\{x_i\}_{i=1}^{n},\{y_i\}_{i=1}^{n}}}\Bigg[\left|\frac{1}{n}\sum_{i=1}^{n}\max\Big\{\min\{y_i,M\},-M\Big\}^q -\mathbb{E}(f_{M}^q(x)) \right|^2\Bigg]\\
&+2\mathlarger{\mathlarger{\mathbb{E}}}_{\substack{\{x_i\}_{i=1}^{n},\{y_i\}_{i=1}^{n}}}\Bigg[\left|\mathbb{E}(f_{M}^q(x)) - I_f^q\right|^2\Bigg]\\    
\end{aligned}    
\end{equation}
where the first and the second term in the last line of (\ref{eqn: App C equation C17, bias variance decomposition}) above denotes the variance and the bias part, respectively. Again, we define $p^\ast = (\max\{\frac{1}{p}-\frac{s}{d},0\})^{-1}$, \emph{i.e,} $p^\ast = \frac{pd}{d-sp}$ when $s<\frac{d}{p}$ and $p^\ast=\infty$ otherwise. Under the assumption that $s<\frac{2dq-dp}{2pq} < \frac{d}{p}$, we have $p^\ast=\frac{pd}{d-sp} \in (p,2q)$. Moreover, from Sobolev Embedding Theorem (Lemma \ref{lem: sobolev embedding}), we have that $f \in W^{s,p}(\Omega) \subseteq L^{p^\ast}(\Omega)$.

On the one hand, since $p < 2q$, we can deduce that $|f(x)|^{2q} \leq M^{2q-p^\ast}|f(x)|^{p^\ast}$ for any $x \in \Omega_{M}^{-}$ and $M^{2q} \leq M^{2q-p^\ast}|f(x)|^{p^\ast}$ for any $x \in \Omega_{M}^+$, which helps us upper bound the variance part as below:
\begin{equation}
\label{eqn: App C equation C18, bound on variance}
\begin{aligned}
&\mathlarger{\mathlarger{\mathbb{E}}}_{\substack{\{x_i\}_{i=1}^{n},\{y_i\}_{i=1}^{n}}}\Bigg[\left|\frac{1}{n}\sum_{i=1}^{n}\max\Big\{\min\{y_i,M\},-M\Big\}^q -\mathbb{E}(f_{M}^q(x)) \right|^2\Bigg]\\
&= \mathlarger{\mathlarger{\mathbb{E}}}_{\substack{\{x_i\}_{i=1}^{n}}}\Bigg[\left|\frac{1}{n}\sum_{i=1}^{n}\Big(f_{M}^q(x_i) -\mathlarger{\mathbb{E}}_{x_i}\Big[f_{M}^q(x_i)\Big]\Big)\right|^2\Bigg]\\
&=\frac{1}{n^2}\sum_{i=1}^{n}\mathlarger{\mathbb{E}}_{x_i}\Big[\Big(f_{M}^q(x_i) -\mathlarger{\mathbb{E}}_{x_i}\Big[f_{M}^q(x_i)\Big]\Big)^2\Big]=\frac{1}{n}\text{Var}(f_{M}^q(x)) \\
&\leq \frac{1}{n}\mathbb{E}(f_{M}^{2q}(x)) =\frac{1}{n}\Big(\int_{\Omega_{M}^{+}}M^{2q}dx + \int_{\Omega_{M}^{-}}f(x)^{2q}dx\Big)\\
&\leq \frac{1}{n}\Big(\int_{\Omega_{M}^{+}}M^{2q-p^\ast}|f(x)|^{p^\ast}dx + \int_{\Omega_{M}^{-}}M^{2q-p^\ast}|f(x)|^{p^\ast}dx\Big)\\
&= \frac{1}{n}\int_{\Omega}M^{2q-p^\ast}|f(x)|^{p^\ast}dx \lesssim \frac{M^{2q-p^\ast}}{n},
\end{aligned}    
\end{equation}
where the last step of (\ref{eqn: App C equation C18, bound on variance}) above follows from the fact that $f \in W^{s,p}(\Omega) \subseteq L^{p^\ast}(\Omega)$.

On the other hand, using the fact that $p^\ast > p > q \Rightarrow |f(x)|^q \leq M^{q-p^\ast}|f(x)|^{p^\ast}$ for any $x \in \Omega_{M}^+$, we may upper-bound the bias part as follows:
\begin{equation}
\label{eqn: App C equation C19, bound on bias}
\begin{aligned}
&\mathlarger{\mathlarger{\mathbb{E}}}_{\substack{\{x_i\}_{i=1}^{n},\{y_i\}_{i=1}^{n}}}\Bigg[\left|\mathbb{E}(f_{M}^q(x)) - I_f^q\right|^2\Bigg]\\
&= \left|\int_{\Omega_{M}^{+}}M^qdx + \int_{\Omega_{M}^{-}}f(x)^qdx -\int_{\Omega_{M}^{+}}f^q(x)dx - \int_{\Omega_{M}^{-}}f(x)^qdx \right|^2\\
&=\left|\int_{\Omega_{M}^{+}}\Big(M^q-f^q(x)\Big)dx \right|^2 \leq \left|\int_{\Omega_{M}^{+}}\Big|M^q-f^q(x)\Big|dx \right|^2 \leq \left|\int_{\Omega_{M}^{+}}\Big(M^q+|f(x)|^q\Big)dx \right|^2 \\
&\leq \left|2\int_{\Omega_{M}^{+}}|f(x)|^qdx \right|^2 \lesssim \left|\int_{\Omega_{M}^+}M^{q-p^\ast}|f(x)|^{p^\ast}dx \right|^2 \leq M^{2q-2p^\ast}\left|\int_{\Omega}|f(x)|^{p^\ast}dx \right|^2\\ &\lesssim M^{2q-2p^\ast},
\end{aligned}    
\end{equation}
where the last step above again follows from the fact that $f \in W^{s,p}(\Omega) \subseteq L^{p^\ast}(\Omega)$. By substituting (\ref{eqn: App C equation C18, bound on variance}) and (\ref{eqn: App C equation C19, bound on bias}) into (\ref{eqn: App C equation C17, bias variance decomposition}), we obtain that
\begin{equation}
\begin{aligned}
&\mathlarger{\mathlarger{\mathbb{E}}}_{\substack{\{x_i\}_{i=1}^{n},\{y_i\}_{i=1}^{n}}}\Bigg[\left|\hat{H}_{M}^q\Big(\{x_i\}_{i=1}^{n},\{y_i\}_{i=1}^{n}\Big) - I_f^q\right|\Bigg]\\
&\leq \sqrt{\mathlarger{\mathlarger{\mathbb{E}}}_{\substack{\{x_i\}_{i=1}^{n},\{y_i\}_{i=1}^{n}}}\Bigg[\left|\hat{H}_{M}^q\Big(\{x_i\}_{i=1}^{n},\{y_i\}_{i=1}^{n}\Big) - I_f^q\right|^2\Bigg]} \lesssim \sqrt{\frac{M^{2q-p^\ast}}{n} + M^{2q-2p^\ast}}.
\end{aligned}
\end{equation}
By balancing the variance part $\frac{M^{2q-p^\ast}}{n}$ and the bias part $M^{2q-2p^\ast}$ above, we may get the optimal choice of $M$ as follows: $\frac{M^{2q-p^\ast}}{n} = M^{2q-2p^\ast} \Rightarrow M = \Theta(n^{\frac{1}{p^\ast}})$. Plugging in the optimal choice of $M$ gives us the final upper bound:
\begin{equation}
\begin{aligned}
\mathlarger{\mathlarger{\mathbb{E}}}_{\substack{\{x_i\}_{i=1}^{n},\{y_i\}_{i=1}^{n}}}\Bigg[\left|\hat{H}_{M}^q\Big(\{x_i\}_{i=1}^{n},\{y_i\}_{i=1}^{n}\Big) - I_f^q\right|\Bigg] \lesssim \sqrt{n^{\frac{2q-2p^\ast}{p^\ast}}} = n^{\frac{q}{p^\ast}-1} = n^{-q(\frac{s}{d}-\frac{1}{p})-1},    
\end{aligned}    
\end{equation}
which finishes our proof of Theorem \ref{thm: upper bound truncate MC}.

\section{Proof of Minimax Lower and Upper Bounds in Section \ref{sec 4: case study}}
\label{app: proof of upper and lower in sec 4}
This section is organized as follows. The first subsection consists of one important lemma used in our proof. In the second subsection, we provide complete proof for the minimax optimal lower bound on the estimation of integrals under any level of noise. In the third subsection, a complete proof for the upper bound on the estimation of integrals is given. 

\subsection{A Key Lemma for Establishing the Upper Bound on Integral Estimation}

\begin{lemma}[Bound on the Expected $k$-Nearest Neighbor Distance: Theorem 2.4, \citep{biau2015lectures}]
\label{lem: upper bound on expected distance of the k-th distance}
Assume that $x_{1},x_{2},\cdots,x_{n}$ are independent and identical samples from the uniform distribution on the domain $\Omega =[0,1]^d$. For any $k \in \{1,2,\cdots,n\}$ and $z \in \Omega$, we use $x_{i^{(z)}_k}$ to denote the $k$-th nearest neighbor of $z$ among $\{x_i\}_{i=1}^{n}$.
When $z$ is also uniformly distributed over the domain $\Omega$, we have the following upper bound on the expected distance between $z$ and $x_{i^{(z)}_k}$:
\begin{equation}
\label{eq: upper bound on expected distance of the k-th distance}
\mathbb{E}_{z,\{x_{i}\}_{i=1}^{n}}\Big[||z-x_{i^{(z)}_k}||^2\Big] \lesssim (\frac{k}{n})^{\frac{2}{d}}.    
\end{equation}
\end{lemma}

\subsection{Proof of Theorem \ref{thm: lower bound for integral} (Lower Bound on Integral Estimation)}
\label{app: proof of lower bound in thm 4.1}
Here we present a comprehensive proof of the two lower bounds given in Theorem \ref{thm: lower bound for integral} above by applying the method of two fuzzy hypotheses (Lemma \ref{lem: method of two fuzzy hypo}). Below we again use $\vec{x}:=(x_1,x_2,\cdots,x_n)$ and $\vec{y}:=(y_1,y_2,\cdots,y_n)$ to denote the two $n$-dimensional vectors formed by the quadrature points and observed function values. Since our lower bound in Theorem \ref{thm: lower bound for integral} consists of two terms, we need to prove the two bounds in the following two separate cases:

(Case I) For the first lower bound in (\ref{eqn: lower bound for integral estimation}), let's consider two constant functions $g_0$ and $g_1$ defined as follows:
\begin{equation}
g_0(x) \equiv 0 \ (\forall \ x \in \Omega), \ g_1(x) \equiv n^{-\gamma-\frac{1}{2}} \ (\forall \ x \in \Omega)
\end{equation}
Clearly we have $g_0,g_1 \in C^{s}(\Omega)$. Then let's take $\mu_k$ to be a Dirac delta measure supported on the set $\{g_j\}$, \emph{i.e,} $\mu_k(\{g_k\}) = 1$, for $k \in \{0,1\}$. By picking $c=\Delta =\frac{1}{2}I_{g_1}=\frac{1}{2}n^{-\gamma-\frac{1}{2}}$ and $\beta_0=\beta_1=0$, we then obtain that 
\begin{equation}
\begin{aligned}
\mu_0(f \in W^{s,p}(\Omega): I_{f} &\leq c-\Delta) = \mu_0(I_{f} \leq 0) =1= 1-\beta_0, \\
\mu_1(f \in W^{s,p}(\Omega): I_{f} &\geq c+\Delta) =\mu_1(I_{f} \geq I_{g_1})=1= 1-\beta_1,     
\end{aligned}    
\end{equation}
which indicates that (\ref{eq: separated hypothesis (two prior)}) holds true. Now let's consider bounding the KL divergence between the two marginal distributions $\mathbb{P}_0,\mathbb{P}_1$ associated with $\mu_0,\mu_1$, respectively. Given that the quadrature points $\{x_i\}_{i=1}^{n}$ and the observational noises $\{\epsilon_i\}_{i=1}^{n}$ are independent and identical samples from the uniform distribution on $\Omega$ and the normal distribution $\mathcal{N}(0,n^{-2\gamma})$, we can write the marginal distributions in an explicit form as follows: 

\begin{equation}
\label{eqn: marginal distributions thm 4.1 case 1}
\begin{aligned}
\mathbb{P}_0(\vec{x},\vec{y}) = \prod_{i=1}^{n}\Big(\frac{1}{\sqrt{2\pi}n^{-\gamma}}e^{-\frac{1}{2n^{-2\gamma}}y_{i}^2}\Big), \  \mathbb{P}_1(\vec{x},\vec{y}) = \prod_{i=1}^{n}\Big(\frac{1}{\sqrt{2\pi}n^{-\gamma}}e^{-\frac{1}{2n^{-2\gamma}}(y_i-n^{-\gamma-\frac{1}{2}})^2}\Big).   
\end{aligned}    
\end{equation}

From (\ref{eqn: marginal distributions thm 4.1 case 1}) we can see that $\mathbb{P}_0$ and $\mathbb{P}_1$ are two $n$-dimensional normal distributions having the same covariance matrix but different mean vectors. Computing the KL divergence between them and applying Pinsker's inequality then give us that
\begin{equation}
\label{thm 4.1 case 1: bound on TV}
TV(\mathbb{P}_0 || \mathbb{P}_1) \leq  \sqrt{\frac{1}{2}KL(\mathbb{P}_0 || \mathbb{P}_1)} = \sqrt{\frac{n(n^{-\gamma-\frac{1}{2}})^2}{4n^{-2\gamma}}} = \frac{1}{2}.   
\end{equation}
Substituting (\ref{thm 4.1 case 1: bound on TV}), $\Delta =\frac{1}{2}I_{g_1} = \frac{1}{2}n^{-\gamma-\frac{1}{2}}$ and $\beta_0=\beta_1=0$ into (\ref{eq: bound on misclassification probability (two prior)}) and applying Markov's inequality yield the final lower bound
\begin{equation}
\label{eqn: thm 4.1 case 1 final bound}
\begin{aligned}
&\inf_{\hat{H} \in \mathcal{H}_{n}^{f}}\sup_{f \in C^{s}(\Omega)}\mathlarger{\mathlarger{\mathbb{E}}}_{\substack{\{x_i\}_{i=1}^{n},\{y_i\}_{i=1}^{n}}}\Bigg[\left|\hat{H}\Big(\{x_i\}_{i=1}^{n},\{y_i\}_{i=1}^{n}\Big) - I_f\right|\Bigg]\\
&\geq \Delta \inf_{\hat{H} \in \mathcal{H}_{n}^{f}}\sup_{f \in C^{s}(\Omega)}\mathlarger{\mathlarger{\mathbb{P}}}_{\substack{\{x_i\}_{i=1}^{n},\{y_i\}_{i=1}^{n}}}\Bigg[\left|\hat{H}\Big(\{x_i\}_{i=1}^{n},\{y_i\}_{i=1}^{n}\Big) - I_f\right| \geq \Delta\Bigg]\\
&\geq \frac{1}{2}I_{g_1}\frac{1-TV(\mathbb{P}_0 || \mathbb{P}_1)-\beta_0-\beta_1}{2} \geq \frac{1}{8}n^{-\gamma-\frac{1}{2}} \gtrsim n^{-\gamma -\frac{1}{2}},
\end{aligned}    
\end{equation}
which is exactly the first term in the RHS of (\ref{eqn: lower bound for integral estimation}).

(Case II) For the second lower bound in (\ref{eqn: lower bound for integral estimation}), our proof is similar to the proof of the second lower bound in Theorem \ref{thm: lower bound q-th moment} presented in Appendix \ref{app: proof of lower bound in sec 2} above. Again,  we select $m= (200n)^{\frac{1}{d}}$ and divide the domain $\Omega$ into $m^d$ small cubes $\Omega_1,\Omega_2,\cdots,\Omega_{m^d}$, each of which has side length $m^{-1}$. For any $1 \leq j \leq m^d$, we use $c_j$ to denote center of the cube $\Omega_j$. Then let's consider the same bump function $K$ defined in (\ref{eqn: def of compact C_infty function K_0 on [-1,1]^d}) and (\ref{eqn: def of K, rescaled K_0 on [-1/2,1/2]^d}) above, which satisfies $\text{supp}(K) \subseteq [-\frac{1}{2},\frac{1}{2}]^d$ and $K \in C^{\infty}([-\frac{1}{2},\frac{1}{2}]^d)$. In an analogous way, for any $1 \leq j \leq m^d$, we associate each cube $\Omega_j$ with a bump function $f_j$ defined as follows:
\begin{equation}
\label{eqn: def of bump functions in thm 4.1}
\begin{aligned}
f_j(x) &= 
\begin{cases}
m^{-s}K(m(x-c_j)) \ (x \in \Omega_j), \\
0 \ (\text{otherwise}),
\end{cases}    
\end{aligned}
\end{equation}
where $\text{supp}(f_j) \subseteq \Omega_j, \ f_j \in C^{\infty}(\Omega)$ and $f_j(x) \geq 0 \ (\forall \ x \in \Omega)$. Then let's consider the following finite set of $2^{m^d}$ functions:
\begin{equation}
\mathcal{S} := \Big\{\sum_{j=1}^{m^d}\eta_j f_j: \eta_j \in \{\pm 1\}, \ \forall \ 1 \leq j \leq m^d\Big\}.    
\end{equation}
We will first verify that $\mathcal{S} \subseteq C^{s}(\Omega)$. Fix any element $f_\ast = \sum_{j=1}^{m^d}\eta_j f_j \in \mathcal{S}$. On the one hand, from our construction of the $f_j$'s given in (\ref{eqn: def of bump functions in thm 4.1}) above, we have 
\begin{equation}
\label{eqn: thm 4.1 bound on holder norm part 1}
\begin{aligned}
\max_{|t| \leq \lfloor s \rfloor}||D^{t}f_\ast||_{L^{\infty}(\Omega)} &= \max_{|t| \leq \lfloor s \rfloor}m^{-s + |t|}||D^{t}K||_{L^{\infty}([-\frac{1}{2},\frac{1}{2}]^d)} \\
&\leq \max_{|t| \leq \lfloor s \rfloor}||D^{t}K||_{L^{\infty}([-\frac{1}{2},\frac{1}{2}]^d)}.
\end{aligned}    
\end{equation}
On the other hand, for any $1 \leq i \neq j \leq m^d$, we consider the function $\psi_i f_i +\psi_j f_j$, where the scalars $\psi_j,\psi_j \in \{0,\pm 1\}$.  Now let's may pick $\beta:= \frac{d}{1-\{s\}}$, where $\{s\} = s -\lfloor s \rfloor \in (0,1)$ denotes the fractional part of $s$. Given that $f_j \in C^{\infty}(\Omega)$, we may upper bound the Sobolev norm $||\cdot||_{W^{1,\beta}}$ of the function $D^{t}(\psi_i f_i +\psi_j f_j)$ for any $t \in \mathbb{N}_{0}^d$ satisfying $|t| = \lfloor s \rfloor$ as follows:

\begin{equation}
\begin{aligned}
&\Big||D^{t}(\psi_i f_i +\psi_j f_j)\Big||_{W^{1,\beta}(\Omega)}^\beta = |\psi_i|^\beta\Big||D^{t}f_i\Big||_{W^{1,\beta}(\Omega_i)}^\beta + |\psi_j|^\beta\Big||D^{t}f_j\Big||_{W^{1,\beta}(\Omega_j)}^\beta\\
&\leq \Big||D^{t}f_i\Big||_{L^{\beta}(\Omega_i)}^\beta+\sum_{r=1}^{d}\Big||\frac{\partial}{\partial x_r}D^{t}f_i\Big||_{L^{\beta}(\Omega_i)}^\beta+\Big||D^{\lfloor s \rfloor}f_j\Big||_{L^{\beta}(\Omega_j)}^\beta+\sum_{r=1}^{d}\Big||\frac{\partial}{\partial x_r}D^{t}f_j\Big||_{L^{\beta}(\Omega_j)}^\beta\\
&= \sum_{l \in \{i,j\}}\int_{\Omega_l}\Big(m^{-s+|t|}D^{t}K(m(x-c_l))\Big)^\beta dx \\
&+  \sum_{l \in \{i,j\}}\sum_{r=1}^{d}\int_{\Omega_l}\Big(m^{-s+|t| + 1}\frac{\partial}{\partial x_r}D^{t}K(m(x-c_l))\Big)^\beta dx.
\end{aligned}    
\end{equation}

From our choice of $\beta$ and assumption on the bump function $K$, we may further upper bound the Sobolev norm $\Big||D^{t}(\psi_i f_i +\psi_j f_j)\Big||_{W^{1,\beta}(\Omega)}$ as below:
\begin{equation}
\label{eqn: thm 4.1 check holder condition}
\begin{aligned}
&\Big||D^{t}(\psi_i f_i +\psi_j f_j)\Big||_{W^{1,\beta}(\Omega)}^\beta \leq \sum_{l \in \{i,j\}}m^{-\beta\{s\}}\int_{[-\frac{1}{2},\frac{1}{2}]^d}\Big(D^{t}K(y)\Big)^\beta\frac{1}{m^d}dy \\
&+\sum_{l \in \{i,j\}} dm^{\beta(1-\{s\})}\sup_{|t'| \leq \lfloor s \rfloor +1}\Bigg(\int_{[-\frac{1}{2},\frac{1}{2}]^d}\Big(D^{t'}K(y)\Big)^\beta\frac{1}{m^d}dy\Bigg)\\
&\leq 2m^{-\beta\{s\}-d}\Big||D^{t}K\Big||_{L^{\beta}([-\frac{1}{2},\frac{1}{2}]^d)}^\beta + 2dm^{\beta(1-\{s\})-d}\cdot \sup_{|t'| \leq \lfloor s \rfloor +1}\Big||D^{t'}K\Big||_{L^{\beta}([-\frac{1}{2},\frac{1}{2}]^d)}^\beta\\
&\lesssim \Big||D^{t}K\Big||_{L^{\beta}([-\frac{1}{2},\frac{1}{2}]^d)}^\beta + \sup_{|t'| \leq \lfloor s \rfloor +1}\Big||D^{t'}K\Big||_{L^{\beta}([-\frac{1}{2},\frac{1}{2}]^d)}^\beta, 
\end{aligned}    
\end{equation}
where the last inequality above follows from our choice of $\beta$. From (\ref{eqn: thm 4.1 check holder condition}) and the second part of the Sobolev Embedding Theorem (Lemma 
 \ref{lem: sobolev embedding}), we can deduce that $D^t(\psi_i f_i + \psi_j f_j)  \in C^1(\Omega) \cap W^{1,\frac{d}{1-\{s\}}}(\Omega) \subseteq C^{\{s\}}(\Omega)$ and the following inequality holds:
\begin{equation}
\label{thm 4.1: sobolev embed in holder}
\begin{aligned}
&\Big||D^{t}(\psi_i f_i +\psi_j f_j)\Big||_{C^{\{s\}}(\Omega)} \lesssim \Big||D^{t}(\psi_i f_i +\psi_j f_j)\Big||_{W^{1,\beta}(\Omega)}\\ 
&\lesssim  \Bigg(\sup_{|t'| = \lfloor s \rfloor}\Big||D^{t'}K\Big||_{L^{\beta}([-\frac{1}{2},\frac{1}{2}]^d)}^\beta + \sup_{|t'| = \lfloor s \rfloor +1}\Big||D^{t'}K\Big||_{L^{\beta}([-\frac{1}{2},\frac{1}{2}]^d)}^\beta\Bigg)^{\frac{1}{\beta}},    
\end{aligned}
\end{equation}
Furthermore, combining (\ref{thm 4.1: sobolev embed in holder}) with our construction of the $f_j$'s given in (\ref{eqn: def of bump functions in thm 4.1}) above gives us that

\begin{equation}
\label{eqn: thm 4.1 bound on holder norm part 2}
\begin{aligned}
&\max_{|t|=\lfloor s \rfloor}\sup_{x,y \in \Omega, x \neq y}\frac{|D^{t}f_\ast(x)-D^{t}f_\ast(y)|}{||x-y||^{s-{\lfloor s\rfloor}}} \\
&\leq \max_{\substack{1 \leq i \neq j \leq k \\ \psi_i,\psi_j \in \{0,\pm 1\}}}\max_{|t|=\lfloor s \rfloor}\sup_{x \neq y \in \Omega}\frac{|D^{t}(\psi_i f_i +\psi_j f_j)(x) - D^{t}(\psi_i f_i +\psi_j f_j)(x)|}{||x-y||^{\{s\}}} \\
&\leq \max_{\substack{1 \leq i \neq j \leq k, |t|=\lfloor s \rfloor \\ \psi_i,\psi_j \in \{0,\pm 1\}}}\Big||D^{t}(\psi_i f_i +\psi_j f_j)\Big||_{C^{\{s\}}(\Omega)}\\ &\lesssim  \Bigg(\sup_{|t'| = \lfloor s \rfloor}\Big||D^{t'}K\Big||_{L^{\beta}([-\frac{1}{2},\frac{1}{2}]^d)}^\beta + \sup_{|t'| = \lfloor s \rfloor +1}\Big||D^{t'}K\Big||_{L^{\beta}([-\frac{1}{2},\frac{1}{2}]^d)}^\beta\Bigg)^{\frac{1}{\beta}}
\end{aligned}    
\end{equation} 
Finally, adding the two inequalities (\ref{eqn: thm 4.1 bound on holder norm part 1}) and (\ref{eqn: thm 4.1 bound on holder norm part 2}) gives us that for any $f_\ast \in \mathcal{S}$, we have
\begin{equation}
\begin{aligned}
&||f_\ast||_{C^{s}(\Omega)} = \max_{|t| \leq \lfloor s \rfloor}||D^{t}f_\ast||_{L^{\infty}(\Omega)} + \max_{|t|=\lfloor s \rfloor}\sup_{x,y \in \Omega, x \neq y}\frac{|D^{t}f_\ast(x)-D^{t}f_\ast(y)|}{||x-y||^{s-\lfloor s \rfloor}} \\
&\lesssim \max_{|t| \leq \lfloor s \rfloor}||D^{t}K||_{L^{\infty}([-\frac{1}{2},\frac{1}{2}]^d)}\\
&+ \Bigg(\sup_{|t'| = \lfloor s \rfloor}\Big||D^{t'}K\Big||_{L^{\beta}([-\frac{1}{2},\frac{1}{2}]^d)}^\beta + \sup_{|t'| = \lfloor s \rfloor +1}\Big||D^{t'}K\Big||_{L^{\beta}([-\frac{1}{2},\frac{1}{2}]^d)}^\beta\Bigg)^{\frac{1}{\beta}} \lesssim 1.
\end{aligned}    
\end{equation}
From the arbitrariness of $f_\ast$, we can then deduce that $\mathcal{S} \subseteq C^{s}(\Omega)$, as desired. For any $p \in (0,1)$, below we again use $w_p$ to denote the discrete random variable satisfying $\mathbb{P}(w_p = -1) = p$ and $\mathbb{P}(w_p = 1)=1-p$. Now let's pick $\kappa = \frac{1}{3}\sqrt{\frac{2}{3n}}$ and take $\{w_{j}^{(0)}\}_{j=1}^{m^d}$ and $\{w_{j}^{(1)}\}_{j=1}^{m^d}$ to be independent and identical copies of $w_{\frac{1+\kappa}{2}}$ and $w_{\frac{1-\kappa}{2}}$ respectively. Then we define $\mu_0,\mu_1$ to be two discrete measures supported on the finite set $\mathcal{S}$ such that the following condition holds for any $\eta_j \in \{\pm 1\} \ (1 \leq j \leq m^d)$:
\begin{equation}
\label{eqn: def of two priors_1 in thm 4.1 case 2}
\mu_k\Big(\Big\{\sum_{j=1}^{m^d}\eta_j f_j\Big\}\Big) = \prod_{j=1}^{m^d}\mathbb{P}(w_{j}^{(k)} = \eta_j), \ k \in \{0,1\}.
\end{equation}
Then we proceed to determine the separation distance $\Delta$ between the two priors $\mu_0$ and $\mu_1$. Similar to what we did in the proof of Theorem \ref{thm: lower bound q-th moment}, we need to first define the following quantity $C:= \int_{\Omega_j}f_j(x)dx$, which remains the same for any $1 \leq j \leq m^d$. Moreover, applying (\ref{eqn: def of bump functions in thm 4.1}) helps us evaluate the quantity $C$ directly as follows
\begin{equation}
\label{eqn: thm 4.1 case 2 calculation of sep dist}
\begin{aligned}
C &= \int_{\Omega_j}f_j(x)dx = \int_{\Omega_j}m^{-s}K(m(x-c_j))dx \\
&= m^{-s}\int_{[-\frac{1}{2},\frac{1}{2}]^d}K(y)\frac{1}{m^d}dy = m^{-s-d}||K||_{L^1([-\frac{1}{2},\frac{1}{2}]^d)}.  
\end{aligned}    
\end{equation}

Moreover, by picking $\lambda = \frac{1}{2}$, we may apply Hoeffding's Inequality (Lemma \ref{lem: hoeffding ineq}) to the bounded random variables $\{w_{j}^{(0)}\}_{j=1}^{m^d}$ and $\{w_{j}^{(1)}\}_{j=1}^{m^d}$ to deduce that
\begin{equation}
\label{eqn: thm 4.1 case 2 hoeffding bound}
\begin{aligned}
\mathbb{P}\Big(\sum_{j=1}^{m^d}w_{j}^{(0)} \geq -(1-\lambda)m^d \kappa\Big) &\leq \exp\Big(-\frac{2(\lambda m^d \kappa)^2}{4m^d}\Big) = \exp\Big(-\frac{1}{2}\lambda^2 \kappa^2 m^d\Big), \\
\mathbb{P}\Big(\sum_{j=1}^{m^d}w_{j}^{(1)} \leq (1-\lambda)m^d \kappa\Big) &\leq \exp\Big(-\frac{2(\lambda m^d \kappa)^2}{4m^d}\Big) = \exp\Big(-\frac{1}{2}\lambda^2 \kappa^2 m^d\Big).
\end{aligned}    
\end{equation}
By taking $c:= 0, \Delta := (1-\lambda)\kappa m^dC$ and $\beta_0 =\beta_1 =\exp\Big(-\frac{1}{2}\lambda^2 \kappa^2 m^d\Big)$, we may use (\ref{eqn: thm 4.1 case 2 hoeffding bound}) justified above to get that
\begin{equation}
\begin{aligned}
&\mu_0(f \in C^{s}(\Omega): I_{f} \leq c-\Delta) = \mathbb{P}\Big(\sum_{j=1}^{m^d}I_{w_{j}^{(0)}f_j} \leq \frac{1-(1-\lambda)\kappa}{2}m^d C - \frac{1+(1-\lambda)\kappa}{2}m^d C\Big) \\
&\geq \mathbb{P}\Big(\sum_{j=1}^{m^d}w_{j}^{(0)} \leq -(1-\lambda)m^d \kappa\Big)=1-\mathbb{P}\Big(\sum_{j=1}^{m^d}w_{j}^{(0)} \geq -(1-\lambda)m^d \kappa\Big)\\
&\geq 1-\exp\Big(-\frac{1}{2}\lambda^2 \kappa^2 m^d\Big) = 1-\beta_0, \\
&\mu_1(f \in C^{s}(\Omega): I_{f} \geq c+\Delta) = \mathbb{P}\Big(\sum_{j=1}^{m^d}I_{w_{j}^{(1)}f_j} \geq \frac{1+(1-\lambda)\kappa}{2}m^d C - \frac{1-(1-\lambda)\kappa}{2}m^d C\Big) \\
&\geq \mathbb{P}\Big(\sum_{j=1}^{m^d}w_{j}^{(1)} \geq (1-\lambda)m^d \kappa\Big)=1-\mathbb{P}\Big(\sum_{j=1}^{m^d}w_{j}^{(0)} \leq (1-\lambda)m^d \kappa\Big)\\
&\geq 1-\exp\Big(-\frac{1}{2}\lambda^2 \kappa^2 m^d\Big) = 1-\beta_1,   
\end{aligned}    
\end{equation}
which indicates that (\ref{eq: separated hypothesis (two prior)}) holds true. Now let's consider bounding the KL divergence between the two marginal distributions $\mathbb{P}_0,\mathbb{P}_1$ associated with $\mu_0,\mu_1$, respectively. Applying the fact that $\{x_1,\cdots,x_n\}$ and $\{\epsilon_1,\cdots,\epsilon_n\}$ are identical and independent samples from the uniform distribution on $\Omega$ and the normal distribution $\mathcal{N}(0,n^{-2\gamma})$ allows us to write the marginal distributions in an explicit form as follows: 
\begin{equation}
\label{eqn: marginal distributions thm 4.1 case 2}
\begin{aligned}
\mathbb{P}_0(\vec{x},\vec{y}) &= \prod_{j=1}^{m^d}\Big(\frac{1-\kappa}{2}\prod_{i:x_i \in \Omega_j}\frac{1}{\sqrt{2\pi}n^{-\gamma}}e^{-\frac{(y_i-f_j(x_i))^2}{2n^{-2\gamma}}} + \frac{1+\kappa}{2}\prod_{i:x_i \in \Omega_j}\frac{1}{\sqrt{2\pi}n^{-\gamma}}e^{-\frac{(y_i+f_j(x_i))^2}{2n^{-2\gamma}}}\Big),\\
\mathbb{P}_1(\vec{x},\vec{y}) &= \prod_{j=1}^{m^d}\Big(\frac{1+\kappa}{2}\prod_{i:x_i \in \Omega_j}\frac{1}{\sqrt{2\pi}n^{-\gamma}}e^{-\frac{(y_i-f_j(x_i))^2}{2n^{-2\gamma}}} + \frac{1-\kappa}{2}\prod_{i:x_i \in \Omega_j}\frac{1}{\sqrt{2\pi}n^{-\gamma}}e^{-\frac{(y_i+f_j(x_i))^2}{2n^{-2\gamma}}}\Big). 
\end{aligned}    
\end{equation}
Furthermore, for any $n$ fixed quadrature points $\vec{x} = (x_1,x_2,\cdots,x_n)$, we use $\mathbb{P}_k(\cdot \ | \ \Vec{x})$ to denote the marginal distribution of the observed function values $\vec{y} = (y_1,y_2,\cdots,y_n)$ conditioned on $\vec{x}$ for $k \in \{0,1\}$. Since $\{x_i\}_{i=1}^{n}$ are identically and independently sampled from the uniform distribution on $\Omega$, we have that the two probability densities $\mathbb{P}_{k}(\Vec{x},\vec{y})$ and $\mathbb{P}_k(\Vec{y} \ | \ \Vec{x})$ have the same mathematical expression for any $k \in \{0,1\}$. Then we may further rewrite the KL divergence between the two marginal distributions $\mathbb{P}_0,\mathbb{P}_1$ as follows:
\begin{equation}
\label{thm 4.1: rewrite KL in terms of conditional density}
\begin{aligned}
&KL(\mathbb{P}_{0} || \mathbb{P}_{1}) = \int_{\Omega} \cdots\int_{\Omega}\Big(\int_{-\infty}^{\infty}\cdots\int_{-\infty}^{\infty}\log\Big(\frac{\mathbb{P}_0(\vec{x},\vec{y})}{\mathbb{P}_1(\vec{x},\vec{y})}\Big)\mathbb{P}_0(\vec{x},\vec{y})dy_1 \cdots dy_n\Big)dx_1 \cdots dx_n\\
&= \int_{\Omega} \cdots\int_{\Omega}\Big(\int_{-\infty}^{\infty}\cdots\int_{-\infty}^{\infty}\log\Big(\frac{\mathbb{P}_0(\vec{y} \ | \ \vec{x})}{\mathbb{P}_1(\vec{y} \ | \ \vec{x})}\Big)\mathbb{P}_0(\vec{y} \ | \ \vec{x})dy_1 \cdots dy_n\Big)dx_1 \cdots dx_n\\
&= \int_{\Omega} \cdots\int_{\Omega}\Bigg(KL\Big(\mathbb{P}_0(\cdot \ | \ \vec{x}) || \mathbb{P}_1(\cdot \ | \ \vec{x})\Big)\Bigg)dx_1 \cdots dx_n.
\end{aligned}    
\end{equation}
It now remains to upper bound the KL divergence between the two conditional distributions $\mathbb{P}_0(\cdot \ | \ \Vec{x})$ and $\mathbb{P}_1(\cdot \ | \ \Vec{x})$ for any fixed $\Vec{x}=(x_1,\cdots,x_n)$. In order to derive such an upper bound, we need to introduce the following notations first. For any $n$ quadrature points $\{x_i\}_{i=1}^{n}$, we use $\mathcal{J}_n$ to denote the set of all indices $j$ satisfying that $\Omega_j$ contains at least one of the points in $\{x_i\}_{i=1}^{n}$, \emph{i.e,}
\begin{equation}
\mathcal{J}_{n}:=\mathcal{J}_n(\Vec{x}) = \Big\{j: 1 \leq j \leq m^d \text{ and } \Omega_j \cap \{x_1,\cdots,x_n\} \neq \varnothing\Big\}. 
\end{equation}
Moreover, we use $\vec{\omega}_{\mathcal{J}_n}^{(k)}$ to denote
$|\mathcal{J}_n|$-dimensional vector formed by the random variables $\{\omega_{j}^{(k)}: j \in \mathcal{J}_n\}$ and $p_{\mathcal{J}_n}^{(k)}(\cdot)$ to denote the probability density function of $\vec{\omega}_{\mathcal{J}_n}^{(k)}$, where $k \in \{0,1\}$. From our assumption on the distribution of the weights $\{w_{j}^{(0)}\}_{j=1}^{n}$ and $\{w_{j}^{(1)}\}_{j=1}^{n}$, we have that for any $\Vec{\omega}_{\mathcal{J}_n} \in \{\pm 1\}^{|\mathcal{J}_n|}$, 
\begin{equation}
\label{eqn: thm 4.1 def of weight distributions}
\begin{aligned}
p_{\mathcal{J}_n}^{(0)}(\vec{\omega}_{\mathcal{J}_n}) &=\prod_{j \in \mathcal{J}_n}\Big(\frac{1+\kappa}{2}\Big)^{\frac{1}{2}(1-\omega_{j})}\Big(\frac{1-\kappa}{2}\Big)^{\frac{1}{2}(1+\omega_{j})}\\
p_{\mathcal{J}_n}^{(1)}(\vec{\omega}_{\mathcal{J}_n}) &=\prod_{j \in \mathcal{J}_n}\Big(\frac{1+\kappa}{2}\Big)^{\frac{1}{2}(1+\omega_{j})}\Big(\frac{1-\kappa}{2}\Big)^{\frac{1}{2}(1-\omega_{j})}
\end{aligned}    
\end{equation}
Furthermore, for any fixed quadrature points $\vec{x} = (x_1,\cdots,x_n)$ and weights $\Vec{\omega}_{\mathcal{J}_n} :=\{\omega_j:j \in \mathcal{J}_n\} \subseteq \{\pm 1\}^{|\mathcal{J}_n|}$, we may define the transition kernel $G(\vec{x},\Vec{\omega}_{\mathcal{J}_n})$ as below
\begin{equation}
\label{eqn: thm 4.1 transition kernel def}
\begin{aligned}
G(\vec{x},\Vec{\omega}_{\mathcal{J}_n}) := \prod_{j \in \mathcal{J}_n}\Big(\prod_{i:x_i \in \Omega_j}\frac{1}{\sqrt{2\pi}n^{-\gamma}}e^{-\frac{(y_i+\omega_{j}f_j(x_i))^2}{2n^{-2\gamma}}}\Big)   
\end{aligned}    
\end{equation}
Combining the expressions in (\ref{eqn: marginal distributions thm 4.1 case 2}),(\ref{eqn: thm 4.1 def of weight distributions}) and (\ref{eqn: thm 4.1 transition kernel def}) allows us to rewrite the two conditional distributions $\mathbb{P}_{k}(\cdot \ | \ \Vec{x})$ as below:
\begin{equation}
\label{eq: thm 4.1 marginal in integral form}
\begin{aligned}
\mathbb{P}_{k}(\Vec{y} \ | \Vec{x}) = \mathbb{P}_k(\Vec{x},\Vec{y}) = \int_{\{\pm 1\}^{|\mathcal{J}_n|}}G(\Vec{x},\Vec{\omega}_{\mathcal{J}_n})p_{\mathcal{J}_n}^{(k)}(\vec{\omega}_{\mathcal{J}_n})d\vec{\omega}_{\mathcal{J}_n} 
\end{aligned}    
\end{equation}
where $k \in \{0,1\}$. Applying the data processing inequality (Lemma \ref{lem: data processing ineq}) to (\ref{eq: thm 4.1 marginal in integral form}) above then enables us to derive the following upper bound on $KL\Big(\mathbb{P}_0(\cdot \ | \ \vec{x}) || \mathbb{P}_1(\cdot \ | \ \vec{x})\Big)$ for any $n$ fixed quadrature points $\Vec{x}=(x_1,\cdots,x_n)$:
\begin{equation}
\label{eqn: thm 4.1 bound on KL}
\begin{aligned}
KL\Big(\mathbb{P}_0(\cdot \ | \ \vec{x}) || \mathbb{P}_1(\cdot \ | \ \vec{x})\Big) &\leq KL\Big(p_{\mathcal{J}_n}^{(0)} \ || \ p_{\mathcal{J}_n}^{(1)}\Big) \\
&= |\mathcal{J}_n|\Big(\log\Big(\frac{1+\kappa}{1-\kappa}\Big)\frac{1+\kappa}{2} + \log\Big(\frac{1-\kappa}{1+\kappa}\Big)\frac{1-\kappa}{2}\Big)\\
&\leq n\kappa\log\Big(\frac{1+\kappa}{1-\kappa}\Big)  
\end{aligned}    
\end{equation}
where the equality in (\ref{eqn: thm 4.1 bound on KL}) above follows from the fact that  $\{w_{j}^{(0)}\}_{j=1}^{m^d}$ and $\{w_{j}^{(1)}\}_{j=1}^{m^d}$ are independent and identical copies of $w_{\frac{1+\kappa}{2}}$ and $w_{\frac{1-\kappa}{2}}$ respectively. The last inequality of (\ref{eqn: thm 4.1 bound on KL}) above, however, is deduced from the fact that $m^d=200n > n$, which implies $|\mathcal{J}_n| \leq n$ for any $n$ quadrature points $\{x_i\}_{i=1}^{n}$. Substituting (\ref{eqn: thm 4.1 bound on KL}) into (\ref{thm 4.1: rewrite KL in terms of conditional density}) and applying Pinkser's inequality yields the final upper bound on the TV distance between $\mathbb{P}_0$ and $\mathbb{P}_1$:

\begin{equation}
\label{thm 4.1 case 2: bound on TV}
\begin{aligned}
TV(\mathbb{P}_0 || \mathbb{P}_1) &\leq  \sqrt{\frac{1}{2}KL(\mathbb{P}_0 || \mathbb{P}_1)} \leq \sqrt{\int_{\Omega} \cdots\int_{\Omega}\frac{n\kappa}{2}\log\Big(\frac{1+\kappa}{1-\kappa}\Big)dx_1 \cdots dx_n} \\
&= \sqrt{\frac{n\kappa}{2}\log\Big(\frac{1+\kappa}{1-\kappa}\Big)} \leq \sqrt{\frac{3n}{2}}\kappa= \frac{1}{3}.     
\end{aligned}  
\end{equation}

Finally, by substituting (\ref{eqn: thm 4.1 case 2 calculation of sep dist}), (\ref{thm 4.1 case 2: bound on TV}), $\Delta =(1-\lambda)\kappa m^d C$ and $\beta_0=\beta_1=\exp\Big(-\frac{1}{2}\lambda^2 \kappa^2 m^d\Big) = \exp(-\frac{50}{27}) < \frac{1}{6}$ into (\ref{eq: bound on misclassification probability (two prior)}) and applying Markov's inequality, we obtain the final lower bound
\begin{equation}
\label{eqn: thm 4.1 case 2 final bound}
\begin{aligned}
&\inf_{\hat{H} \in \mathcal{H}_{n}^{f}}\sup_{f \in C^{s}(\Omega)}\mathlarger{\mathlarger{\mathbb{E}}}_{\substack{\{x_i\}_{i=1}^{n},\{y_i\}_{i=1}^{n}}}\Bigg[\left|\hat{H}\Big(\{x_i\}_{i=1}^{n},\{y_i\}_{i=1}^{n}\Big) - I_f\right|\Bigg]\\
&\geq \Delta \inf_{\hat{H} \in \mathcal{H}_{n}^{f}}\sup_{f \in C^{s}(\Omega)}\mathlarger{\mathlarger{\mathbb{P}}}_{\substack{\{x_i\}_{i=1}^{n},\{y_i\}_{i=1}^{n}}}\Bigg[\left|\hat{H}\Big(\{x_i\}_{i=1}^{n},\{y_i\}_{i=1}^{n}\Big) - I_f\right| \geq \Delta\Bigg]\\
&\geq (1-\lambda)\kappa m^d C \frac{1-TV(\mathbb{P}_0 || \mathbb{P}_1)-\beta_0-\beta_1}{2} \geq \frac{1}{2}\frac{\sqrt{2}}{3\sqrt{3n}}\cdot(200n) \cdot \frac{C}{6} \\
&\gtrsim \sqrt{n}C \gtrsim \sqrt{n}(200n)^{-\frac{s+d}{d}}||K||_{L^1([-\frac{1}{2},\frac{1}{2}]^d)} \gtrsim n^{-\frac{s}{d}-\frac{1}{2}}, 
\end{aligned}    
\end{equation}
which is exactly the second term in the RHS of (\ref{eqn: lower bound for integral estimation}). Combining the two lower bounds proved in (\ref{eqn: thm 4.1 case 1 final bound}) and (\ref{eqn: thm 4.1 case 2 final bound}) concludes our proof of Theorem \ref{thm: lower bound for integral}

\subsection{Proof of Theorem \ref{thm: upper bound for integral} (Upper Bound on Integral Estimation)}
\label{app: proof of upper bound in thm 4.2}
Before proving the upper bound on integral estimation, we need to derive an upper bound on the expected error of the $k$-nearest neighbor estimator $\hat{f}_{k\text{-NN}}$, which is built based on the first half of the given dataset $\{(x_i,y_i)\}_{i=1}^{n}$, with respect to the $L^2$ norm. From our construction of $\hat{f}_{k\text{-NN}}$ given in Section \ref{sec 4.2: k nearest neighbor is optimal for regression adjusted CV}, we have that for any fixed $\frac{n}{2}$ quadrature points $\{x_i\}_{i=1}^{\frac{n}{2}}$, $z \in \Omega$ and $k \in \{1,2,\cdots,\frac{n}{2}\}$, the expected value of $\hat{f}_{k\text{-NN}}(z)$ with respect to the observational noises $\{\epsilon_i\}_{i=1}^{\frac{n}{2}}$ is given by
\begin{equation}
\label{eqn: thm 4.2 expected value of k-NN}
\mathlarger{\mathbb{E}}_{\{\epsilon_i\}_{i=1}^{\frac{n}{2}}}\Big[\hat{f}_{k\text{-NN}}(z)\Big] = \frac{1}{k}\sum_{j=1}^{k}\mathlarger{\mathbb{E}}_{\{\epsilon_i\}_{i=1}^{\frac{n}{2}}}\Big[f(x_{i^{(z)}_j})+\epsilon_{i^{(z)}_j}\Big] = \frac{1}{k}\sum_{j=1}^{k}f(x_{i^{(z)}_j}), 
\end{equation}
where $\{x_{i^{(z)}_{j}}\}_{j=1}^{k}$ above are the $k$ nearest neighbors of $z$ among $\{x_{i}\}_{i=1}^{\frac{n}{2}}$. Now let's consider using the bias-variance decomposition to upper bound the error $||\hat{f}_{k\text{-NN}}(z) -f(z)||_{L^2(\Omega)}^2$. Based on the expected value computed in (\ref{eqn: thm 4.2 expected value of k-NN}) above, we may decompose the function $\hat{f}_{k\text{-NN}} -f$ as a sum of the bias part and the variance part as follows:
\begin{equation}
\label{eqn: thm 4.2 bias part}
B(z) := \mathlarger{\mathbb{E}}_{\{\epsilon_i\}_{i=1}^{\frac{n}{2}}}\Big[\hat{f}_{k\text{-NN}}(z)\Big] - f(z) =  \frac{1}{k}\sum_{j=1}^{k}f(x_{i^{(z)}_j})-f(z) = \frac{1}{k}\sum_{j=1}^{k}\Big(f(x_{i^{(z)}_j})-f(z)\Big),
\end{equation}
\begin{equation}
\label{eqn: thm 4.2 variance part}
V(z) := \hat{f}_{k\text{-NN}}(z) - \mathlarger{\mathbb{E}}_{\{\epsilon_i\}_{i=1}^{\frac{n}{2}}}\Big[\hat{f}_{k\text{-NN}}(z)\Big] = \hat{f}_{k\text{-NN}}(z)-\frac{1}{k}\sum_{j=1}^{k}f(x_{i^{(z)}_j})  = \frac{1}{k}\sum_{j=1}^{k}\epsilon_{i_{j}^{(z)}},
\end{equation}
where the function $B$ corresponds to the bias part and the function $V$ corresponds to the variance part. Using the decomposition $\hat{f}_{k\text{-NN}} -f=B+V$ allows us to upper bound the expected error of $\hat{f}_{k\text{-NN}}$ with respect to the $L^2$ norm as below:
\begin{equation}
\label{eqn: thm 4.2 bias variance decomp}
\begin{aligned}
&\mathlarger{\mathbb{E}}_{\substack{\{x_i\}_{i=1}^{\frac{n}{2}},\{y_i\}_{i=1}^{\frac{n}{2}}}}\Big[||\hat{f}_{k\text{-NN}} - f||_{L^2(\Omega)}^2\Big] = \mathlarger{\mathbb{E}}_{\{x_i\}_{i=1}^{\frac{n}{2}},\{y_i\}_{i=1}^{\frac{n}{2}}}\Big[||B+V||_{L^2(\Omega)}^2\Big]\\
&\leq \mathlarger{\mathbb{E}}_{\{x_i\}_{i=1}^{\frac{n}{2}},\{y_i\}_{i=1}^{\frac{n}{2}}}\Bigg[\Big(||B||_{L^2(\Omega)}+||V||_{L^2(\Omega)}\Big)^2\Bigg] \\
&\leq \mathlarger{\mathbb{E}}_{\{x_i\}_{i=1}^{\frac{n}{2}},\{y_i\}_{i=1}^{\frac{n}{2}}}\Big[2||B||_{L^2(\Omega)}^2+2||V||_{L^2(\Omega)}^2\Big]\\
&\lesssim \mathlarger{\mathbb{E}}_{\{x_i\}_{i=1}^{\frac{n}{2}},\{y_i\}_{i=1}^{\frac{n}{2}}}\Big[||V||_{L^2(\Omega)}^2\Big] + \mathlarger{\mathbb{E}}_{\{x_i\}_{i=1}^{\frac{n}{2}},\{y_i\}_{i=1}^{\frac{n}{2}}}\Big[||B||_{L^2(\Omega)}^2\Big]\\
&= \mathlarger{\mathbb{E}}_{z,\{x_i\}_{i=1}^{\frac{n}{2}},\{\epsilon_i\}_{i=1}^{\frac{n}{2}}}\Big[|V(z)|^2\Big] + \mathlarger{\mathbb{E}}_{z,\{x_i\}_{i=1}^{\frac{n}{2}},\{\epsilon_i\}_{i=1}^{\frac{n}{2}}}\Big[|B(z)|^2\Big],
\end{aligned}    
\end{equation}
where $z$ above is uniformly distributed over the domain $\Omega$ and independent of $x_i$ for any $1 \leq i \leq \frac{n}{2}$. On the one hand, using the expression of the variance part $V$ derived in (\ref{eqn: thm 4.2 variance part}) above and the fact that $\{\epsilon_i\}_{i=1}^{\frac{n}{2}}$ are independent and identical distributed noises, we may compute the first term in (\ref{eqn: thm 4.2 bias variance decomp}) above as follows:
\begin{equation}
\label{eqn: thm 4.2 bound on variance}
\begin{aligned}
\mathlarger{\mathbb{E}}_{z,\{x_i\}_{i=1}^{\frac{n}{2}},\{\epsilon_i\}_{i=1}^{\frac{n}{2}}}\Big[|V(z)|^2\Big] &= \mathlarger{\mathbb{E}}_{z}\Bigg[\mathlarger{\mathbb{E}}_{\{x_i\}_{i=1}^{\frac{n}{2}},\{\epsilon_i\}_{i=1}^{\frac{n}{2}}}\Bigg[\Big|\frac{1}{k}\sum_{j=1}^{k}\epsilon_{i_{j}^{(z)}}\Big|^2\Bigg]\Bigg] \\
&= \mathlarger{\mathbb{E}}_{z}\Bigg[\frac{1}{k^2}\mathlarger{\mathbb{E}}_{\{x_i\}_{i=1}^{\frac{n}{2}},\{\epsilon_i\}_{i=1}^{\frac{n}{2}}}\Bigg[\sum_{j=1}^{k}\epsilon_{i_{j}^{(z)}}^2\Bigg]\Bigg] \\
&= \mathbb{E}_{z}\Big[\frac{n^{-2\gamma}k}{k^2}\Big] = \frac{n^{-2\gamma}}{k}.
\end{aligned}   
\end{equation}

On the other hand, since $s \in (0,1)$ and the given function $f$ is $s$-H\"older smooth, we have that the inequality $|f(x)-f(y)| \lesssim ||x-y||^s$ holds true for any $x,y \in \Omega$. Combining this inequality with the expression of the bias part $B$ derived in (\ref{eqn: thm 4.2 variance part}) above helps us upper bound the second term in (\ref{eqn: thm 4.2 bias variance decomp}) as below:
\begin{equation}
\label{eqn: thm 4.2 bound on bias}
\begin{aligned}
\mathlarger{\mathbb{E}}_{z,\{x_i\}_{i=1}^{\frac{n}{2}},\{\epsilon_i\}_{i=1}^{\frac{n}{2}}}\Big[|B(z)|^2\Big] &= \mathlarger{\mathbb{E}}_{z}\Bigg[\mathlarger{\mathbb{E}}_{\{x_i\}_{i=1}^{\frac{n}{2}},\{\epsilon_i\}_{i=1}^{\frac{n}{2}}}\Bigg[\Big|\frac{1}{k}\sum_{j=1}^{k}\Big(f(x_{i^{(z)}_j})-f(z)\Big)\Big|^2\Bigg]\Bigg] \\
&\leq \frac{1}{k}\mathlarger{\mathbb{E}}_{z}\Bigg[\mathlarger{\mathbb{E}}_{\{x_i\}_{i=1}^{\frac{n}{2}}}\Bigg[\sum_{j=1}^{k}\Big|f(x_{i^{(z)}_j})-f(z)\Big|^2\Bigg]\Bigg]\\
&\lesssim \frac{1}{k}\mathlarger{\mathbb{E}}_{z}\Bigg[\mathlarger{\mathbb{E}}_{\{x_i\}_{i=1}^{\frac{n}{2}}}\Bigg[\sum_{j=1}^{k}\Big|x_{i^{(z)}_j}-z\Big|^{2s}\Bigg]\Bigg]\\
&\leq \mathlarger{\mathbb{E}}_{z,\{x_i\}_{i=1}^{\frac{n}{2}}}\Bigg[\Big|x_{i^{(z)}_k}-z\Big|^{2s}\Bigg]\\
&\leq \Bigg(\mathlarger{\mathbb{E}}_{z,\{x_i\}_{i=1}^{\frac{n}{2}}}\Bigg[\Big|x_{i^{(z)}_k}-z\Big|^{2}\Bigg]\Bigg)^s \lesssim \Big(\frac{k}{n}\Big)^\frac{2s}{d}.
\end{aligned}    
\end{equation}
The second least inequality follows from the fact that $\omega(t):=t^s$ is a concave function when $s\in (0,1)$, while the last inequality is obtained by plugging in (\ref{eq: upper bound on expected distance of the k-th distance}) given in Lemma \ref{lem: upper bound on expected distance of the k-th distance}. Substituting (\ref{eqn: thm 4.2 bound on variance}) and (\ref{eqn: thm 4.2 bound on bias}) into (\ref{eqn: thm 4.2 bias variance decomp}) then yields that for any $k \in \{1,2,\cdots,\frac{n}{2}\}$, the expected error of $\hat{f}_{k\text{-NN}}$ with respect to the $L^2$ norm can be upper bounded as follows:
\begin{equation}
\label{eqn: thm 4.2 bound on expected L2 norm}
\begin{aligned}
\mathlarger{\mathbb{E}}_{\substack{\{x_i\}_{i=1}^{\frac{n}{2}},\{y_i\}_{i=1}^{\frac{n}{2}}}}\Big[||\hat{f}_{k\text{-NN}} - f||_{L^2(\Omega)}^2\Big] \lesssim  \frac{n^{-2\gamma}}{k} + \Big(\frac{k}{n}\Big)^\frac{2s}{d}. 
\end{aligned}    
\end{equation}
Furthermore, from our construction of the integral estimator $\hat{H}_{k\text{-NN}}$ given in Section \ref{sec 4.2: k nearest neighbor is optimal for regression adjusted CV}, we may upper bound the expectation of the estimator $\hat{H}_{k\text{-NN}}$'s squared error via the expected error of $\hat{f}_{k\text{-NN}}$ with respect to the $L^2$ norm as below:
\begin{equation}
\label{eqn: thm 4.2 second last step of final upper}
\begin{aligned}
&\mathlarger{\mathlarger{\mathbb{E}}}_{\substack{\{x_i\}_{i=1}^{n},\{y_i\}_{i=1}^{n}}}\Bigg[\left|\hat{H}_{k\text{-NN}}\Big(\{x_i\}_{i=1}^{n},\{y_i\}_{i=1}^{n}\Big) - I_f\right|^2\Bigg] \\
&= \mathlarger{\mathlarger{\mathbb{E}}}_{\substack{\{x_i\}_{i=1}^{n},\{y_i\}_{i=1}^{n}}}\Bigg[\left|\int_{\Omega}\hat{f}_{k\text{-NN}}(x)dx + \frac{2}{n}\sum_{i=\frac{n}{2}+1}^{n}\Big(y_i-\hat{f}_{k\text{-NN}}(x_i)\Big)- \int_{\Omega}f(x)dx\right|^2\Bigg]\\
&\lesssim  \mathlarger{\mathlarger{\mathbb{E}}}_{\substack{\{x_i\}_{i=1}^{\frac{n}{2}}, \\ \{y_i\}_{i=1}^{\frac{n}{2}}}}\Bigg[\mathlarger{\mathlarger{\mathbb{E}}}_{\substack{\{x_i\}_{i=\frac{n}{2}+1}^{n}, \\ \\ \{y_i\}_{i=\frac{n}{2}+1}^{n}}}\Bigg[\Big|\frac{1}{\frac{n}{2}}\sum_{i=\frac{n}{2}+1}^{n}\Big(f(x_i)-\hat{f}_{k\text{-NN}}(x_i)-\int_{\Omega}(f(x)-\hat{f}_{k\text{-NN}}(x))dx\Big)\Big|^2\Bigg]\Bigg]\\
&+ \mathlarger{\mathlarger{\mathbb{E}}}_{\substack{\{x_i\}_{i=1}^{n},\{y_i\}_{i=1}^{n}}}\Bigg[\Big|\frac{2}{n}\sum_{i=\frac{n}{2}+1}^{n}\epsilon_i\Big|^2\Bigg]\\
&= \mathlarger{\mathlarger{\mathbb{E}}}_{\substack{\{x_i\}_{i=1}^{\frac{n}{2}}, \{y_i\}_{i=1}^{\frac{n}{2}}}}\Bigg[\frac{4}{n^2}\sum_{i=\frac{n}{2}+1}^{n}\mathlarger{\mathlarger{\mathbb{E}}}_{x_i}\Bigg[\Big|\Big(f(x_i)-\hat{f}_{k\text{-NN}}(x_i)-\int_{\Omega}(f(x)-\hat{f}_{k\text{-NN}}(x))dx\Big)\Big|^2\Bigg]\Bigg]\\
&+ \frac{4}{n^2}\sum_{i=\frac{n}{2}+1}^{n}\mathlarger{\mathlarger{\mathbb{E}}}_{x_i,y_i}\Big[\epsilon_i^2\Big] \lesssim \frac{1}{n}\Bigg(\mathlarger{\mathbb{E}}_{\substack{\{x_i\}_{i=1}^{\frac{n}{2}},\{y_i\}_{i=1}^{\frac{n}{2}}}}\Big[||\hat{f}_{k\text{-NN}} - f||_{L^2(\Omega)}^2\Big] + n^{-2\gamma}\Bigg)\\
&\lesssim \frac{1}{n}\Big(\frac{n^{-2\gamma}}{k} + \Big(\frac{k}{n}\Big)^\frac{2s}{d}\Big) + n^{-2\gamma -1}.
\end{aligned}    
\end{equation}
Based on the magnitude of the noises, we have the following two cases for the final upper bound:\\
When $\gamma \in [0,\frac{s}{d})$, the optimal $k$ is determined by balancing the two terms $\frac{n^{-2\gamma}}{k}$ and $\Big(\frac{k}{n}\Big)^{\frac{2s}{d}}$ in (\ref{eqn: thm 4.2 second last step of final upper}), which yields $\frac{n^{-2\gamma}}{k} = \Big(\frac{k}{n}\Big)^{\frac{2s}{d}} \Rightarrow k=\Theta(n^{\frac{2(s-\gamma d)}{d+2s}})$. The corresponding upper bound is given by
\begin{equation}
\label{eqn: thm 4.2 high noise case bound on squared expected error}
\begin{aligned}
\frac{1}{n}\Big(\frac{n^{-2\gamma}}{k} + \Big(\frac{k}{n}\Big)^\frac{2s}{d}\Big) + n^{-2\gamma -1} &\lesssim \frac{1}{n}n^{-2\gamma - \frac{2(s-\gamma d)}{d+2s}} + n^{-1-2\gamma} = n^{-\frac{2s(1+2\gamma)}{2s+d}-1} + n^{-2\gamma - 1}\\
&\lesssim \max\{n^{-\frac{2s(1+2\gamma)}{2s+d}-1}, n^{-2\gamma - 1}\} = n^{-2\gamma - 1}. 
\end{aligned}
\end{equation}
When $\gamma \in [\frac{s}{d},\infty]$, we note that $k \in \{1,2,\cdots,\frac{n}{2}\}$ must be of at least constant level. Therefore, the optimal $k$ is determined by balancing the two terms $\frac{n^{-2\gamma -1}}{k}$ and $n^{-2\gamma -1}$, which yields that $k=\Theta(1)$ is of constant level. The corresponding upper bound is given by
\begin{equation}
\label{eqn: thm 4.2 low noise case bound on squared expected error}
\begin{aligned}
\frac{1}{n}\Big(\frac{n^{-2\gamma}}{k} + \Big(\frac{k}{n}\Big)^\frac{2s}{d}\Big) + n^{-2\gamma -1} &\lesssim n^{-\frac{2s}{d}-1} +n^{-2\gamma - 1} \\
&\lesssim \max\{n^{-\frac{2s}{d}-1}, n^{-2\gamma - 1}\} = n^{-\frac{2s}{d}-1}.
\end{aligned}
\end{equation}
Finally, substituting (\ref{eqn: thm 4.2 high noise case bound on squared expected error}) and  (\ref{eqn: thm 4.2 low noise case bound on squared expected error}) into (\ref{eqn: thm 4.2 second last step of final upper}) gives us the final upper bound:
\begin{equation}
\begin{aligned}
&\mathlarger{\mathlarger{\mathbb{E}}}_{\substack{\{x_i\}_{i=1}^{n},\{y_i\}_{i=1}^{n}}}\Bigg[\left|\hat{H}_{k\text{-NN}}\Big(\{x_i\}_{i=1}^{n},\{y_i\}_{i=1}^{n}\Big) - I_f\right|\Bigg]\\
&\leq \sqrt{\mathlarger{\mathlarger{\mathbb{E}}}_{\substack{\{x_i\}_{i=1}^{n},\{y_i\}_{i=1}^{n}}}\Bigg[\left|\hat{H}_{k\text{-NN}}\Big(\{x_i\}_{i=1}^{n},\{y_i\}_{i=1}^{n}\Big) - I_f\right|^2\Bigg]}  \\
&\lesssim \sqrt{\frac{1}{n}\Big(\frac{n^{-2\gamma}}{k} + \Big(\frac{k}{n}\Big)^\frac{2s}{d}\Big) + n^{-2\gamma -1}} \lesssim \sqrt{\max\{n^{-\frac{2s}{d}-1}, n^{-2\gamma - 1}\}}\\
&= n^{\max\{-\frac{1}{2}-\gamma, -\frac{1}{2}-\frac{s}{d}\}},  
\end{aligned}    
\end{equation}
which concludes our proof of Theorem \ref{thm: upper bound for integral}.

\end{document}